

\errorcontextlines999
\newif\ifllncs
\llncstrue
\llncsfalse
\ifllncs
\documentclass{llncs}
\usepackage{cc4llncs}
\makeatletter
\@longarticletrue
\AtBeginDocument{\@adjustnumbering}
\makeatother
\usepackage[T1]{fontenc}
\usepackage[latin1]{inputenc}
\usepackage[german,polutonikogreek,english]{babel}
\else
\documentclass[noccpublish,balancedmargins,
11pt,textarea=cc1,
multilingual,german,polutonikogreek,english]{cc}
\fi



\newif\ifdraft%
\drafttrue%
\draftfalse%

\newif\ifshortversion%
\shortversiontrue%
\shortversionfalse%

\newif\iffullversion%
\fullversionfalse%
\fullversiontrue%

\newif\ifwatermark%
\watermarktrue%
\watermarkfalse%

\newif\ifhyper
\hyperfalse
\hypertrue

\newif\ifanonymous%
\anonymoustrue%
\anonymousfalse%

\newif\ifbw
\bwtrue%
\bwfalse%

\ifdraft

\usepackage[draft,showlabels]{develop}
\ifhyper\makeatletter\AtBeginDocument{\def\oops@label{\label}}\makeatother\fi
\def\rightflag#1{}

\def\latertodo{\todo}

\else

\def\latertodo#1{}

\fi

\usepackage{xcolor}
\ifwatermark
\usepackage{watermark}
\definecolor{watermarkcolor}{cmyk}{.08,.052,0,0}
\watermark{\rput[tl]{60}(0,-\textheight){%
      \psscaleboxto(22cm,15mm){\Huge\bf\color{watermarkcolor}
    DRAFT
    ---
    Do not distribute!
  }}}

\fi

\makeatletter
\newenvironment{shortversion}[1][]{\ignorespaces}{\ignorespaces}
\newenvironment{fullversion}[1][]{\ignorespaces}{\ignorespaces}
\usepackage{skipenv}\newenvironment{skipme}{}{}\skipenv{skipme}
\ifdraft
\renewenvironment{shortversion}[1][]{%
  \@bsphack\small\blue\@esphack
}{}
\definecolor{darkgreen}{hsb}{0.333333,1,0.5}{}
\renewenvironment{fullversion}[1][]{%
  \@bsphack\color{darkgreen}\small\@esphack
}{}
\fi
\iffullversion\else
\renewenvironment{fullversion}[1][]{%
  \@bsphack\shortversion\@esphack
  #1\endshortversion\skipme}{\endskipme}
\fi
\ifshortversion\else
\iffullversion
\renewenvironment{shortversion}[1][]{%
  \@bsphack\fullversion\@esphack
  #1\endfullversion\skipme}{\endskipme}
\else
\renewenvironment{shortversion}[1][]{\@bsphack\skipme}{\endskipme\@esphack}%
\fi\fi
\makeatother

\let\Cite\Citet
\ifdraft
\else
\fi


\ifhyper
\usepackage[breaklinks=true,dvips]{hyperref}
\usepackage{breakurl}
\fi

\usepackage{multicol}
\usepackage{supertabular}
\usepackage{graphicx}
\usepackage{pstricks}[2008/07/31]
\usepackage{pst-node}[2006/01/01]
\usepackage{pst-plot}[2005/05/20]
\usepackage{pst-xkey}[2005/11/25]
\usepackage{pst-all}

\usepackage{tilde}
\usepackage{maps}
\usepackage{mathell}
\usepackage{set}
\usepackage{spacedno}
\usepackage{rounded}
\usepackage{bigOetc}

\def\indicate#1{\mathchoice
  {\left[\!\!\middle[#1\middle]\!\!\right]}
  {\left[\!\middle[#1\middle]\!\right]}
  {\left[\!\middle[#1\middle]\!\right]}
  {\left[\!\middle[#1\middle]\!\right]}}%

\ifdraft
\makeatletter\def\thefootnote{\@fnsymbol\c@footnote}\makeatother
\usepackage{endnotes}
\long\def\FOOTNOTE#1{ {\footnotesize$\lceil$#1$\rfloor$}}
\else
\long\def\endnote#1{}%
\long\def\endnotetext#1{}%
\fi

\listfiles

\def\glnB{4.827869664 }%
\def\lnormoneatone{0.2054729137 }

\def\swap#1#2{\let\swaptmp#1\let#1#2\let#2\swaptmp}
\swap\epsilon\varepsilon
\swap\theta\vartheta
\swap\rho\varrho

\newcommand{\defiff}{\;:\Longleftrightarrow\;}

\newcommand{\N}{\mathbb{N}}

\newcommand{\R}{\mathbb{R}}

\newcommand{\floor}[1]{\left\lfloor#1\right\rfloor}


\makeatletter
\def\treatparameter{%
  \@ifnextchar({\treatparameterx}{%
    \@ifnextchar<{\treatparameterxi}{}}}
\makeatother
\def\treatparameterx(#1){\left\lgroup#1\right\rgroup}
\def\treatparameterx(#1){\left(#1\right)}
\def\treatparameterxi<#1>{\left\ldbrack#1\right\rdbrack}
\def\treatparameterxi<#1>{\left\langle#1\right\rangle}

\def\braket#1#2{\Left\langle #1;#2\Right\rangle}

\newcommand{\Qclean}[2][B,C]{\pi_{#1}^{#2}\treatparameter}
\newcommand{\Qcleanapprox}[2][B,C]{\widetilde{\pi}_{#1}^{#2}\treatparameter}
\newcommand{\Qcleanerror}[2][B,C]{\widehat{\pi}_{#1}^{#2}\treatparameter}
\newcommand{\Q}[2][B,C]{\kappa_{#1}^{#2}\treatparameter}
\newcommand{\Qapprox}[2][B,C]{\widetilde{\kappa}_{#1}^{#2}\treatparameter}
\newcommand{\Qerror}[2][B,C]{\widehat{\kappa}_{#1}^{#2}\treatparameter}

\newcommand{\Qboth}[2][B,C]{(\widetilde{\kappa}_{#1}^{#2}+\widehat{\kappa}_{#1}^{#2})\treatparameter}

\newcommand{\Kapprox}[2][]{\widetilde{\lambda}_{#1}^{#2}\treatparameter}
\newcommand{\Kerror}[2][]{\widehat{\lambda}_{#1}^{#2}\treatparameter}

\newcommand{\Kboth}[2][]{(\widetilde{\lambda}_{#1}^{#2}+\widehat{\lambda}_{#1}^{#2})\treatparameter}

\newcommand{\Napprox}[2][]{\widetilde{\nu}_{#1}^{#2}\treatparameter}
\newcommand{\Eapprox}[2][]{\widetilde{\eta}_{#1}^{#2}\treatparameter}

\ifcase1
\newcommand{\capprox}[1]{c_{\Kapprox{#1}}}
\newcommand{\cerror}[1]{c_{\Kerror{#1}}}
\newcommand{\clower}[1]{\delta_{\Kapprox{#1}}}
\or
\newcommand{\capprox}[1]{\widetilde{c}_{#1}}
\newcommand{\cerror}[1]{\widehat{c}_{#1}}
\newcommand{\clower}[1]{\widebreve{c}_{#1}}
\or
\newcommand{\capprox}[1]{[\newalpha^{-(#1)} c_{\Kapprox{#1}}]}
\newcommand{\cerror}[1]{c_{\Kerror{#1}}}
\newcommand{\clower}[1]{\delta_{\Kapprox{#1}}}
\else
\newcommand{\capprox}[1]{\widehat{\Kapprox[~norm]{#1}}}
\newcommand{\cerror}[1]{\widehat{\Kerror[~norm]{#1}}}
\fi
\newcommand{\cboth}[1]{\left(\capprox{#1}+\cerror{#1}\right)}

\def\myexp#1{~exp\treatparameter({#1})}
\def\myexp#1{~e^{#1}}

\usepackage{sumint}



\renewcommand{\P}[1]{{\mathbb{P}\mathbin{\cap}\left]#1\right]}}
\def\primesum#1#2#3#4{\sum_{#2 \in \P{#3,#4}} #1}%

\renewcommand{\limsup}{\mathop{~limsup}\limits}

\def\Frac#1#2{{\mathpalette{\let\cursize}{%
      \raise0.5ex\hbox{$\cursize#1$}/\lower0.5ex\hbox{$\cursize#2$}}}}%
\def\transpose#1{#1^{\mathsf{T}}}%

\iftrue
\fi\iftrue
\newcommand{\weaklyisotone}{increasing}
\newcommand{\weaklyantitone}{decreasing}
\newcommand{\strictlyisotone}{strictly increasing}

\else
\newcommand{\weaklyisotone}{isotone}
\newcommand{\weaklyantitone}{antitone}
\newcommand{\strictlyisotone}{strictly isotone}

\fi

\newcommand{\intrangeleft}{B}
\newcommand{\intrangeright}{C}

\def\eulernumber{\textup{e}}%
\def\Tau{T}%
\def\tausplit{\tau_{\text{split}}}

\makeatletter
\newtoks\@smashedheights \@smashedheights{}%
\def\smashedunderbrace#1_#2{%
  \vphantom{#1}\smash{\underbrace{#1}_{#2}}%
  \expandafter\global\expandafter\@smashedheights\expandafter{\the\@smashedheights
    \vphantom{\underbrace{#1}_{#2}}}%
}
\def\smashedheights{\the\@smashedheights\global\@smashedheights{}}
\makeatother

\makeatletter
\def\casesincreaseheight#1{%
  \expandafter\global\expandafter\setbox\csname
  @arstrutbox\endcsname\hbox{%
    \vrule height#1\ht\strutbox depth#1\dp\strutbox width0pt}}
\makeatother
\def\allowpagebreak{\noalign{\penalty\postdisplaypenalty}}

\let\newalpha\alpha

\ifdraft
\AtEndDocument{%
  \clearpage\appendix
  \section{List of symbols}
  \label{sec:labels}
  \begin{multicols}{2}\tiny
    \begin{tabular}{@{}lp{0.7\hsize}}
      $a$ & Coefficient of monomial, old attempt \\
      $b$ & Constant, OLD \\
      $c$ & Constant \\
      $d$ & Differential operator \\
      $e$ & Euler's constant \\
      $f$ & Arbitrary function, monomial $\xi^nB^{k+\xi\newalpha}$ \\
      $g$ & Arbitrary function after one recursive approximation \\
      $h$ & Constant part in the solution \\
      $i$ & Index \\
      $j$ & Twig \\
      $k$ & Branch \\
      $l$ & Non-constant part in the solution; index OLD\\
      $m$ & Hill polynomials \\
      $n$ & Integer, degree of monomial OLD \\
      $o$ & Landau symbol\\
      $p$ & Primes \\
      $q$ & \\
      $r$ & Index NOTES \\
      $s$ & Bound s.t. $C < B^s$; index OLD \\
      $t$ & Index OLD \\
      $u$ & Index OLD \\
      $v$ & Index\\
      $w$ & Index \\
      $x$ & Bound \\
      $y$ & Bound CRAMER \\
      $z$ & Variable, OLD \\
      $\alpha$ & General variable $B^{1+\alpha\newalpha}$; function CRAMER \\
      $\beta$ & General variable $B^{1+\beta\newalpha}$ \\
      $\gamma$ & Real number, s.t. $C=B^{1+\newalpha}$ \\
      $\delta$ & Constant in the approximation of the error\\
      $\epsilon$ & \\
      $\zeta$ & \\
      $\eta$ & \\
      $\theta$ & \\
      $\iota$ & \\
      $\kappa$ & Cleaned aim of study \\
      $\lambda$ & Approximation of $\kappa$ \\
      $\mu$ & \\
      $\nu$ & \\
      $\xi$ & Log-scale $x$, i.e. $x = B^{k+\xi\newalpha}$ \\
      $\pi$ & Aim of study; prime counting function; the number $\pi$ \\
      $\rho$ & Log-scale $p$, i.e. $p = B^{k+\rho\newalpha}$ \\
      $\sigma$ & Permuatation, NOTES; standard deviation CRAMER \\
      $\tau$ & Function $z \mapsto \frac{(z^n)^{(i)}}{(-\newalpha/2)^i}$ \\
      $\upsilon$ & \\
      $\phi$ & \\
      $\chi$ & \\
      $\psi$ & \\
      $\omega$ & \\
    \end{tabular}

    \begin{tabular}{@{}lp{0.7\hsize}}
      $A$ & Constant in PNT; Partition set NOTES\\
      $B$ & Lower bound\\
      $C$ & Upper bound; Decker-Moree RSA-counting function; Complement CRAMER\\
      $D$ & \\
      $E$ & Error term, Event CRAMER\\
      $F$ & \\
      $G$ & \\
      $H$ & \\
      $I$ & \\
      $J$ & \\
      $K$ & Random variable CRAMER \\
      $L$ & \\
      $M$ & \\
      $N$ & \\
      $O$ & \\
      $P$ & Tuple NOTES \\
      $Q$ & \\
      $R$ & \\
      $S$ & $B^{k-1-j}C^j$; sorting NOTES\\
      $T$ & Type NOTES\\
      $U$ & \\
      $V$ & \\
      $W$ & \\
      $X$ & Random variable CRAMER\\
      $Y$ & \\
      $Z$ & Random variable CRAMER\\
      $\Gamma$ & \\
      $\Delta$ & \\
      $\Theta$ & Landau symbol \\
      $\Lambda$ & \\
      $\Xi$  & \\
      $\Pi$ & Cram\'{e}r prime counting function \\
      $\Sigma$ & \\
      $\Upsilon$ & \\
      $\Phi$ & \\
      $\Psi$ & \\
      $\Omega$ & \\
      $\mathcal{D}_{\xi}$ & Differential operator \\
      $\mathcal{I}$ & Integral operator \\
      $\mathcal{L}$ & $\frac{-1}{\newalpha ~ln B}\frac{~d}{~d \xi}$ \\
      $\mathcal{M}$ & Integral operator \\
      $\mathcal{O}$ & Landau symbol \\
      $\mathcal{S}$ & Shift operator \\
      $\mathbb{N}$ & Set of nonnegative integers \\
      $\mathbb{P}$ & Set of primes \\
      $\mathbb{R}$ & Set of reals \\
      $\mathbb{Z}$ & Set of integers \\
    \end{tabular}
  \end{multicols}
  \clearpage
  \section{TODO}
  \label{sec:todo}
  
  \begin{itemize}\def\done{\item[\checkmark]}

    \done Final aim: approximate $\Q{k}(x)$ by a piecewise analytic
    function and estimate the error term.  The error term must be at
    least by a factor of order $\Frac{1}{~ln x}$ smaller than the main
    term. REACHED!

    \done Prove $\Kapprox{k}(x) \in \Omega \left( \frac{x}{~ln B}
    \right)$ in a nice form.

    \done Approximate $\Q{2}(x)$ and $\Qerror{2}(x)$ in detail.
    (We have to redo all the work, since we used a wrong form of
    \ref{sec:using-estimates}\dots)
    
    In \texttt{BroughCsmoothCaser5-2D.mn} we need to check the order of
    individual terms under (NEW!) an extra condition like $B<C<B^{2}$.
    So we have essentially one asymptotic direction.
    
    Automate the worksheet using ideas from \ref{sec:firstattempt}.
    
    Continue with case $x \in \left]BC,C^{2}\right]$.

    Estimate the asymptotic order of integrals like $\int_{B}^{x/B}
    \int_{B}^{x/p_{2}} \frac{~d p_{1}}{~ln p_{1}} \frac{~d p_{2}}{~ln
      p_{2}}$ for $B^{2} \leq x < B C$ rewriting $C=B^{1+\newalpha}$, $x=B^{2}
    B^{\xi}$ under a condition like $0\leq \xi \leq \newalpha \leq 1$ (ie.\ $B < C
    < B^{2}$).  The condition implies that $~ln B < ~ln(p_{2}) < ~ln C
    = (1+\newalpha) ~ln B$ and so we can replace logarithmic terms like
    $~ln(p_{2})$ by a constant losing no more than a possible factor
    $1+\newalpha$.
    
    Example:
    \begin{align*}
      \Qapprox{2}(x)
      &=
      \int_{B}^{x/B} \int_{B}^{x/p_{2}} \frac{~d p_{1}}{~ln p_{1}}
      \frac{~d p_{2}}{~ln p_{2}}
      \\
      &=
      \frac{1}{c_{1}c_{2}}
      \int_{B}^{x/B} \int_{B}^{x/p_{2}} \frac{~d p_{1}}{~ln B}
      \frac{~d p_{2}}{~ln B}
      \\
      &=
      \frac{1}{c^{2} ~ln^{2} B}
      \frint{\frint{1}{p_{2}}{B}{x/p_{2}}}{p_{1}}{B}{x/B}
      \\
    \end{align*}
    for some $c_{1},c_{2}, c\in\left[1,1+\newalpha\right]$.
    
    \done Describe introduction of constants:
    $c_{i} \cdot \frint{\frac{f_{i}(p)}{~ln p}}{p}{a_{i}}{b_{i}} = 
    \frint{\frac{f_{i}(p)}{~ln B}}{p}{a_{i}}{b_{i}}$\dots{}
    Repair/improve the implementation so that $c_{i}$ depending on an
    outer variable (like $p_{3}$) is treated correctly, ie.\
    \begin{align*}
      c_{i} &\cdot \frint{g_{i}(p_{2})\frint{\frac{f_{i}(p_{1},p_{2})}{~ln
            p_{1}}}{p_{1}}{a_{i}(p_{2})}{b_{i}(p_{2})}}{p_{2}}{a_{2}}{b_{2}}
      \\&= 
      \frint{g_{i}(p_{2})\frint{\frac{f_{i}(p_{1},p_{2})}{~ln
            B}}{p_{1}}{a_{i}(p_{2})}{b_{i}(p_{2})}}{p_{2}}{a_{2}}{b_{2}}
    \end{align*}
    \dots{}
    
    \done Prove $k=2$ case, maybe slightly weakened, without
    endless calculations.

    \done Complete \ref{sec:solve-Kerror}. DECIDED not to do this.

  \item Understand \citeauthor{cra35}.  Apply the \citeauthor{cra35}
    methods to~$\Q{2}(x)$ to understand the expected size of the
    error.  Start with determining the order with respect to~$B$
    (under condition that $C=B^{1+\newalpha}$, $0\leq \newalpha \leq 1$).

    \done Continue with higher cases.

    \done Look up \cite{sch76a} to see whether there is a conclusion
    like: if the Riemann hypothesis is correct for $|t| < B$ then the
    bound \ref{pnt-schoenfeld} holds for $x < f(B)$\dots{} Numerical
    calculations have shown that the Riemann hypthesis does hold for
    the first $2 \cdot 10^{11}$ critical zeros.
    
    This is asymptotically(!) impossible: if there is a single zero
    with real part $t$ strictly larger than $\frac{1}{2}$ then the
    error bound is at best $\mathcal{O}( x^{t} )$.
    
  \item Check previous again!  Write down Riemann formula for $\pi$
    somewhere.

    \done Rewrite $\frint{f(\rho)}{\rho}{\xi-1}{\xi}$ to
    $\frint{f(\xi-\rho)}{\rho}{0}{1}$ to unify the presentation.

    \done Create a symbol list.

    \done NFS application!?

  \item Verify \ref{pnt-schoenfeld} numerically.
    \begin{itemize}
    \item Use the compute server.
    \item Compute $\pi(x_{i})$ and determine for which interval this
      guarantees the inequality.
    \item Parallelize: each cpu gets an interval to cover.
    \item We can use forward and backward estimates, so we need good
      guesses to obtain a cover with only small overlaps.
    \end{itemize}
    Check
    \url{http://numbers.computation.free.fr/Constants/Primes/Pix/results.html}
    for possible programs\dots

    23 January 2009: Check done up to $10^{9}$ based on available
    values for $\pi$.  For larger arguments we need a denser set of
    values beyond $10^{9}$.  (We have values for $10^{9}+k\cdot
    10^{6}$ with $k\in\N_{<9001}$.  We would need something like
    $10^{9}+0.35\cdot10^{6}$\dots)

    4 February 2009: Check done up to $10^{11}$.  Started to build
    denser $\pi$ table to continue up to $10^{12}$.  Plan to use
    sieving in C to determine maxima and minima of
    $\frac{\pi(x)-~Li(x)}{\frac{1}{8\pi} \sqrt{x} ~ln(x)}$ for $x\leq
    10^{15}$?

  \item Choose (for \ref{sec:last-words}) the right pictures\dots

  \item What happens with our approximations of $\Q{k}(x)$ if $\newalpha
    \to 0$ (ie.\ $C \to B$)?  There should be a problem since we
    otherwise would make statements about primes in very short
    intervals.  More generally, which condition on $\newalpha$ must be
    fulfilled such that the results are meaningful?

  \item Consider special cases:
    \begin{enumerate}
    \item integers with $k$ prime factors: $B=2$, $C=x$.
    \item smooth integers: $B=2$, $C=~const$.
    \item rough integers: $B=~const$, $C=x$.
    \end{enumerate}
    
    Behaviour of $\Qclean{}(x) := \presum{k}{0}{~log_{2}x}
    \Qclean{k}(x)$?

    Dependence of $\Qclean{k}$ on $k$?

  \item Corollaries?

  \item \showrefsfalse\showcitesfalse\the\oopslist.

  \item Post to arXiv.

    Send draft for comments (and publication suggestions) to
    \begin{itemize}
    \item Igor Shparlinski.
    \item Peter B\"urgisser.
    \item Peter Moree.
    \item Alexander May.
    \item Schindler?
    \end{itemize}
    Submit extended abstract  to
    \begin{itemize}
    \item ANTS?
    \end{itemize}
    Submit paper to
    \begin{itemize}
    \item Journal of Number Theory?
    \item Mathematics of Computation?
    \item Acta mathematica?
    \item Crelle (Journal f{\"{u}}r die reine und angewandte Mathematik)?
    \end{itemize}
  \item Future directions:
    \begin{enumerate}
    \item Exact formula for $\Qclean{k}$ or $\Q{k}$?
    \item Fast counting algorithm for $\Qclean{k}$ or $\Q{k}$?
      (Gourdon et al\@.  have counted $\Qclean[0,\infty]{2}(x)$ for
      some $x$ values)
    \end{enumerate}
  \end{itemize}
  \global\oopslist{\csname @gobble\endcsname}%
  \clearpage
  \begingroup
  \def\enoteformat{\rightskip0pt \leftskip0pt \parindent=1.8em
    \def\xxx{1}\expandafter\ifx\csname @theenmark\endcsname\xxx\else
    \addvspace{4ex}%
    \fi
    \leavevmode\llap{\makeenmark}}
  \theendnotes
  \endgroup
  \clearpage
  \tableofcontents
}
\fi

\title{Coarse-grained integers}
\subtitle{Smooth? Rough? Both!
  \ifshortversion
  \\(Extended abstract)
  \fi}

\ifanonymous
\author{
  Anonymous.\\
  Address.}
\else
\ifllncs
\author{Daniel Loebenberger
  \and
  Michael N{\"{u}}sken}
\institute{%
  b-it\\
  Universit{\"{a}}t Bonn\\
  Dahlmannstr. 2\\
  D53012 Bonn\\
  \email{daniel@bit.uni-bonn.de}, \email{nuesken@bit.uni-bonn.de}\\
  \homepage{http://cosec.bit.uni-bonn.de/}
}
\else
\author{Daniel Loebenberger\\
  b-it\\
  Universit{\"{a}}t Bonn\\
  Dahlmannstr. 2\\
  D53012 Bonn\\
  \email{daniel@bit.uni-bonn.de}\\
  \homepage{http://cosec.bit.uni-bonn.de/}
\and
  Michael N{\"{u}}sken\\
  b-it\\
  Universit{\"{a}}t Bonn\\
  Dahlmannstr. 2\\
  D53012 Bonn\\
  \email{nuesken@bit.uni-bonn.de}\\
  \homepage{http://cosec.bit.uni-bonn.de/}
}
\fi
\fi

\ifllncs\begin{document}\maketitle\fi

\begin{abstract}
  We count $\left] B, C \right]$-grained, $k$-factor integers which
  are simultaneously $B$-rough and $C$-smooth and have a fixed number
  $k$ of prime factors.  Our aim is to exploit explicit versions of
  the prime number theorem as much as possible to get good explicit
  bounds for the count of such integers.  This analysis was inspired
  by certain inner procedures in the general number field sieve.  The
  result should at least provide some insight in what happens there.

  We estimate the given count in terms of some recursively defined
  functions.  Since they are still difficult to handle, only another
  approximation step reveals their orders.

  Finally, we use the obtained bounds to perform numerical experiments
  that show how good the desired count can be approximated for the
  parameters of the general number field sieve in the mentioned
  inspiring application.
\end{abstract}

\begin{keywords}
  Smooth numbers, rough numbers, counting, prime number theorem,
  general number field sieve, RSA.
\end{keywords}

\begin{subject}
  11Axx,
  11N05,
  11N25
\end{subject}

\ifllncs\else\begin{document}\fi
\thispagestyle{empty}

Let us call an integer \emph{$\mathcal{A}$-grained} iff all its prime
factors are in the set $\mathcal{A}$.  Then an integer is
\emph{$C$-smooth} iff it is \hbox{$\left]0,C\right]$-grained}, and
\emph{$B$-rough} iff it is \hbox{$\left]B,\infty\right[$-grained}.
You may want to call an integer family \emph{coarse-grained} if it is
\hbox{$\left]B,C\right]$-grained} and the number of factors is
bounded.  We consider the number of integers up to a real positive
bound~$x$ that are $\left]B,C\right]$-grained and have a given
number~$k$ of factors, in formulae:
\begin{equation}
  \label{def:Qclean}
  \Qclean{k}(x)
  :=
  \#\Set{\vphantom{\left(\P{B,C}\right)^{k}}
    n\leq x;
    \exists p_{1}, \dots, p_{k} \in \P{B,C}\colon 
    n=p_{1}\cdots p_{k}
  }.
\end{equation}
We always assume that $C > B$, since otherwise the set under
consideration is empty.

Such numbers occur for example in an intermediate step in the number
field sieve when trying to factor large numbers.  During the sieving
integers are constructed as random values of carefully chosen
small-degree polynomials and made $B$-rough by dividing out all
smaller factors.  The remaining number is fed into an elliptic curve
factoring algorithm with a time limit that should allow to find
factors up to~$C$.  To tune the overall algorithm it is vital to know
the probability that the remaining number is $C$-smooth.  Assuming
that the polynomials output true random numbers the counting task we
deal with is the missing link since estimates for counting $B$-smooth
numbers are known, see for example the overview by \cite{gra08}. It is
difficult, however, to compare our results to existing ones, since we
consider integers with a \emph{fixed} number of factors.  This
restriction is completely sufficient for the application we have in
mind. To our knowledge these kind of investigations cannot be found in the 
literature. Our main result implies the following
\begin{corollary*}
  Fix $k\in\N_{\geq2}$, $\newalpha>0$ and $\epsilon>0$.  Then for
  large $B$ and $C=B^{1+\alpha}$ we have uniformly for $x \in
  [B^{k}(1+\epsilon), C^{k}(1 - \epsilon)]$ that
  \begin{align*}
    \Qclean{k}(x) &\in \bigTheta{ \frac{x}{~ln B}}.
    \qed
  \end{align*}
\end{corollary*}
Actually, we can describe a piecewise smooth function that
approximates $\Qclean{k}$ up to an additive error of order
$\bigO{\frac{x}{\sqrt{B}}}$ and the hidden constant depends on
$\epsilon$, $k$ and $\alpha$ only, see \ref{res:cleanapprox},
\ref{res:final-count-lambda}.

To avoid difficulties with non-squarefree numbers, instead of numbers we
count lists of primes
\begin{equation}
  \label{def:Q}
  \Q{k}(x)
  :=
  \#\Set{
    ( p_{1}, \dots, p_{k} ) \in \left(\P{B,C}\right)^{k};
    p_{1} \cdots p_{k} \leq x
  }.
\end{equation}
It turns out that there are anyways only a few non-squarefree numbers
counted by~$\Qclean{k}$, namely $k! \cdot \Qclean{k}(x) \approx
\Q{k}(x)$.  This would actually be an equality if there were no
non-squarefree numbers in the count.  We defer a precise treatment
\begin{fullversion}[until later.]
  until \ref{sec:squareful-negligible}.
\end{fullversion}

We head for determining precise bounds $\Q{k}(x)$ or $\Qclean{k}(x)$
that can be used in practical situations.  However, to understand
these bounds we additionally consider the asymptotical behaviours.
This is tricky since we have to deal with the three parameters $B$,
$C$ and $x$ simultaneously.  To guide us in considering different
asymptotics, we usually write $B = x^{\beta}$, $C = x^{\gamma}$ and
$\gamma = \beta (1+\newalpha)$.  In particular, $C = B^{1+\newalpha}$.
So we replace $(B,C,x)$ with $(x,\beta,\gamma)$ or $(x,\alpha,\beta)$,
similar to the considerations when counting smooth or rough numbers in
the literature.  Alternatively, it seems also natural to fix $x$
somehow in the interval $]B^{k},C^{k}]$ by introducing a parameter
$\xi$ by
$$
x = B^{k-\xi} C^{\xi}
  = B^{k+\xi \alpha}.
$$
Now the parameters are $(B,C,\xi)$ or $(B,\alpha,\xi)$.

As a first observation, note that $\Q{k}(x)$ is constantly $0$ for
$x<B^{k}$ and constantly $(\pi(C)-\pi(B))^{k}$ for $x\geq C^{k}$
and grows
monotonically when $x$ goes through $[B^{k},C^{k}]$.  
Here $\pi$ denotes the prime counting function.
For `middle' $x$-values the asymptotics can be derived from our main
result:
\medskip
\begin{corollary*}
  Fix $k\in\N_{\geq2}$, $\newalpha>0$ and $\epsilon>0$.  Then for
  large $B$ and $C=B^{1+\alpha}$ we have uniformly for $x \in
  [B^{k}(1+\epsilon), C^{k}(1 - \epsilon)]$ that
  \begin{align*}
    \Q{k}(x) &\in \bigTheta{ \frac{x}{~ln B}} = \bigTheta{
      \frac{x}{\beta ~ln x}}.
  \end{align*}
\end{corollary*}
Later, we specify a piecewise smooth function $\Qapprox{k}$ which
approximates $\Q{k}$ with an error of order
$\bigO{\frac{x}{\sqrt{B}}}$ where the hidden constants depend on
$\epsilon$, $k$ and $\alpha$ only, see \ref{res:approx},
\ref{res:order-basis}, and subsequent considerations.

The preceeding result is in contrast to the asymptotics at $x =
C^{k}$:
$$
\Q{k}(C^{k}) \approx \frac{C^{k}}{~ln^{k} C}
= \frac{B^{k(1+\newalpha)}}{(1+\newalpha)^{k}~ln^{k} B}
\in \bigTheta{ \frac{x}{\beta^{k} ~ln^{k} x}}.
$$
This observation is explained as follows: Note that roughly half of
the numbers up to~$x$ are in the interval $[\frac12 x, x]$ and
similarly for primes.  Thus the behaviour of those candidates largely
rule $\Frac{\Q{k}(x)}{x}$.  For an $x = B^{k-\xi}C^{\xi}$ with a fixed
$\xi \in \left]0,k\right[$, it is mostly determined by the requirement
that the counted numbers are $B$-rough, and we thus observe a
comparatively large fraction of $]B,C]$-grained numbers.  In the
extreme case $x = C^{k}$, most candidates are ruled out by the
requirement to be $C$-smooth and thus we see a much smaller fraction
of $]B,C]$-grained numbers.

The case that the intervals $]B^{k},C^{k}]$ are disjoint for
considered values~$k$ is especially nice, as then the number of prime
factors of a $B$-smooth and $C$-rough number $n$ can be derived from
the number~$n$.  So we assume in the entire paper that $C < B^{s}$ for
some fixed $s>1$.  (Clearly, we cannot have a fixed $s$ that grants
disjointness for all intervals.  But for the first few we can.)  In
the inspiring number field sieve application we have $C<
B^{1.\rounded{231}{902}}$.  This ensures that the intervals $\left]
  B^{k},C^{k} \right]$ for $k\leq 5$ are disjoint.

A further application is related to RSA.  \cite{decmor08} give
estimates for the number of RSA-integers.  They use one ad-hoc
definition proposed by B.~de~Weger.  However, it is not so clear which
numbers we should call RSA-integers.  A discussion and further
calculations to adapt our results to the different shape are needed.
We have treated these issues in \ifanonymous[anonymized
citation]\else\citet{loenus11a}\fi.

As our basic field of interest is cryptography and there the largest
occurring numbers are actually small in the number theorist's view, we
assume the Riemann hypothesis throughout the entire paper.  We use the
following version of the prime number theorem\endnote{%
  We use the following versions of the prime number theorem.
  \begin{theorem}[Prime number theorem]
    We list a few variants:
    \begin{enumerate}
    \item
      \label{pnt-legendre}
      \cite{che52}, conjectured by \cite{leg98}:
      $$
      \pi(x) \approx \frac{x}{~ln x}.
      $$
    \item
      \label{pnt-li}
      \cite{had96}, \cite{val96}, \cite{wal63}, conjectured by
      \cite{gau49}:
      $$
      \pi(x) \in ~Li(x) + \mathcal{O}\left(x \, \exp
        \left( -\frac{A(\ln x)^{3/5}}{(\ln \ln x)^{1/5}} \right) \right)
      \subset ~Li(x) + \mathcal{O}\left(\frac{x}{~ln^{k}x}\right)
      $$
      for any $k$.  The presently best known value for $A$ is
      $A=0.2098$ \citep[p.~566]{for02}\nocite{pin84,for02a}.  The only
      estimate that we found with a given `constant' has
      $A=\frac{1}{57}$ and $a = 11.88 ~ln^{3/5} x$ (well, that's not
      constant, but almost\dots), \cite{for02a} attributes this to
      Y. Cheng, \textit{Explicit estimates on prime numbers
        (pre-print)}. [Though there are several preprints on arXiV
      (\url{http://arxiv.org/abs/0810.2113v1},
      \url{http://arxiv.org/abs/0810.2102v3},
      \url{http://arxiv.org/abs/0810.2103v5}), we couldn't find that
      paper\dots]

      Here, the logarithmic integral $~Li$ is given by $~Li(x) :=
      \int_{2}^{x} \frac{\mathrm{d}t}{~ln t}$.
    \item
      \label{pnt-dusart} 
      \cite{dus98}:
      For $x>\spacednumber{355991}$ we have
      $$
      \frac{x}{~ln x} + \frac{x}{~ln^{2} x} < \pi(x)
      < \frac{x}{~ln x} + \frac{x}{~ln^{2} x} + 2.51 \frac{x}{~ln^{3} x}.
      $$
    \item
      \Cite{koc01}, \cite{sch76a}:
      If (and only if) the Riemann hypothesis holds then for $x\geq2657$
      $$
      \left| \pi(x) - ~Li(x) \right| < \frac{1}{8\pi} \sqrt{x} ~ln x.
      $$
    \item\label{pnt-cramer}
      \cite{cra35,cra37}:
      If primes were random then
      $$
      \left| \Pi(x) - ~Li(x) \right| < \sqrt{ \frac{2 x ~ln ~ln x}{~ln
          x}}
      $$
      asymptotically surely, and even more:
      $$
      ~prob\left(
        \limsup_{x\to\infty}
        \frac{\left|\Pi(x)-~Li(x)\right|}{\sqrt{2 ~ln ~ln x}
          \sqrt{\frac{x}{~ln x}}} = 1 
      \right) = 1.
      $$
    \end{enumerate}
  \end{theorem}
  \begin{corollary}
    For $x>\spacednumber{355991}$ we have
    \begin{enumerate}
    \item\label{pnt-dusart1} $\frac{x}{~ln x} < \pi(x) < 1.1
      \frac{x}{~ln x}$, and
    \item\label{pnt-dusart2} $\frac{x}{~ln x} + \frac{x}{~ln^{2} x} <
      \pi(x) < \frac{x}{~ln x} + 1.2 \frac{x}{~ln^{2} x}$.
    \end{enumerate}
  \end{corollary}
  \begin{proof}
    This follows directly from \ref{pnt-dusart}.
  \end{proof}
  
  Similar bounds hold for the logarithmic integral:
  \begin{lemma}
    \label{li-bounds}
    \begin{enumerate}
    \item\label{li-bound1} For $x\geq 705$ we have
      $
      \frac{x}{~ln x} < ~Li(x) < 1.1253 \frac{x}{~ln x}.
      $
    \item\label{li-bound1a} For $x\geq 18986$ we have
      $
      \frac{x}{~ln x} < ~Li(x) < 1.1 \frac{x}{~ln x}.
      $
    \item\label{li-bound2}
      For $x\geq 1302$ we have
      $
      \frac{x}{~ln x} + \frac{x}{~ln^{2} x} 
      <
      ~Li(x)
      <
      \frac{x}{~ln x} + 1.1691 \frac{x}{~ln^{2} x}.
      $
    \item\label{li-bound3}
      For $x\geq 40088$ we have
      $
      \frac{x}{~ln x} + \frac{x}{~ln^{2} x} + 2 \frac{x}{~ln^{3} x}
      <
      ~Li(x)
      <
      \frac{x}{~ln x} + \frac{x}{~ln^{2} x} + 2.3854 \frac{x}{~ln^{3} x}.
      $
    \end{enumerate}
  \end{lemma}
  Except for \ref{li-bound1a} the requirement on $x$ is only needed
  for the lower bound and the given upper bound holds for $x\geq 2$.
  \begin{proof}
    It is actually easy to see by $k$-fold integration by parts that
    $~Li(x) = \presum{j}{1}{k} \frac{j! x}{~ln^{j}x} +
    o\left(\frac{x}{~ln^{k} x}\right)$.  The remainder is a little
    calculus, see our file \url{li-bounds.mn}.
  \end{proof}
}:
\vspace*{.5\baselineskip}
\begin{namedtheorem}{Prime number theorem}[\Citealt{koc01}, \cite{sch76a}]
  \label{pnt-schoenfeld}
  \leavevmode\hfil
  If (and only if) the Riemann hypothesis holds then for $x\geq2657$
  $$
  \left| \pi(x) - ~li(x) \right| < \frac{1}{8\pi} \sqrt{x} ~ln x,
  $$
  where $~li(x) = \frint{\frac{1}{~ln t}}{t}{0}{x}$.
\end{namedtheorem}
We have numerically verified this inequality for $x \leq 2^{40}
\approx 1.1 \cdot 10^{12}$ based on \citeauthor{kul08}'s tables
\citep{kul08} and extensions built using the segmented siever
implemented by \citet{oli03}.
We are confident that we can extend this verification much further.
In the inspiring application we have $x<2^{37}$ and so for those $x$
we can take this theorem for granted even if the Riemann hypothesis
should not hold.

\condbreak{1cm}%
We arrive at the following description of the desired count:
\begin{theorem*}
  Let $B<C=B^{1+\alpha}$ with $\alpha \geq \frac{~ln B}{\sqrt{B}}$ and
  fix $k\geq 2$.  Then for any (small) $\epsilon>0$ and $B$ tending to
  infinity we have for $x \in \left[ B^{k} (1+\epsilon), C^{k}
    (1-\epsilon) \right]$ a value $\widetilde{a} \in \left[
    \frac{\alpha^{k-1} \clower{k}}{k! (1+\alpha)^k}, \frac{1}{k!}
  \right]$ with $\clower{k} = ~min\left(2^{-4} \frac{\epsilon^{k}}{k!},
    2^{-k} \frac{\epsilon^{k-1}}{(k-1)!}\right) $ such that
  $$
  \left|\Qclean{k}(x) - \widetilde{a} \frac{x}{~ln B}\right| \leq
  (2^k-1) \alpha^{k-2} (1+\alpha) \cdot \frac{x}{\sqrt{B}} + 2^{k-1}
  \frac{x}{B}.\qed
  $$
\end{theorem*}

Also without assuming the Riemann hypothesis we can achieve meaningful
results provided we use a good unconditional version of the prime
number theorem whose error estimate is at least in $\mathcal{O} \left(
  \frac{x}{~ln^{3} x} \right)$.  The famous work by
\cite{rossch62,rossch75} is not sufficient.  Yet, \cite{dus98}
provides an explicit error bound of order $\mathcal{O} \left(
  \frac{x}{~ln^{3} x} \right)$, and \cite{for02} provides explicit
error bounds of order $\mathcal{O}\left(x \, \exp \left( -\frac{A(\ln
      x)^{3/5}}{(\ln \ln x)^{1/5}} \right) \right)$ though this only
applies for $x$ beyond $10^{171}$ or even much later depending on $A$
and the $\mathcal{O}$-constant%
\begin{fullversion}[.]%
  , see \ref{res:pnt-ford}.
\end{fullversion}

\begin{shortversion}
  \ifanonymous[anonymized citation] \else\cite{loenus12a} \fi is the
  full version corresponding to this extended abstract.
\end{shortversion}

\section{The recursion}
\label{sec:recursion}

The essential basis for the analysis of the counting functions $\Q{k}$
is the following simple description.
\begin{lemma}
  \label{res:recursion}
  For all $k \in \N_{>0}$ we have the recursion
  $$
  \Q{k}(x) = \primesum{ \Q{k-1} (x/p_{k}) }{p_{k}}{B}{C}
  $$
  based on
  $$
  \Q{0}(x)
  =
  \begin{cases}
    0 & \text{if $x \in \left[0,1\right[$},\\
    1 & \text{if $x \in \left[1,\infty\right[$}.
  \end{cases}
  \iffullversion\else\qed\fi
  $$
\end{lemma}

\begin{fullversion}
  \begin{proof}
    In case $k>0$ we have
    \begin{align*}
      \Q{k}(x) &= \#\Set{ (p_{1}, \ldots, p_{k}) \in
        \left(\P{B,C}\right)^{k}; p_{1} \cdots p_{k} \leq x}
      \\
      &= \# \biguplus_{p \in \P{B,C}} \Set{ (p_{1}, \ldots, p_{k}) \in
        \left(\P{B,C}\right)^{k};
        \begin{aligned}
          p_{1} \cdots p_{k-1} &\leq x/p_{k},
          \\
          p_{k}&=p
        \end{aligned}
      }
      \\
      &= \primesum{\Q{k-1} (x/p_{k})}{p_{k} }{B}{C} .
    \end{align*}
    The case $k=0$ is immediate from the definition.
  \end{proof}

  From the definition \ref{def:Q} or from \ref{res:recursion}, it is clear that
  \begin{gather*}
    \Q{1}(x)
    =
    \begin{cases}
      0 & \text{if $x \in \left[0,B\right[$},\\
      \pi(x)-\pi(B) & \text{if $x \in \left[B,C\right[$},\\
      \pi(C)-\pi(B) & \text{if $x \in \left[C,\infty\right[$}.\\
    \end{cases}
  \end{gather*}
  This reveals that the case distinction in $\Q{0}$ leaves its traces on
  higher $\Q{k}$.  For further calculations it is vital that we make
  this precise.  This will enable us later to do our estimates.
\end{fullversion}
For $k\in\N$ we distinguish $k+2$~cases:
\begin{align*}
  \text{$x$ is in case $(k,-1)$}
  &\defiff x \in \left[0,B^{k}\right[,
  \\
  \text{$x$ is in case $(k,j)$}
  &\defiff x \in \left[B^{k-j}C^{j}, B^{k-1-j} C^{j+1}\right[,
  \\
  \text{$x$ is in case $(k,k)$} &\defiff x \in \left[C^{k},
    \infty\right[,
\end{align*}
where $j\in\N_{<k}$.  Note that most cases are characterized by the
exponent of $C$ at the left end of the interval.
\begin{lemma}
  \label{res:cases}
  For $k, j \in \N$, $0 \leq j < k$ and $x$ in case $(k,j)$, we have
  $$
  \Q{k}(x) =
  \primesum{\Q{k-1}(x/p_{k})}{p_{k} }{\frac{x}{B^{k-1-j}C^{j}}}{C}  + 
  \primesum{\Q{k-1}(x/p_{k})}{p_{k} }{B}{\frac{x}{B^{k-1-j}C^{j}}}
  $$
  where in the first sum $x/p_{k}$ is in case $(k-1,j-1)$ and in the
  second in case $(k-1,j)$.  For $j=-1$, and $j=k$ we do not split the
  sum, as then all $x/p_{k}$ are in one case anyways.  For $j=0$ the
  left part is zero, so that the splitting there is less visible.
  \iffullversion\else\qed\fi
\end{lemma}
\begin{fullversion}
  \begin{proof}
    We only have to verify that $x/p_{k} \in \left[B^{k-1-j}C^{j},
      B^{k-j-2}C^{j+1}\right[$ for $p_{k} \in
    \P{B,\Frac{x}{B^{k-1-j}C^{j}}}$ and $x \in \left[B^{k-j}C^{j},
      B^{k-1-j} C^{j+1}\right[$.  Similarly, the statement for the
    second sum is established.
  \end{proof}
  For example, we obtain
  \begin{align}
    \label{res:Q2}
    \Q{2}(x) &=
    \begin{cases}
      0 & \text{if $x\in \left[0, B^{2}\right[$},
      \\
      \primesum{ \primesum{ 1 }{p_{1}}{B}{\frac{x}{p_{2}}}
      }{p_{2}}{B}{\frac{x}{B}} & \text{if $x\in \left[B^{2}, B
          C\right[$},
      \\
      \begin{array}{@{}l@{}}
        \primesum{
          \primesum{
            1
          }{p_{1}}{B}{\frac{x}{p_{2}}}
        }{p_{2}}{\frac{x}{C}}{C}
        \\
        +
        \primesum{
          \primesum{
            1
          }{p_{1}}{B}{C}
        }{p_{2}}{B}{\frac{x}{C}}
      \end{array}
      & \text{if $x\in \left[B C, C^{2}\right[$},
      \\
      \primesum{ \primesum{ 1 }{p_{1}}{B}{C} }{p_{2}}{B}{C} & \text{if
        $x\in \left[C^{2}, \infty\right[$}.
      \\
    \end{cases}
  \end{align}
  So for $\Q{2}(x)$ we have four cases with four $2$-fold sums.  In
  general $\Q{k}(x)$ has $k+2$ cases with $2^{k}$ $k$-fold sums.
  Well, we better stop unfolding here.%
\endnote{Actually,
  \begin{align}
    \label{res:Q3}
    \Q{3}(x)
    &=
    \footnotesize
    \begin{cases}
      0
      &
      \text{if $x\in \left[0, B^{3}\right[$},
      \\&\\
      \primesum{
        \primesum{
          \primesum{
            1
          }{p_{1}}{B}{\frac{x}{p_{2} p_{3}}}
        }{p_{2}}{B}{\frac{x}{B p_{3}}}
      }{p_{3}}{B}{\frac{x}{B^{2}}}
      &
      \text{if $x\in \left[B^{3}, B^{2} C\right[$},
      \\&\\
      \begin{array}{@{}l@{}}
        \primesum{
          \primesum{
            \primesum{
              1
            }{p_{1}}{B}{C}
          }{p_{2}}{B}{\frac{x}{C p_{3}}}
        }{p_{3}}{B}{\frac{x}{B C}}
        \\
        +
        \primesum{
          \primesum{
            \primesum{
              1
            }{p_{1}}{B}{\frac{x}{p_{2} p_{3}}}
          }{p_{2}}{B}{\frac{x}{B p_{3}}}
        }{p_{3}}{\frac{x}{B C}}{C}
        \\
        +
        \primesum{
          \primesum{
            \primesum{
              1
            }{p_{1}}{B}{\frac{x}{p_{2} p_{3}}}
          }{p_{2}}{\frac{x}{C p_{3}}}{C}
        }{p_{3}}{B}{\frac{x}{B C}}
      \end{array}
      &
      \text{if $x\in \left[B^{2} C, B C^{2}\right[$},
      \\&\\
      \begin{array}{@{}l@{}}
        \primesum{
          \primesum{
            \primesum{
              1
            }{p_{1}}{B}{C}
          }{p_{2}}{B}{C}
        }{p_{3}}{B}{\frac{x}{C^{2}}}
        \\
        +
        \primesum{
          \primesum{
            \primesum{
              1
            }{p_{1}}{B}{C}
          }{p_{2}}{B}{\frac{x}{C p_{3}}}
        }{p_{3}}{\frac{x}{C^{2}}}{C}
        \\
        +
        \primesum{
          \primesum{
            \primesum{
              1
            }{p_{1}}{B}{\frac{x}{p_{2} p_{3}}}
          }{p_{2}}{\frac{x}{C p_{3}}}{C}
        }{p_{3}}{\frac{x}{C^{2}}}{C}
      \end{array}
      &
      \text{if $x\in \left[B C^{2}, C^{3}\right[$},
      \\&\\
      \primesum{
        \primesum{
          \primesum{
            1
          }{p_{1}}{B}{C}
        }{p_{2}}{B}{C}
      }{p_{3}}{B}{C}
      &
      \text{if $x\in \left[C^{3}, \infty\right[$}.
      \\
    \end{cases}
  \end{align}
  and
  \begin{align}
    \label{res:Q4}
    \Q{4}(x)
    &=
    \tiny
    \begin{cases}
      0
      &
      \text{if $x\in \left[0, B^{4}\right[$},
      \\&\\
      \primesum{
        \primesum{
          \primesum{
            \primesum{
              1
            }{p_{1}}{B}{\frac{x}{p_{2} p_{3} p_{4}}}
          }{p_{2}}{B}{\frac{x}{B p_{3} p_{4}}}
        }{p_{3}}{B}{\frac{x}{B^{2} p_{4}}}
      }{p_{4}}{B}{\frac{x}{B^{3}}}
      &
      \text{if $x\in \left[B^{4}, B^{3} C\right[$},
      \\&\\
      \begin{array}{@{}l@{}}
        \primesum{
          \primesum{
            \primesum{
              \primesum{
                1
              }{p_{1}}{B}{\frac{x}{p_{2} p_{3} p_{4}}}
            }{p_{2}}{B}{\frac{x}{B p_{3} p_{4}}}
          }{p_{3}}{B}{\frac{x}{B^{2} p_{4}}}
        }{p_{4}}{\frac{x}{B^{2} C}}{C}
        \\
        +
        \primesum{
          \primesum{
            \primesum{
              \primesum{
                1
              }{p_{1}}{B}{C}
            }{p_{2}}{B}{\frac{x}{C p_{3} p_{4}}}
          }{p_{3}}{B}{\frac{x}{B C p_{4}}}
        }{p_{4}}{B}{\frac{x}{B^{2} C}}
        \\
        +
        \primesum{
          \primesum{
            \primesum{
              \primesum{
                1
              }{p_{1}}{B}{\frac{x}{p_{2} p_{3} p_{4}}}
            }{p_{2}}{B}{\frac{x}{B p_{3} p_{4}}}
          }{p_{3}}{\frac{x}{B C p_{4}}}{C}
        }{p_{4}}{B}{\frac{x}{B^{2} C}}
        \\
        +
        \primesum{
          \primesum{
            \primesum{
              \primesum{
                1
              }{p_{1}}{B}{\frac{x}{p_{2} p_{3} p_{4}}}
            }{p_{2}}{\frac{x}{C p_{3} p_{4}}}{C}
          }{p_{3}}{B}{\frac{x}{B C p_{4}}}
        }{p_{4}}{B}{\frac{x}{B^{2} C}}
      \end{array}
      &
      \text{if $x\in \left[B^{3} C, B^{2} C^{2}\right[$},
      \\&\\
      \begin{array}{@{}l@{}}
        \primesum{
          \primesum{
            \primesum{
              \primesum{
                1
              }{p_{1}}{B}{C}
            }{p_{2}}{B}{\frac{x}{C p_{3} p_{4}}}
          }{p_{3}}{\frac{x}{C^{2} p_{4}}}{C}
        }{p_{4}}{B}{\frac{x}{B C^{2}}}
        \\
        +
        \primesum{
          \primesum{
            \primesum{
              \primesum{
                1
              }{p_{1}}{B}{\frac{x}{p_{2} p_{3} p_{4}}}
            }{p_{2}}{\frac{x}{C p_{3} p_{4}}}{C}
          }{p_{3}}{\frac{x}{C^{2} p_{4}}}{C}
        }{p_{4}}{B}{\frac{x}{B C^{2}}}
        \\
        +
        \primesum{
          \primesum{
            \primesum{
              \primesum{
                1
              }{p_{1}}{B}{C}
            }{p_{2}}{B}{C}
          }{p_{3}}{B}{\frac{x}{C^{2} p_{4}}}
        }{p_{4}}{B}{\frac{x}{B C^{2}}}
        \\
        +
        \primesum{
          \primesum{
            \primesum{
              \primesum{
                1
              }{p_{1}}{B}{C}
            }{p_{2}}{B}{\frac{x}{C p_{3} p_{4}}}
          }{p_{3}}{B}{\frac{x}{B C p_{4}}}
        }{p_{4}}{\frac{x}{B C^{2}}}{C}
        \\
        +
        \primesum{
          \primesum{
            \primesum{
              \primesum{
                1
              }{p_{1}}{B}{\frac{x}{p_{2} p_{3} p_{4}}}
            }{p_{2}}{B}{\frac{x}{B p_{3} p_{4}}}
          }{p_{3}}{\frac{x}{B C p_{4}}}{C}
        }{p_{4}}{\frac{x}{B C^{2}}}{C}
        \\
        +
        \primesum{
          \primesum{
            \primesum{
              \primesum{
                1
              }{p_{1}}{B}{\frac{x}{p_{2} p_{3} p_{4}}}
            }{p_{2}}{\frac{x}{C p_{3} p_{4}}}{C}
          }{p_{3}}{B}{\frac{x}{B C p_{4}}}
        }{p_{4}}{\frac{x}{B C^{2}}}{C}
      \end{array}
      &
      \text{if $x\in \left[B^{2} C^{2}, B C^{3}\right[$},
      \\&\\
      \begin{array}{@{}l@{}}
        \primesum{
          \primesum{
            \primesum{
              \primesum{
                1
              }{p_{1}}{B}{C}
            }{p_{2}}{B}{\frac{x}{C p_{3} p_{4}}}
          }{p_{3}}{\frac{x}{C^{2} p_{4}}}{C}
        }{p_{4}}{\frac{x}{C^{3}}}{C}
        \\
        +
        \primesum{
          \primesum{
            \primesum{
              \primesum{
                1
              }{p_{1}}{B}{\frac{x}{p_{2} p_{3} p_{4}}}
            }{p_{2}}{\frac{x}{C p_{3} p_{4}}}{C}
          }{p_{3}}{\frac{x}{C^{2} p_{4}}}{C}
        }{p_{4}}{\frac{x}{C^{3}}}{C}
        \\
        +
        \primesum{
          \primesum{
            \primesum{
              \primesum{
                1
              }{p_{1}}{B}{C}
            }{p_{2}}{B}{C}
          }{p_{3}}{B}{C}
        }{p_{4}}{B}{\frac{x}{C^{3}}}
        \\
        +
        \primesum{
          \primesum{
            \primesum{
              \primesum{
                1
              }{p_{1}}{B}{C}
            }{p_{2}}{B}{C}
          }{p_{3}}{B}{\frac{x}{C^{2} p_{4}}}
        }{p_{4}}{\frac{x}{C^{3}}}{C}
      \end{array}
      &
      \text{if $x\in \left[B C^{3}, C^{4}\right[$},
      \\&\\
      \primesum{
        \primesum{
          \primesum{
            \primesum{
              1
            }{p_{1}}{B}{C}
          }{p_{2}}{B}{C}
        }{p_{3}}{B}{C}
      }{p_{4}}{B}{C}
      &
      \text{if $x\in \left[C^{4}, \infty\right[$}.
      \\
    \end{cases}
  \end{align}
}
\end{fullversion}

Based on the intuition that a randomly selected integer~$n$ is prime
with probability~$\frac{1}{~ln n}$ we can replace
$\primesum{f(p)}{p}{B}{C}$ with $\frint{\frac{f(p)}{~ln p}}{p}{B}{C}$.
This directly leads to the approximation function.  We prefer however
to follow a better founded way to them which will also give
information about the error term.

\section{Using estimates}
\label{sec:using-estimates}

From the recursion for $\Q{k}(x)$ it is clear that we have to compute
terms like
$$
\primesum{f(p)}{p}{\intrangeleft}{\intrangeright}
\quad\text{or}\quad
\primesum{f(x/p)}{p}{\intrangeleft}{\intrangeright}
.
$$
To get good estimates for such a sum we follow the classic path, as
\cite{rossch62}: we rewrite the sum as a Lebesque-Stieltjes-integral
over the prime counting function $\pi(p)$.  Then we substitute $\pi(t)
= ~Li(t) + E(t)$, keeping in mind that we know good bounds on the
error term $E(t)$ by the \ref{pnt-schoenfeld}.\endnote{%
  Under the Riemann hypothesis we have $|E(t)| <
  \widehat{E}(t):=\frac{1}{8\pi} \sqrt{t} ~ln t$.  $\widehat{E}'(t) =
  \frac{1}{16\pi} \cdot \frac{2+~ln t}{\sqrt{t}}$, $\frint{
    \frac{\widehat{E}'(x)}{x} }{x}{}{t} = -\frac{1}{8\pi} \cdot
  \frac{4+~ln t}{\sqrt{t}}$.}  Finally, we integrate by parts,
estimate and integrate by parts back%
\begin{fullversion}[.]%
  :
  \begin{align*}
    &\primesum{f(p)}{p}{\intrangeleft}{\intrangeright}
    \\
    &= \frint{f(t)}{\pi(t)}{\intrangeleft}{\intrangeright}
    \\
    &= \frint{f(t)}{\,~Li(t)}{\intrangeleft}{\intrangeright} +
    \frint{f(t)}{E(t)}{\intrangeleft}{\intrangeright}
    \\
    &= \frint{\frac{f(t)}{~ln t}}{t}{\intrangeleft}{\intrangeright} +
    f(\intrangeright) E(\intrangeright) - f(\intrangeleft)
    E(\intrangeleft) -
    \frint{E(t)}{f(t)}{\intrangeleft}{\intrangeright}.
  \end{align*}
  The existence of all integrals follow from the existence of the
  first.  If the sum kernel $f$ is differentiable with respect to~$t$
  we can rewrite $\frint{E(t)}{f(t)}{\intrangeleft}{\intrangeright} =
  \frint{f'(t) E(t)}{t}{\intrangeleft}{\intrangeright}$.  Now we can
  use the estimate on the error term $E(t)$:
  \begin{align*}
    & \left| \primesum{f(p)}{p}{\intrangeleft}{\intrangeright} -
      \int_{\intrangeleft}^{\intrangeright} \frac{f(t)}{~ln t} ~d\! t
    \right|
    \\
    & \leq |f(\intrangeright)| \widehat{E}(\intrangeright) +
    |f(\intrangeleft)| \widehat{E}(\intrangeleft) +
    \int_{\intrangeleft}^{\intrangeright} |f'(t)| \widehat{E}(t) ~d\!
    t.
  \end{align*}
  What remains is, given the concrete $f$, to determine the occurring
  integrals.  For the counting functions $\Q{k}$ ---as one would
  guess--- this task is more and more complicated the larger $k$ is.
  Clearly, smoothness properties of~$f$ must be considered carefully.
\end{fullversion}


During all this we make sure that the involved functions stay
sufficiently smooth:
\begin{lemma}[Prime sum approximation]\label{res:recapprox}
  Let $f$, $\widetilde{f}$, $\widehat{f}$ be functions $\shortmap{
    \R_{>0} }{\R_{\geq0}}$ such that $\widetilde{f}$ and
  $\widehat{f}$ are piecewise continuous, $\widetilde{f}+\widehat{f}$
  is \weaklyisotone{}, and
  \begin{equation*}
    \left| f(x) - \widetilde{f}(x) \right| \leq \widehat{f}(x)
  \end{equation*}
  for $x \in \R_{>0}$.  Further, let
  $\widehat{E}(p)$ 
  be a positive valued, \weaklyisotone{}, smooth function of~$p$
  bounding $|\pi(p)-~Li(p)|$ on $[\intrangeleft,\intrangeright]$.
  (For example, under the Riemann hypothesis we can take
  $\widehat{E}(p) = \frac{1}{8\pi} \sqrt{p} ~ln p$ provided $p \geq
  2657$.)  Then for $x\in\R_{> 0}$
  \begin{equation*}
    \Biggl|
    \primesum{f(x/p)}{p}{\intrangeleft}{\intrangeright}
    -
    \quad\widetilde{g}(x)\quad
    \Biggr|
    \;\leq\;
    \widehat{g}(x)
  \end{equation*}
  where
  \begin{align*}
    \widetilde{g}(x) &= \frint{ \frac{\widetilde{f}(x/p)}{~ln p}
    }{p}{\intrangeleft}{\intrangeright},
    \\
    \widehat{g}(x)
    &= \frint{ \frac{\widehat{f}(x/p)}{~ln p}
    }{p}{\intrangeleft}{\intrangeright}
    + 2 (\widetilde{f}+\widehat{f})(x/\intrangeleft)
    \widehat{E}(\intrangeleft)
    + \frint{ \left( \widetilde{f} + \widehat{f} \right)(x/p) \widehat{E}'(p)
    }{p}{\intrangeleft}{\intrangeright}.
  \end{align*}
  Moreover, $\widetilde{g}$ and $\widehat{g}$ are piecewise
  continuous, and $\widetilde{g}+\widehat{g}$ is \weaklyisotone{}. \iffullversion\else\qed\fi
\end{lemma}
\begin{fullversion}[The proof will appear in the full version.]
  \begin{proof}
    The assumption immediately implies
    $$
    \left| \primesum{f(x/p)}{p}{\intrangeleft}{\intrangeright} -
      \primesum{\widetilde{f}(x/p)}{p}{\intrangeleft}{\intrangeright}
    \right| \leq
    \primesum{\widehat{f}(x/p)}{p}{\intrangeleft}{\intrangeright}.
    $$
    Using the techniques just sketched we obtain
    $\primesum{\widetilde{f}(x/p)}{p}{\intrangeleft}{\intrangeright} =
    \frint{\frac{\widetilde{f}(x/p)}{~ln
        p}}{p}{\intrangeleft}{\intrangeright} +\linebreak[3]
    \frint{\widetilde{f}(x/p)}{E(p)}{\intrangeleft}{\intrangeright}$
    with $E(p) = \pi(p) - ~Li(p)$.  Shifting the second term to the
    correspondingly transformed error bound, we now have
    \begin{align*}
      \left| \primesum{f(x/p)}{p}{\intrangeleft}{\intrangeright} -
        \underbrace{\frint{\frac{\widetilde{f}(x/p)}{~ln
              p}}{p}{\intrangeleft}{\intrangeright}}_{=\widetilde{g}(x)}
      \right| &\leq \frint{\frac{\widehat{f}(x/p)}{~ln
          p}}{p}{\intrangeleft}{\intrangeright}
      + \frint{(\widehat{f}+\widetilde{f})(x/p)}
      {E(p)}{\intrangeleft}{\intrangeright}.
    \end{align*}
    This bound always holds, yet we still have to estimate $E(p)$ in it.
    Abbreviating $h(p) = (\widetilde{f} + \widehat{f})(x/p)$ we estimate
    the last integral:
    \begin{align*}
      \frint{h(p)}{E(p)}{\intrangeleft}{\intrangeright} &=
      h(\intrangeright) E(\intrangeright) - h(\intrangeleft)
      E(\intrangeleft) -
      \frint{E(p)}{h(p)}{\intrangeleft}{\intrangeright}
      \\
      &\leq
      h(\intrangeright) \widehat{E}(\intrangeright) + h(\intrangeleft)
      \widehat{E}(\intrangeleft) -
      \frint{\widehat{E}(p)}{h(p)}{\intrangeleft}{\intrangeright}
      \\
      &=
      \begin{aligned}[t]
        &+
        h(\intrangeright) \widehat{E}(\intrangeright) + h(\intrangeleft)
        \widehat{E}(\intrangeleft)
        \\
        &- h(\intrangeright) \widehat{E}(\intrangeright) +
        h(\intrangeleft) \widehat{E}(\intrangeleft)
        \\
        &+ \frint{h(p)}{\widehat{E}(p)}{\intrangeleft}{\intrangeright}.
      \end{aligned}
    \end{align*}
    For the inequality we use that $h(p)$ is \emph{\weaklyantitone{}}
    in~$p$.  Collecting gives the claim.

    Finally, $\widetilde{g}$ and $\widehat{g}$ being obviously piecewise
    continuous it remains to show that $\widetilde{g} + \widehat{g}$ is
    \weaklyisotone{}.  By the (defining) equation
    \begin{align*}
      (\widetilde{g} + \widehat{g})(x) &= \frint{
        (\widetilde{f}+\widehat{f})(x/p) \left( \frac{1}{~ln p} +
          \widehat{E}'(p) \right) }{p}{\intrangeleft}{\intrangeright}
      +
      2 (\widetilde{f}+\widehat{f})(x/\intrangeleft)
      \widehat{E}(\intrangeleft)
    \end{align*}
    this follows from $\widetilde{f}+\widehat{f}$ and $\widehat{E}$
    being \weaklyisotone{}.
  \end{proof}
\end{fullversion}

\section{Approximations}
\label{sec:approximations}

Based on \ref{res:recapprox} we recursively define approximation
functions~$\Qapprox{k}$ and error bounding functions~$\Qerror{k}$.
Recursive application of \ref{res:recapprox} will lead to a good
estimate
in \ref{res:approx}.  For the understanding we further need to
determine the asymptotic order of the functions defined here, which we
start in the remainder of this section.  \ref{res:order-basis} relates
the functions~$\Qapprox{k}$ and $\Qerror{k}$ to some easier manageable
functions.  These in turn are computed or estimated%
\begin{fullversion}[.]
  , respectively, in \ref{sec:solve-Kapprox} and
  \ref{sec:estimate-Kerror}.
\end{fullversion}

\begin{definition}
  \label{def:Qapprox}%
  For $x\geq 0$ we define
  $$
  \Qapprox{0}(x) := \Q{0}(x), \qquad   \Qerror{0}(x) := 0
  $$
  \label{def:Qerror}%
  and recursively for $k>0$
  \begin{align*}
    \Qapprox{k}(x) &:= \frint{\frac{\Qapprox{k-1}(x/p_{k})}{~ln
        p_{k}}}{p_{k}}{B}{C},
    \\
    \Qerror{k}(x) &:= \frint{\frac{\Qerror{k-1}(x/p_{k})}{~ln
        p_{k}}}{p_{k}}{B}{C}
    \\
    &\phantom{:=\;}+ 2 \left(\Qapprox{k-1}+\Qerror{k-1}\right)(x/B)
    \widehat{E}(B)
    \\
    &\phantom{:=\;}+ \frint{ \left( \Qapprox{k-1}+\Qerror{k-1}
      \right)(x/p_{k}) \widehat{E}'(p_{k})}{p_{k}}{B}{C}.
  \end{align*}
\end{definition}
These functions now describe the behavior of $\Q{k}$ nicely:
\begin{theorem}
  \label{res:approx}
  Given $x \in \R_{> 0}$ and $k \in \N$. Then the inequality
  $$
  \left|\Q{k}(x) - \Qapprox{k}(x)\right| \leq \Qerror{k}(x)
  $$
  holds.
\end{theorem}
\begin{proof}
  Using \ref{res:recapprox} the claim together with the fact that
  $\Qapprox{k}+\Qerror{k}$ is \weaklyisotone{} follow simultaneously
  by induction on~$k$ based on $\Qapprox{0} = \Q{0}$ and $\Qerror{0} =
  0$.
\end{proof}

\begin{fullversion}
  In order to give a first impression we calculate $\Qapprox{1}$ and
  $\Qerror{1}$.  Analogous to \ref{res:cases}, we split the
  integration at $\Frac{x}{B^{k-1-j}C^{j}}$ so that the parts fall
  entirely into case $(k-1,j-1)$ or into case $(k-1,j)$:
  $$
  \Qapprox{k}(x) =
  \frint{\Qapprox{k-1}(x/p_{k})}{p_{k}}{\Frac{x}{B^{k-1-j}C^{j}}}{C} +
  \frint{\Qapprox{k-1}(x/p_{k})}{p_{k}}{B}{\Frac{x}{B^{k-1-j}C^{j}}}.
  $$
  Also for $\Qerror{k}$ this can be done, simply split the occurring
  integrals at $\Frac{x}{B^{k-1-j}C^{j}}$.\endnote{Writing $\Qboth{k} :=
    \Qapprox{k} + \Qerror{k}$, we obtain analogously to \ref{res:cases}
    \begin{align*}
      \Qerror{k}(x) = &
      \frint{\frac{\Qerror{k-1}(x/p_{k})}{~ln p_{k}}}{p_{k}}{\Frac{x}{B^{k-1-j}C^{j}}}{C}
      +\frint{\frac{\Qerror{k-1}(x/p_{k})}{~ln p_{k}}}{p_{k}}{B}{\Frac{x}{B^{k-1-j}C^{j}}}\\
      &+ \frint{ \Qboth{k-1}(x/p_{k}) \widehat{E}'(p_{k})}{p_{k}}{\Frac{x}{B^{k-1-j}C^{j}}}{C}\\
      &+ \frint{ \Qboth{k-1}(x/p_{k}) \widehat{E}'(p_{k})}{p_{k}}{B}{\Frac{x}{B^{k-1-j}C^{j}}}\\
      &+ 2 \Qboth{k-1}(B) \widehat{E}(B)
    \end{align*}
  }
  Now unfolding the recursive definition of $\Qapprox{1}$ gives:
  \begin{align*}
    \Qapprox{1}(x) &=
    \begin{cases}
      0 & \text{if $x\in \left[0, B\right[$},
      \\
      \frint{ \frac{1}{~ln p_{1}} }{p_{1}}{B}{x} & \text{if $x\in
        \left[B, C\right[$},
      \\
      \frint{ \frac{1}{~ln p_{1}} }{p_{1}}{B}{C} & \text{if $x\in
        \left[C, \infty\right[$},
      \\
    \end{cases}
    \\
    \Qerror{1}(x) &=
    \begin{cases}
      \begin{aligned}[t]
        & 0 
      \end{aligned}
      & \text{if $x \in \left[0,B\right[$},
      \\
      \widehat{E}(x) + \widehat{E}(B)
      & \text{if $x \in \left[B,C\right[$},
      \\
      \widehat{E}(C) + \widehat{E}(B)
      & \text{if $x \in \left[C,\infty\right[$},
      \\
    \end{cases}
  \end{align*}
  which corresponds exactly to the approximation of the prime counting
  function by the logarithmic integral $~Li$.  In case $k=2$ we obtain
\end{fullversion}
\begin{shortversion}
  Unfolding the definition and splitting occuring integrals at
  $\Frac{x}{B^{k-1-j}C^{j}}$ in the spirit of \ref{res:cases} we
  obtain
\end{shortversion}
\begin{align}
  \label{res:Qapprox2}
  \Qapprox{2}(x)
  &=
  \begin{cases}
    0
    &
    \text{if $x\in \left[0, B^{2}\right[$},
    \\
    \frint{
      \frint{
        \frac{1}{~ln p_{1} ~ln p_{2}}
      }{p_{1}}{B}{\frac{x}{p_{2}}}
    }{p_{2}}{B}{\frac{x}{B}}
    &
    \text{if $x\in \left[B^{2}, B C\right[$},
    \\
    \begin{array}{@{}l@{}}
      \frint{
        \frint{
          \frac{1}{~ln p_{1} ~ln p_{2}}
        }{p_{1}}{B}{\frac{x}{p_{2}}}
      }{p_{2}}{\frac{x}{C}}{C}
      \\
      +
      \frint{
        \frint{
          \frac{1}{~ln p_{1} ~ln p_{2}}
        }{p_{1}}{B}{C}
      }{p_{2}}{B}{\frac{x}{C}}
    \end{array}
    &
    \text{if $x\in \left[B C, C^{2}\right[$},
    \\
    \frint{
      \frint{
        \frac{1}{~ln p_{1} ~ln p_{2}}
      }{p_{1}}{B}{C}
    }{p_{2}}{B}{C}
    &
    \text{if $x\in \left[C^{2}, \infty\right[$}.
  \end{cases}
\end{align}
\begin{fullversion}
  This is now exactly the transformed version of \ref{res:Q2}, as
  announced there.
  You may ask where you can find a display of the error term
  corresponding to \ref{res:Qapprox2}.  Well, we have computed it.
  But the resulting terms are so complex that we didn't really learn
  much from it.  Here is an expression for $x \in
  \left[B^{2},BC\right[$:
  \begin{align*}
    \Qerror{2}(x
    ) &=
    \begin{aligned}[t]
      &\preint{p}{B}{\frac{x}{B}} \preint{q}{B}{\frac{x}{p}} \left(
        \frac{\widehat{E}'(q)}{~ln p} + \frac{\widehat{E}'(p)}{~ln q}
        + \widehat{E}'(p) \widehat{E}'(q) \right) \postint \postint
      \\
      &+ 4 \widehat{E}(B) \preint{p}{B}{\frac{x}{B}} \frac{1}{~ln p}
      \postint
      \;+\; 4 
      \widehat{E}(B)
      \widehat{E}(\Frac{x}{B})
    \end{aligned}
  \end{align*}
\end{fullversion}
Though this term is still handleable, it becomes apparent that things
get more and more complicated with increasing $k$.  We escape from
this issue by loosening the bonds and weakening our bounds slightly.
The first aim will be to obtain easily computable terms while
retaining the asymptotic orders, the second aim will be to still
retain meaningful bounds for the fixed values $B=1100 \cdot 10^{6}$,
$C=2^{37}-1$, $k\in\Set{2,3,4}$ from our inspiring application.

Our next task is to describe the orders of $\Qapprox{k}$ and
$\Qerror{k}$.  The main problem in an exact calculation is that most
of the time we cannot elementary integrate a function with a logarithm
occurring in the denominator.  But using $B\leq p \leq C$ we can obtain
a suitably good approximation instead by replacing $\frac{1}{~ln p}$
with $\frac{1}{~ln B}$ in the integrals.  At this point we start using
$C \leq B^{s}$ and rewrite $C = B^{1+\newalpha}$ where $\newalpha$ is a new
parameter (bounded by $s-1$).  For the time being you can consider
$\newalpha$ as a constant, but actually we make no assumption on it.
This leads to the following
\begin{definition}
  \label{def:Kapprox}
  For $x\geq 0$ we let
  $$
  \Kapprox{0}(x) := \Q{0}(x), \qquad   \Kerror{0}(x) := 0,
  $$
  and recursively for $k>0$
    \label{def:Kerror}
  \begin{align*}
    \Kapprox{k}(x) &:= \frint{\frac{\Kapprox{k-1}(x/p_{k})}{~ln
        B}}{p_{k}}{B}{C},
    \\
    \Kerror{k}(x) &:= \frint{\frac{\Kerror{k-1}(x/p_{k})}{~ln
        B}}{p_{k}}{B}{C}
    \\
    &\phantom{:=\;}+ 2
    \left(\Kapprox{k-1}+\Kerror{k-1}\right)\treatparameter(x/B)
    \widehat{E}(B)
    \\
    &\phantom{:=\;}+ \frint{ \left( \Kapprox{k-1}+\Kerror{k-1}
      \right)\treatparameter(x/p_{k}) \widehat{E}'(p_{k})}{p_{k}}{B}{C}.
  \end{align*}
\end{definition}

\noindent We observe that $~ln B \leq ~ln p_{k} \leq ~ln C = (1+\newalpha) ~ln B$.
Thus we obtain $\frint{\frac{f(p)}{~ln p}}{p}{B}{C} \in \left[
  \frac{1}{1+\newalpha} , 1 \right] \frint{\frac{f(p)}{~ln B}}{p}{B}{C}$
for any positive integrable function $f$.  By induction on $k$ we
obtain\endnote{%
  For $k=0$ the assertion holds by definition.  Thus consider $k>0$.
  We obtain
  \begin{align*}
    \Qapprox{k}(x)
    &=
    \frint{ \frac{\Qapprox{k-1}(x/p)}{~ln p} }{p}{B}{C}
    \\
    &\in
    \left[\frac{1}{1+\newalpha},1\right]
    \frac{1}{~ln B} \frint{ \Qapprox{k-1}(x/p) }{p}{B}{C}
    \\
    &\subset%
    \left[\frac{1}{(1+\newalpha)^{k}},1\right]
    \underbrace{
      \frac{1}{~ln B}
      \frint{ \Kapprox{k-1}(x/p) }{p}{B}{C}
    }_{=
      \Kapprox{k}(x)
    }
  \end{align*}
  and similarly
  \begin{align*}
    \Qerror{k}(x)
    &=
    \begin{aligned}[t]
      &\frint{ \frac{\Qerror{k-1}(x/p)}{~ln p} }{p}{B}{C}
      \\
      &+ 2 \left(\Qapprox{k-1}+\Qerror{k-1}\right)(x/B) \widehat{E}(B)
      \\
      &+ \frint{ \left( \Qapprox{k-1}+\Qerror{k-1} \right)(x/p_{k})
        \widehat{E}'(p_{k})}{p_{k}}{B}{C}
    \end{aligned}
    \\
    &\in
    \left[\frac{1}{1+\newalpha},1\right]
    \begin{aligned}[t]
      \Biggl(
      &\frac{1}{~ln B} \frint{ \Qerror{k-1}(x/p) }{p}{B}{C}
      \\
      &+ 2 \left(\Qapprox{k-1}+\Qerror{k-1}\right)(x/B)
      \widehat{E}(B)
      \\
      &+ \frint{ \left( \Qapprox{k-1}+\Qerror{k-1} \right)(x/p_{k})
        \widehat{E}'(p_{k})}{p_{k}}{B}{C}
      \Biggr)
    \end{aligned}
    \\
    &\subset
    \left[\frac{1}{(1+\newalpha)^{k}},1\right]
    \underbrace{
      \begin{aligned}[t]
        \Biggl(
        &\frac{1}{~ln B} \frint{ \Kerror{k-1}(x/p) }{p}{B}{C}
        \\
        &+ 2 \left(\Kapprox{k-1}+\Kerror{k-1}\right)\treatparameter(x/B)
        \widehat{E}(B)
        \\
        &+ \frint{ \left( \Kapprox{k-1}+\Kerror{k-1}
          \right)\treatparameter(x/p_{k}) \widehat{E}'(p_{k})}{p_{k}}{B}{C}
        \Biggr)
      \end{aligned}
    }_{=\Kerror{k}(x)}
  \end{align*}
}
\begin{theorem}
  \label{res:order-basis}
  Write $C=B^{1+\newalpha}$ and fix $k \in \N_{>0}$.  Then for $x \in
  \R_{>0}$ we have
  \begin{align*}
    \Qapprox{k}(x) &\in \left[\frac{1}{(1+\newalpha)^{k}},1\right]
    \Kapprox{k}(x),
    \\
    \Qerror{k}(x) &\in \left[\frac{1}{(1+\newalpha)^{k}},1\right]
    \Kerror{k}(x).
\qed
  \end{align*}
\end{theorem}

\noindent In order to determine at least the asymptotic orders of $\Qapprox{k}$
and $\Qerror{k}$ (along with precise estimates) it remains to solve
the recursions for $\Kapprox{k}$ and $\Kerror{k}$, or at least to
estimate these functions.
As a first step, we rewrite all integrals in \ref{def:Kapprox} in
terms of $\rho$ defined by $p_{k} = B^{1+\rho\newalpha}$.
\begin{namedtheorem*}{\protect\ref{def:Kapprox} continued}
  Rewrite $x = B^{k+\xi\newalpha}$ with a new parameter $\xi$ and let
  $\Kapprox{k}<\xi> := \Kapprox{k}(B^{k+\xi\newalpha})$ and
  $\Kerror{k}<\xi> := \Kerror{k}(B^{k+\xi\newalpha})$.  (Actually, you
  may think of $f^{k}\treatparameter<\xi> := f^{k}(B^{k+\xi\newalpha})$
  for any family of functions $f^{k}$.)  
\end{namedtheorem*}
\begin{lemma}
  For $k>0$ we now have the recursion
  \label{def:Lapprox}
  \label{def:Lerror}
  \begin{align*}
    \Kapprox{k}<\xi> &= \newalpha \frint{ \Kapprox{k-1}<\xi-\rho>
      B^{1+\rho\newalpha} }{\rho}{0}{1},
    \\
    \Kerror{k}<\xi> &=
    \begin{aligned}[t]
      &\newalpha \frint{ \Kerror{k-1}<\xi-\rho>
        B^{1+\rho\newalpha}}{\rho}{0}{1}
      \\
      &+ 2 \left( \Kapprox{k-1}+\Kerror{k-1}
      \right)\treatparameter<\xi> \widehat{E}(B)
      \\
      &+ \newalpha ~ln B \frint{ \left( \Kapprox{k-1}+\Kerror{k-1}
        \right)\treatparameter<\xi-\rho>
        \widehat{E}'(B^{1+\rho\newalpha}) B^{1+\rho\newalpha}}{\rho}{0}{1}.
    \end{aligned}
  \end{align*}
  \null\vskip0pt\vskip-2\baselineskip\vskip0pt\qed
\end{lemma}

\begin{shortversion}
  \iffullversion\def\next{\section*}\else\let\next\section\fi
  \next{About the proofs}

  The remainder of the full version gives a detailed analysis of the
  behavior of the functions~$\Kapprox{k}$ and~$\Kerror{k}$.  Here we
  only describe the main results.

  After renormalizing by $\Kapprox[~norm]{k}<\xi> =
  \frac{1}{\newalpha^{k} B^{k+\xi\newalpha}} \Kapprox{k}<\xi>$ the recursion
  simplifies considerably and we can read off the upper bound
  $\Kapprox[~norm]{k}<\xi> \leq \frac{1}{\newalpha ~ln B}$.  Note that
  $B^{k+\xi\newalpha} = x$ so we essentially consider a function related
  to the ratio of $\left] B, C \right]$-grained integers.
  \begin{figure}[hb]
    \centering
    \psset{xunit=1cm,yunit=10cm}
\begin{pspicture}(-3em,-3ex)(7.6,.23)
  {\footnotesize
    \psset{linewidth=.4pt,tickstyle=bottom}
    \psaxes{->}(0,0)(-0.6,0)(7.6,.23)
    \psline(0,\lnormoneatone)(-3pt,\lnormoneatone)
    \uput[180](-3pt,0.2054729137){%
      $\frac{1-B^{-\gamma}}{\gamma ~ln B}$}
    \uput[90](5.6,3pt){$\xi$}
  }%
  \psset{linecolor=blue}
  \def\cutexpdef{
    /cutexp {
      dup 1 lt {
        pop 0 exch
      } {
        0 exch 
        1 exch 
        -1 add 1 exch 1 exch 
        {
          4 1 roll 
          dup 
          3 1 roll add 
          4 1 roll 
          exch dup 
          5 1 roll 
          mul exch div 
        } for 
        add exch 
      } ifelse 
      2.718281828 exch exp 
      exch sub 
    } def }%
  \def\lambdanormdef{%
    /lambdanorm {
      0 1 4 2 roll 
      exch dup floor 
      0 exch 1 exch 
      { 
        4 copy 
        3 2 roll 
        2 copy 
        4 1 roll exch 
        sub exch 1 add div 
        4 -1 roll mul 
        8 1 roll 
        7 2 roll 
        sub \glnB mul neg exch cutexp 
        mul add 
        4 1 roll exch 
      } for 
      pop exch pop 
      \glnB neg exch exp div 
    } def }%
  \psplot[linecolor=cyan]{-0.5}{1.5}{\cutexpdef \lambdanormdef x 1 lambdanorm}
  \psplot[linecolor=green]{-0.5}{2.5}{\cutexpdef \lambdanormdef x 2 lambdanorm}
  \psplot[linecolor=brown]{-0.5}{3.5}{\cutexpdef \lambdanormdef x 3 lambdanorm}
  \psplot[linecolor=orange]{-0.5}{4.5}{\cutexpdef \lambdanormdef x 4 lambdanorm}
  \psplot[linecolor=red]{-0.5}{5.5}{\cutexpdef \lambdanormdef x 5 lambdanorm}
  \psplot[linecolor=magenta]{-0.5}{6.5}{\cutexpdef \lambdanormdef x 6 lambdanorm}
  \psplot[linecolor=blue]{-0.5}{7.5}{\cutexpdef \lambdanormdef x 7 lambdanorm}
\end{pspicture}
    \caption{$\Kapprox[~norm]{k}<\xi>$ for $\xi \in
      \left[-\frac12,k+\frac12\right]$ and $k=1,2,3,4,5,6,7$}
    \label{fig:lambdanorm}
  \end{figure}
  Moreover, we are able to give explicit formulae for $\Kapprox{k}$.
  They allow us to also derive a corresponding lower bound, which we
  need to make statements about relative errors, if we avoid the
  boundaries.  As we prove in the full version and
  \ref{fig:lambdanorm} asserts, the functions
  $\Kapprox[~norm]{k}<\xi>$ increase on an interval
  $[0,\xi^{k}_{\frac12}]$ and decrease on the remaining interval
  $[\xi^{k}_{\frac12},k]$.
  \begin{theorem}
    \label{res:Kapprox-bounds}
    For any $k \geq 1$ we have
    $$
    \Kapprox[~norm]{k}<\xi> \leq \capprox{k} :=
    \begin{cases}
      \frac{1}{\newalpha ~ln B} & \text{if $k>0$},
      \\
      1 & \text{if $k=0$}.
    \end{cases}
    $$
    Assume $\newalpha ~ln B \geq ~max\treatparameter( ~ln 2, \frac{~ln
      16}{k} )$.  Then for any $\epsilon \in \left] 0,1 \right]$ with
    $\xi^{k}_{\frac12} \in \left[ \epsilon, k-\epsilon \right]$ or
    $k<3$ there is a $\clower{k}>0$ such that for $\xi \in
    \left[ \epsilon, k-\epsilon \right]$
    $$
    \Kapprox[~norm]{k}<\xi> \geq \frac{\clower{k}}{\newalpha ~ln B}.
    $$
    Here we can choose 
    \begin{align*}
      \clower{k} = ~min\left(
        \quad
        2^{-4}
        \cdot
        \frac{\epsilon^{k}}{k!},
        \quad
        2^{-k}
        \cdot
        \frac{\epsilon^{k-1}}{(k-1)!}
        \quad
      \right).
    \end{align*}
  \end{theorem}
  Summing up we obtain:
  \begin{corollary}
    For any $\epsilon \in \left]0,1\right]$ and any $k\geq 1$ we have
    uniformly for $\xi \in \left[\epsilon,k-\epsilon\right]$
    $$
    \Kapprox[~norm]{k}<\xi> \in
    \Theta\treatparameter(\frac{1}{\newalpha ~ln B}).
    \qed
    $$
  \end{corollary}

  Studying the behavior of the error bounds $\Kerror{k}$ is
  considerably more difficult and we finally abstained from proving
  explicit formulae and only proved upper bounds on them.
  \begin{theorem}\label{res:Kerror-upperbound}
    For any $k$ and $\xi \in \R$ we have
    $$
    \Kerror[~norm]{k}<\xi> \leq \cerror{k}
    $$
    where the values $\cerror{k}$ are defined recursively by
    $$
    \cerror{k}
    := 
    \cerror{k-1}
    +
    \frac{4+3~ln B}{8\pi\newalpha\sqrt{B}}
    \cboth{k-1}
    $$
    for $k\geq 3$ based on $\cerror{0} := 0$, $\cerror{1} :=
    \frac{(2+\newalpha) ~ln B}{8\pi\newalpha\sqrt{B}}$, and
    \begin{align*}
      \cerror{2} &:= 
      \frac{6+3\newalpha}{8\pi\newalpha^{2}}
      \frac{1}{\sqrt{B}}
      +
      \frac{4+3~ln B}{8\pi\newalpha\sqrt{B}}
      \cboth{1}
      \\
      &=
      \frac{9+3\newalpha}{8\pi\newalpha^{2}}
      \frac{1}{\sqrt{B}}
      + \frac{1}{2\pi\newalpha^{2}\sqrt{B} ~ln B}
      +
      \frac{(2+\newalpha) (4+3~ln B) ~ln B}{64\pi^{2}\newalpha^{2}B}
      .
    \end{align*}
    If $\newalpha \geq \frac{~ln B}{\sqrt{B}}$ we have for $k\geq 2$
    and large $B$ the inequality
    $$
    \cerror{k} \leq \frac{(2^k-1) (1+\newalpha)}{\newalpha^{2} \sqrt{B}}.
    $$
    For $k=1$ the order of $\cerror{1}$ is necessarily slightly
    larger. More precisely we have
    $$
    \cerror{1} \leq \frac{(1+\alpha) ~ln B}{\newalpha \sqrt{B}}. \qed
    $$
  \end{theorem}
  Instead of defining $\cerror{2}$ and $\cerror{1}$ we could have left
  that to the recursion.  But the given values are smaller than the ones
  derived from the recursion based on $\cerror{0}$ only.  For $k=1$ the
  recursion would give $\frac{4+3~ln B}{8\pi\newalpha\sqrt{B}}$.  For
  $k\geq 2$ however the improvement due to these explicit settings is a
  factor of order $~ln B$.

  The full version also gives unconditional error bounds.  However, as
  this was not our primary interest we did not optimize this and only
  obtain meaningful results when the used error bound on $\left|
    \pi(x) - ~Li(x) \right|$ is of order $\bigO{ \frac{x}{~ln^{3} x}}$
  or better.
  \begin{theorem}
    \iffullversion\else\label{res:order-error-woRiemann}\fi%
    Assume that $\widehat{E}(x)$ bounds $ \left| \pi(x) - ~Li(x)
    \right|$ and $\widehat{E}(x)/x$ is \weaklyantitone{} for $x>x_{0}$
    and the relative error decreases fast, ie.
    $$
    \frac{
      \widehat{E}(x^{1+\newalpha}) \frac{~ln x^{1+\newalpha}}{x^{1+\newalpha}}
    }{
      \widehat{E}(x) \frac{~ln x}{x}
    }
    =
    \frac{
      (1+\newalpha) \widehat{E}(x^{1+\newalpha})
    }{
      x^{\newalpha} \widehat{E}(x)
    }
    $$
    is bounded under the chosen behavior of $\newalpha$ (which is allowed
    to depend on $x$ here).  Then for any $k\geq 2$ and $B > x_{0}$
    large we have
    $$
    \Kerror[~norm]{k}<\xi> \in
    \mathcal{O}\treatparameter({
      \left( 1+\frac{k-1}{\newalpha ~ln B} \right)
      \left( 1+\frac{1}{\newalpha ~ln B} \right)
      \frac{\widehat{E}(B) ~ln B}{B}
    })
    $$
    for $k\geq 2$.
    \qed
  \end{theorem}
  This is close to optimal, we only loose a factor of order $~ln B$ in
  the relative error compared to the used error bound in the prime
  number theorem.  Due to the lost $~ln B$ the result is only
  meaningful if the used error bound is at least of order
  $\bigO{\frac{x}{~ln^{l} x}}$ with $l\geq 3$, so \cite{rossch62} does
  not suffice.  From \cite{dus98} we can use $\widehat{E}(x) = 2.3854
  \frac{x}{~ln^{3}x}$ for $x>\spacednumber{355991}$, and obtain (for
  constant $\newalpha$)
  $$
  \frac{\Kerror[~norm]{k}<\xi>}{\Kapprox[~norm]{k}<\xi>} \in
  \mathcal{O}\treatparameter(\frac{1}{~ln B}).
  $$
  
  That could be the entire story.  However, we still have the
  inspiring application in mind.  So we have a look at the quality of
  \ref{res:order-basis} when applied to it.  There
  $B=1100\cdot10^{6}$, $C=2^{37}-1$, $\newalpha = ~ln(C)/~ln(B)-1 =
  0.\rounded{231}{902}$, and the largest $k$ of interest is $k=4$.
  For these parameters we find
  $$
  \left[ \frac{1}{(1+\newalpha)^{k}}, 1 \right]
  \subset
  \left[ 0.434, 1 \right].
  $$
  That is a great loss when we try to enclose the function of interest
  in a small interval.  Actually, we can improve \ref{res:order-basis}
  for the price of a slightly more complicated recursion.  Instead of
  approximating $\mathring{\kappa}(\rho) = \frac{1}{~ln p}$ with
  $\mathring{\lambda}(\rho) =\frac{1}{~ln B}$ in the recursion for
  $\Kapprox{k}$, we use $\mathring{\eta}(\rho) = \myexp{ -~ln ~ln B -
    \rho ~ln(1+\newalpha) } = \frac{1}{(1+\newalpha)^{\rho} ~ln B}$, define
  $\Eapprox{k}$ analogous to $\Kapprox{k}$ and obtain
  \begin{theorem}
    \label{res:Eorder-basis}
    Write $C=B^{1+\newalpha}$ and fix $k \in \N_{>0}$.  Then for $x \in
    \R_{>0}$ we have
    \begin{align*}
      \Qapprox{k}(x) &\in 
      \left[
        \left(
          \frac{~ln(1+\newalpha)}{\newalpha}
          (1+\newalpha)^{\frac{1}{~ln(1+\newalpha)} - \frac{1}{\newalpha}}
        \right)^{k},
        1
      \right]
      \Eapprox{k}(x).
      \qed
    \end{align*}
  \end{theorem}
  In the light of the inspiring application we now find
  $$
  \left[
    \left(
      \frac{~ln(1+\newalpha)}{\newalpha}
      (1+\newalpha)^{\frac{1}{~ln(1+\newalpha)} - \frac{1}{\newalpha}}
    \right)^{k},
    1
  \right]
  \subset
  [0.978, 1].
  $$

  Finally, we have to deal with the difference between $\Qclean{k}$
  and $\frac{1}{k!}\Q{k}$ caused by integers that are not squarefree.
  We obtain
  \begin{lemma}
    \label{res:pi-kappa-comparison}
    We have 
    $
    \left| \Qclean{k}(x) - \frac{1}{k!}\Q{k}(x) \right| 
    < 2^{k-1} \frac{x}{B}.
    $\qed
  \end{lemma}
  Compared to the error bound in $\left| \Q{k}(x)-\Qapprox{k}(x)
  \right| \leq \Qerror{k}(x) \in \mathcal{O}(\frac{x}{\sqrt{B}})$ this
  is negligible for large $B$.  Here we assume $k\geq 2$ since the
  present observations are irrelevant for $k=1$, namely $\Qclean{1} =
  \Q{1}$.
  Now let $\Qcleanapprox{k}(x) := \frac{1}{k!} \Qapprox{k}(x)$,
  $\Qcleanerror{k}(x) := \frac{1}{k!} \Qerror{k}(x) + 2^{k-1}
  \frac{x}{B}$.  Combining with \ref{res:Kapprox-bounds} and
  \ref{res:Kerror-upperbound} we finally arrive at
  \ref{res:Brough+Csmooth}, the overall summary of our results.

\end{shortversion}

\begin{fullversion}
  
  \section{Solving the recursion for $\Kapprox{k}$}
  \label{sec:solve-Kapprox}

  To construct a useful description of the function $\Kapprox{k}$ we
  make a small excursion and consider the following family of piecewise
  polynomial functions.
  \begin{definition}[Polynomial hills]
    \label{def:hillpolynomials}
    Initially, define the integral operator $\mathcal{M}$ by
    $$
    (\mathcal{M}f)(\xi)
    =
    \frint{f(\xi-\rho)}{\rho}{0}{1}
    =
    \frint{f(\rho)}{\rho}{\xi-1}{\xi}
    $$
    for any integrable function $\shortmap[f]{\R}{\R}$.
    \label{def:m-recursion}%
    Now for $k\in\N_{>1}$ let the \emph{$k$-th polynomial hill} be
    $$
    \widetilde{m}^{k} := \mathcal{M} \widetilde{m}^{k-1}
    $$
    based on the rectangular function $\widetilde{m}^{1}$ given by
    $\widetilde{m}^{1}(\xi) = 1$ for $\xi\in \left[0,1\right[$ and
    $\widetilde{m}^{1}(\xi) = 0$ otherwise. 
  \end{definition}
  Actually, $\widetilde{m}^{1} = \mathcal{M} \widetilde{m}^{0}$ if we
  let $\widetilde{m}^{0} = \delta$ be the `left lopsided' Dirac delta
  distribution defined by its integral
  $\frint{\delta(t)}{t}{-\infty}{\xi}$ being $0$ for $\xi<0$ and $1$ for
  $\xi\geq 0$.  Contrastingly, the standard Dirac delta function is
  balanced and has $\frint{\delta(t)}{t}{-\infty}{0} =
  \frac{1}{2}$.\linebreak[3] However, we will stick to the lopsided
  variant throughout the entire paper.
  By $\mathcal{D}_{\xi}$ we denote the differential operator with
  respect to $\xi$ .

  \begin{lemma}[Polynomial hills]
    \label{res:polyhill}
    \begin{enumerate}
    \item\label{res:polyhill-outside} $\widetilde{m}^{k}(\xi) = 0$ for
      $\xi< 0$ or $\xi\geq k$.
    \item\label{res:polyhill-base} $\widetilde{m}^{k}(\xi) =
      \frac{1}{(k-1)!}\xi^{k-1}$ for $\xi\in \left[0,1\right[$ and
      $\widetilde{m}^{k}(\xi) = \frac{1}{(k-1)!}(k-\xi)^{k-1}$ for $\xi
      \in \left[k-1,k\right[$.
    \item\label{res:polyhill-form} $\widetilde{m}^{k}$ restricted to
      ${[j,j+1[}$ is a polynomial function of degree $k-1$ for $j \in
      \Set{0,\dots,k-1}$.  In particular, $\widetilde{m}^{k}$ is smooth
      for $\xi \in \R\setminus\Set{0,\dots,k}$.
    \item\label{res:polyhill-smooth} $\widetilde{m}^{k}$ is $(k-2)$-fold
      continuously differentiable.
    \item\label{res:polyhill-uniqueness} Conversely, the conditions
      \ref{res:polyhill-outside} through \ref{res:polyhill-smooth}
      uniquely determine $\widetilde{m}^{k}$.
    \item\label{res:polyhill-smooth-value} $\mathcal{D}_{\xi}^{i}
      \widetilde{m}^{k} = \mathcal{M} \mathcal{D}_{\xi}^{i}
      \widetilde{m}^{k-1}$ as long as $i\leq k-2$ and even for $i=k-1$
      when read for distributions.
    \item\label{res:polyhill-symmetry} The function $\widetilde{m}^{k}$
      is symmetric to $\frac{k}{2}$: $\widetilde{m}^{k}(\xi) =
      \widetilde{m}^{k}(k-\xi)$.
    \item\label{res:polyhill-growth} For $\xi \in \left[j,j+1\right[$
      (corresponding to case $(k,j)$) we have
      $$
      \mathcal{D}_{\xi}^{k-1} \widetilde{m}^{k}(\xi) = (-1)^{j} \cdot
      \binom{k-1}{j}.
      $$
    \item\label{res:polyhill-growth2} The next derivate can only be
      correctly described as a linear combination of Dirac delta
      distributions:
      $$
      \mathcal{D}_{\xi}^{k} \widetilde{m}^{k}(\xi) = \presum{j}{0}{k}
      (-1)^{j} \binom{k}{j} \cdot \delta(\xi-j).
      $$
    \item\label{res:polyhill-explicit} For any $0\leq l < k$ we obtain
      the following explicit description:
      \begin{align*}
        \mathcal{D}_{\xi}^{l} \widetilde{m}^{k} (\xi) &=
        \frac{1}{(k-1-l)!}  \presum{i}{0}{\floor{\xi}} \binom{k}{i}
        (-1)^{i} (\xi-i)^{k-1-l}
      \end{align*}
    \item\label{res:polyhill-convexcombination} Further,
      $\widetilde{m}^{k}(\xi) = \frac{\xi}{k-1} \widetilde{m}^{k-1}(\xi)
      + \frac{k-\xi}{k-1} \widetilde{m}^{k-1}(\xi-1)$ holds.
    \end{enumerate}
  \end{lemma}
  \begin{figure}[!t]
    \centering%
    \input{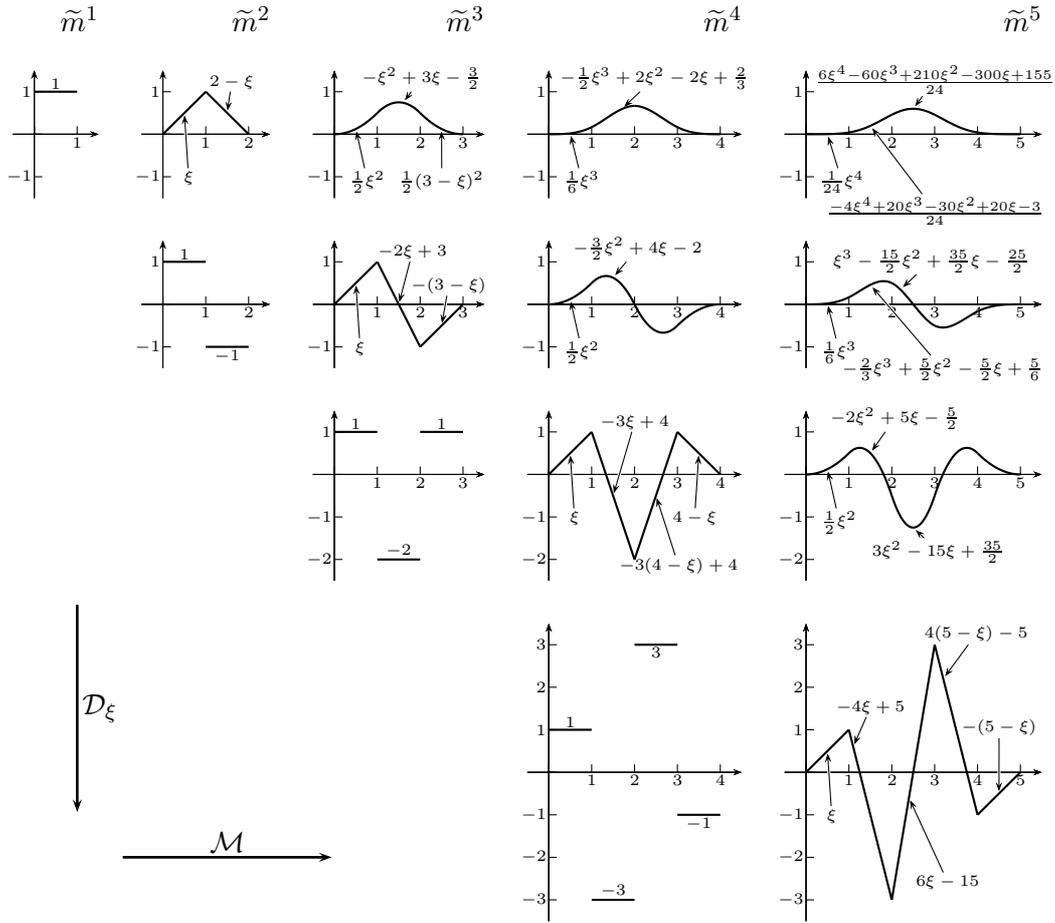}
    \caption{Graphs of the polynomial hills $\widetilde{m}^{k}$ for $k
      \leq 5$ and their derivatives}
    \label{fig:hillpoly}
  \end{figure}
  Curiosity: The function $\sqrt{k}\cdot \widetilde{m}^{k}$ can also
  be described as  the volume of a slice through a $(k-1)$-dimensional
  unit hypercube of thickness $\frac{1}{\sqrt{k}}$ orthogonal to a main
  diagonal.\endnote{Compare
    \protect\url{http://www.research.att.com/~njas/sequences/A011818}.}
  That's the same as saying it is $\frac{1}{\sqrt{k}}$ times the
  $(k-1)$-volume of the cut between a $k$-dimensional hypercube and
  the hyperplane $\presum{i}{1}{k} \rho_{i} = \xi$.\endnote{%
    We have \begin{align*} \widetilde{m}^{k}(\xi) &= \frint{
        \widetilde{m}^{k-1}(\xi-\rho_{k}) }{\rho_{k}}{0}{1}
      \\
      &= \frint{ \frint{ \widetilde{m}^{k-2}(\xi-\rho_{k}-\rho_{k-1})
        }{\rho_{k-1}}{0}{1} }{\rho_{k}}{0}{1}
      \\
      &= \int_{0}^{1} \cdots \int_{0}^{1}
      \widetilde{m}^{1}\treatparameter(\vphantom{\sum}\xi-\smash{\presum{i}{2}{k}}
      \rho_{i}) ~d\!\rho_{2} \dots
      ~d\!\rho_{k} 
      \\
      &= \int_{0}^{1} \cdots \int_{0}^{1} \int_{0}^{1}
      \widetilde{m}^{0}\treatparameter(\vphantom{\sum}\xi-\smash{\presum{i}{1}{k}}
      \rho_{i}) ~d\!\rho_{1} ~d\!\rho_{2} \dots ~d\!\rho_{k} =
      \int_{[0,1]^k}
      \delta\treatparameter(\vphantom{\sum}\xi-\braket{1}{\rho}) ~d\!
      \rho.  \end{align*} We now perform the orthogonal change of
    variables $A$ that maps the standard basis $(e_1, \ldots, e_k)$ to
    the basis $(z_1, \ldots, z_k)$, with $z_1 = \frac{1}{\sqrt{k}}
    \cdot \transpose{\left[1, \ldots, 1\right]}$ and $z_2$, \dots,
    $z_k$ such that the resulting basis is orthonormal.  Clearly, $A$
    is volume preserving, that is we can assume $~det A = 1$.  (If
    needed replace $z_{k}$ by $-z_{k}$.)  Then we can write any vector
    $v = \sum_{1 \leq i \leq k} \rho_i e_i = \sum_{1 \leq j \leq k}
    \zeta_j z_j$ with $z_j = \sum_{1 \leq i \leq k} A_{ij} e_i$. Thus
    $\rho_i = \sum_{1 \leq j \leq k} A_{ij} \zeta_j$, i.e. $\rho = A
    \zeta$.  
    Denote by $\indicate{\phi}$ the indicator function for the formula
    $\phi$. Substituting above gives \begin{align*}
      \widetilde{m}^{k}(\xi) &= \preint{\rho}{[0,1]^k}{}
      \delta\treatparameter(\vphantom{\sum}\xi-\braket{1}{\rho})
      \postint
      \\
      &= \preint{\zeta}{}{} \indicate{A\zeta \in [0,1]^{k}}
      \delta\treatparameter(\vphantom{\sum}\xi-\braket{1}{A \zeta})
      \postint
      \\
      &= \preint{\zeta}{}{} \indicate{A\zeta \in [0,1]^{k}}
      \delta\treatparameter(\vphantom{\sum}\xi-\braket{ A^{-1}
        \transpose{\left[1,\dots,1\right]}}{ \zeta}) \postint
      \\
      &= \preint{\zeta}{}{} \indicate{A\zeta \in [0,1]^{k}}
      \delta\treatparameter(\vphantom{\sum}\xi-\braket{ \sqrt{k} \cdot
        \transpose{[1,0,\dots,0]}}{ \zeta}) \postint
      \\
      &= \preint{\zeta}{}{} \indicate{A\zeta \in [0,1]^{k}}
      \delta\treatparameter(\vphantom{\sum}\xi-\sqrt{k} \cdot \zeta_1)
      \postint
      \\
      &= \preint{(\zeta_2, \ldots, \zeta_k)}{}{} \preint{\zeta_1}{}{}
      \indicate{A\zeta \in [0,1]^{k}}
      \delta\treatparameter(\vphantom{\sum}\xi-\sqrt{k} \cdot \zeta_1)
      \postint \postint
      \\
      &= \frac{1}{\sqrt{k}} \preint{\,(\zeta_2, \ldots, \zeta_k)}{}{}
      \indicate{A\transpose{\left[\frac{\xi}{\sqrt{k}}, \zeta_2,
            \ldots, \zeta_k\right]} \in [0,1]^k}
      \postint.  \end{align*} For the last step note that
    $\preint{\zeta_{1}}{}{} \delta\treatparameter(a-b\zeta_{1})
    \postint = \frac{1}{b} \preint{\eta}{}{}
    \delta\treatparameter(\eta) \postint$ by substituting
    $b\zeta_{1}-a = \eta$.  This is applicable if
    $A\transpose{\left[\frac{\xi}{\sqrt{k}}, \zeta_2, \ldots,
        \zeta_k\right]} \in ]0,1[^k$.  If to the contrary
    $A\transpose{\left[\frac{\xi}{\sqrt{k}}, \zeta_2, \ldots,
        \zeta_k\right]} \notin [0,1]^k$ then the inner integral
    vanishes.  The intermediate marginal cases are ignored by the
    outer integral.

    Finally, the last integral is exactly the $(k-1)$-volume of the
    intersection of the hypercube with the hyperplane $\xi =
    \braket{1}{\rho}$.
  }
  This last interpretation makes it obvious that
  $$
  \frint{\widetilde{m}^{k}(\xi)}{\xi}{0}{k} = 1,
  $$
  since this is the volume of the $k$-hypercube.\endnote{%
    As a corollary of \ref{res:polyhill-explicit} we obtain an explicit
    formula for the `normalized center slice volume' of a
    $k$-dimensional hypercube as described in
    \todo{\cite{slo-A011818}\url{http://www.research.att.com/~njas/sequences/A011818}}:
    \begin{align}
      &2^{n} n! ~vol\Set{x\in\left[0,1\right]^{n}; \presum{i}{1}{n}
        x_{i} \in \left[\frac{n}{2},\frac{n+1}{2}\right]}
      \nonumber\\
      &= 2^{n} n! \cdot \widetilde{m}^{n+1}\treatparameter( \frac{n}{2}
      )
      \nonumber\\
      &= \presum{i}{0}{\floor{n/2}} \binom{n+1}{i} (-1)^{i} (n-2i)^{n}.
    \end{align}
  }

  \begin{proof}
    From the definition \ref{res:polyhill-outside} through
    \ref{res:polyhill-symmetry} follow directly.
    Further, note that $\left(\mathcal{M} \mathcal{D}_{\xi}
      f\right)(\xi) = f(\xi)-f(\xi-1)$.  This is obvious from
    $(\mathcal{M}f)(\xi) = \frint{f(\rho)}{\rho}{\xi-1}{\xi}$.
    
    We prove \ref{res:polyhill-growth} by induction on $k$.  For $k=1$
    this is true by definition.  So consider $k>1$.  The recursion for
    $\widetilde{m}$ differentiated $k-1$~times yields for $\xi \in
    \left[j,j+1\right[$
    \begin{align*}
      \mathcal{D}_{\xi}^{k-1} \widetilde{m}^{k}(\xi) &=
      \mathcal{D}_{\xi} \mathcal{M} \mathcal{D}_{\xi}^{k-2}
      \widetilde{m}^{k-1}(\xi)
      \\
      &= \mathcal{D}_{\xi}^{k-2} \widetilde{m}^{k-1}(\xi) -
      \mathcal{D}_{\xi}^{k-2} \widetilde{m}^{k-1}(\xi-1)
      \\
      &= (-1)^{j} \binom{k-2}{j} - (-1)^{j-1} \binom{k-2}{j-1} =
      (-1)^{j} \binom{k-1}{j}.
    \end{align*}
    \ref{res:polyhill-growth2} follows similarly.\endnote{%
      \begin{align*}
        \mathcal{D}_{\xi}^{k} \widetilde{m}^{k}(\xi) &=
        \mathcal{D}_{\xi} \mathcal{M} \mathcal{D}_{\xi}^{k-1}
        \widetilde{m}^{k-1}(\xi)
        \\
        &= \mathcal{D}_{\xi}^{k-1} \widetilde{m}^{k-1}(\xi) -
        \mathcal{D}_{\xi}^{k-1} \widetilde{m}^{k-1}(\xi-1)
        \\
        &= \presum{j}{0}{k} \left( (-1)^{j} \binom{k-1}{j} \cdot
          \delta(\xi-j) - (-1)^{j-1} \binom{k-1}{j-1} \cdot
          \delta(\xi-j) \right)
        \\
        &= \presum{j}{0}{k} (-1)^{j} \binom{k}{j} \cdot \delta(\xi-j).
      \end{align*}
    }
    To prove \ref{res:polyhill-explicit} consider $h(\xi) =
    \frac{1}{(k-1)!}  \presum{i}{0}{\floor{\xi}} \binom{k}{i} (-1)^{i}
    (\xi-i)^{k-1}$.  Then $\mathcal{D}_{\xi}^{k-1} h =
    \mathcal{D}_{\xi}^{k-1} \widetilde{m}^{k}$ using
    \ref{res:polyhill-growth}. And obviously $\mathcal{D}_{\xi}^{l} h(0)
    = 0 = \mathcal{D}_{\xi}^{l} \widetilde{m}^{k}(0)$ for $0\leq l<k-1$,
    so that inductively (with falling $l$) we get $\mathcal{D}_{\xi}^{l}
    h =\mathcal{D}_{\xi}^{l} \widetilde{m}^{k}$.\endnote{Alternatively,
      we also find the following:
      \begin{itemize}
      \item $h(\xi) = g(\xi) = 0$ for $\xi>0$. (Note that the sum in $h$
        has only zero term for $i>k$.)  And clearly $h(\xi) = 0$ for
        $\xi\leq 0$.
      \item $h(\xi) = \frac{1}{(k-1)!} \xi^{k-1}$ for $\xi \in
        \left[0,1\right[$, and $h(\xi) = \frac{1}{(k-1)!} (k-\xi)^{k-1}$
        for $\xi\in\left[k-1,k\right[$.
      \item $h(\xi)$ restricted to $[j,j+1[$ is a polynomial function of
        degree at most $k-1$ for $j\in \Set{0,\dots,k-1}$.
      \item $h$ is $(k-2)$-fold continuously differentiable.
      \end{itemize}
      By \ref{res:polyhill-uniqueness} we thus obtain $\widetilde{m}^{k}
      = h$.  }

    As we do not need \ref{res:polyhill-uniqueness} and
    \ref{res:polyhill-convexcombination}, we leave these proofs to the
    interested reader.
  \end{proof}

  Most of the following is easier if we first renormalize
  $\Kapprox{k}$.
  So we let
  \begin{align}
    \Kapprox[~norm]{k}<\xi> &:= \frac{1}{\newalpha^{k} B^{k+\xi\newalpha}}
    \Kapprox{k}<\xi>.
  \end{align}
  The recursion for $\Kapprox{k}$ now turns into
  $$
  \Kapprox[~norm]{k} = \mathcal{M} \Kapprox[~norm]{k-1}.
  $$

  \begin{theorem}[Approximation order]
    \label{res:order-approx}
    For any $\xi \in \R$ we have
    \begin{align*}
      \Kapprox[~norm]{k}<\xi> &= \frint{B^{-\rho\newalpha}
        \widetilde{m}^{k}(\xi-\rho)}{\rho}{0}{\xi} =
      B^{-\xi\newalpha}\frint{B^{\rho\newalpha}
        \widetilde{m}^{k}(\rho)}{\rho}{0}{\xi},
      \\
      \Kapprox[~norm]{k}<\xi> &= \frac{1}{\left(-\newalpha ~ln
          B\right)^{k}}
      \begin{aligned}[t]
        & \left( \presum{i}{0}{\floor{\xi}} \binom{k}{i} (-1)^{i}
          B^{-(\xi-i)\newalpha} \right.
        \\
        &\hphantom{\left(\vphantom{\presum{i}{0}{\floor{\xi}}}\right.}\left.-
          \presum{l}{0}{k-1} \left( -\newalpha ~ln B \right)^{l} \cdot
          \mathcal{D}_{\xi}^{k-l-1}\widetilde{m}^{k}(\xi) \right),
      \end{aligned}
      \\
      \Kapprox[~norm]{k}<\xi> &=
      \begin{aligned}[t]
        & \presum{i}{0}{\floor{\xi}} \binom{k}{i} (-1)^{i} \frac{
          ~cutexp_{k} \left( -(\xi-i)\newalpha ~ln B \right)
        }{\left(-\newalpha ~ln B\right)^{k}},
      \end{aligned}
    \end{align*}
    where $~cutexp_{k}(\zeta) = ~exp(\zeta) - \presum{l}{0}{k-1} \frac{
      \zeta^{l} }{l!} = \sum_{l\geq k} \frac{ \zeta^{l} }{l!}$.  We can
    also express $~cutexp_{k}$ using the incomplete Gamma function
    $\Gamma(k,\zeta) = \frint{\eulernumber^{-\rho}
      \rho^{k-1}}{\rho}{\zeta}{\infty}$ by $~cutexp_{k}(\zeta) =
    ~exp(\zeta) - \frac{\Gamma(k,\zeta)}{~exp(-\zeta) \Gamma(k,0)}$.
  \end{theorem}

  \begin{proof}
    The definition for $\Kapprox{0}$ turns into
    $$
    \Kapprox[~norm]{0}<\xi> = \frint{ B^{-\rho\newalpha} \delta(\xi-\rho)
    }{\rho}{0}{\infty}.
    $$
    Now, since $\mathcal{M}$ commutes with this integration this
    immediately implies that
    $$
    \Kapprox[~norm]{k}<\xi> = \frint{ B^{-\rho\newalpha}
      \widetilde{m}^{k}(\xi-\rho) }{\rho}{0}{\infty}
    $$
    which is the first stated equality noting that $\widetilde{m}^{k}$
    is zero outside $[0,k]$.  By partial integration we obtain
    \begin{align*}
      \Kapprox[~norm]{k}<\xi>
      &=
      \frac{1}{\newalpha ~ln B} \widetilde{m}^{k}(\xi)
      - \frac{1}{\newalpha ~ln B} \frint{ B^{-\rho\newalpha}
        \mathcal{D}_{\xi}\widetilde{m}^{k}(\xi-\rho) }{\rho}{0}{\infty}
      \\\allowpagebreak
      &=
      \presum{i}{1}{k} \frac{(-1)^{i-1}}{ \left( \newalpha ~ln B
        \right)^{i}} \mathcal{D}_{\xi}^{i-1}\widetilde{m}^{k}(\xi)
      + \frac{(-1)^{k}}{ \left( \newalpha ~ln B \right)^{k}} \frint{
        B^{-\rho\newalpha} \mathcal{D}_{\xi}^{k}\widetilde{m}^{k}(\xi-\rho)
      }{\rho}{0}{\infty}.
    \end{align*}
    Using the description of $\mathcal{D}_{\xi}^{k} \widetilde{m}^{k}$
    from \ref{res:polyhill-growth2} the last integral turns into the
    claimed sum of the second stated equality.
    Expressing $\mathcal{D}_{\xi}^{k-l-1}\widetilde{m}^{k}(\xi-\rho)$
    using \ref{res:polyhill-explicit} and rearranging slightly yields
    the third equality.
  \end{proof}

  \endnotetext{%
    \label{en:lambdatilde}
    Explicit expressions for $\Kapprox{0}$, $\Kapprox{1}$,
    $\Kapprox{2}$:
    \begin{align*}
      \Kapprox{0}<\xi>/x &=
      \begin{cases}
        0 & \text{if $\xi\in \left]- \infty, 0\right[$},
        \\
        \frac{B^{- \xi \newalpha}}{1} & \text{if $\xi\in \left[0,
            \infty\right[$},
        \\
      \end{cases}
      \\
      \Kapprox{1}<\xi>/x &=
      \begin{cases}
        0 & \text{if $\xi\in \left]- \infty, 0\right[$},
        \\
        \frac{1}{~ln B} - \frac{B^{- \xi \newalpha}}{~ln B} & \text{if
          $\xi\in \left[0, 1\right[$},
        \\
        \frac{B^{\newalpha - \xi \newalpha}}{~ln B} - \frac{B^{- \xi
            \newalpha}}{~ln B} & \text{if $\xi\in \left[1, \infty\right[$},
        \\
      \end{cases}
      \\
      \Kapprox{2}<\xi>/x &=
      \begin{cases}
        0 & \text{if $\xi\in \left]- \infty, 0\right[$},
        \\
        \left(\xi \newalpha\right) \frac{1}{~ln B} - \frac{1}{~ln^{2} B} +
        \frac{B^{- \xi \newalpha}}{~ln^{2} B} & \text{if $\xi\in \left[0,
            1\right[$},
        \\
        \left(2 \newalpha - \xi \newalpha\right) \frac{1}{~ln B} +
        \frac{1}{~ln^{2} B} - 2 \frac{B^{\newalpha - \xi \newalpha}}{~ln^{2}
          B} + \frac{B^{- \xi \newalpha}}{~ln^{2} B} & \text{if $\xi\in
          \left[1, 2\right[$},
        \\
        \frac{B^{2 \newalpha - \xi \newalpha}}{~ln^{2} B} - 2 \frac{B^{\newalpha
            - \xi \newalpha}}{~ln^{2} B} + \frac{B^{- \xi \newalpha}}{~ln^{2}
          B} & \text{if $\xi\in \left[2, \infty\right[$},
        \\
      \end{cases}
    \end{align*}
  }


  We are going to estimate the estimation of the error in the next
  section.  To that aim we first need to
  estimate~$\Kapprox[~norm]{k}$.\endnote{We have
    $$
    \left.\mathcal{D}_{\xi}^{i} \Kapprox[~norm]{k}<\xi> \right|_{\xi=0}
    = 0
    $$
    as long as $i<k$.}
  If $k>0$ then, based on \ref{res:order-approx} and
  $\widetilde{m}^{k}(\xi)\leq 1$, we obtain the upper bound
  \begin{align*}
    \Kapprox[~norm]{k}<\xi> &= \frint{B^{-\rho\newalpha}
      \widetilde{m}^{k}(\xi-\rho)}{\rho}{0}{\infty}
    \leq \frac{1}{\newalpha ~ln B}.
  \end{align*}
  For $k=0$ we have $\Kapprox[~norm]{0}<\xi> = B^{-\xi\newalpha}$ for
  $\xi\geq0$ and so $\Kapprox[~norm]{0}<\xi> \leq 1$ will do for all
  $\xi \in \R$.

  We first describe the qualitative behaviour of $\Kapprox[~norm]{k}$.
  Actually its graph looks like a slighty biaswise hill.
  \begin{lemma}
    \label{res:Kapprox-qualitativebehavior}
    The function $\Kapprox[~norm]{k}<\xi> = B^{-\xi\newalpha} \frint{
      B^{\rho\newalpha} \widetilde{m}^{k}(\rho)}{\rho}{0}{\xi}$ is zero at
    $\xi=0$, positive at $\xi=k$, more precisely
    $$
    \Kapprox[~norm]{k}<k> = \left( \frac{1-B^{-\newalpha}}{\newalpha ~ln B}
    \right)^{k},
    $$
    and there is a position $\xi^{k}_{\frac12} \in \left]0,k\right]$ such
    that it is \weaklyisotone{} on $]0,\xi^{k}_{\frac12}[$ and
    \weaklyantitone{} on $]\xi^{k}_{\frac12},\xi[$.  Further, $\xi^{k}_{\frac12}
    \geq \frac{k}{2}$.
    \endnote{$\xi^{k}_{\frac12} \leq k-1$ is equivalent to a stronger
      version of \ref{res:Kapprox-rightlowerbound}. To see that note
      that this is equivalent to
      $\mathcal{D}_{\xi}\Kapprox[~norm]{k}<k-1> \leq 0$ and then use
      \ref{res:D-Kapproxnorm}\dots }
  \end{lemma}

  \begin{proof}
    First, inspecting
    \begin{align*}
      \Kapprox[~norm]{1}<\xi> &= \frint{ B^{-\rho\newalpha}
        \widetilde{m}^{1}(\xi-\rho)}{\rho}{0}{\xi}
      \\
      &=
      \begin{cases}
        0 & \text{if $\xi<0$},
        \\
        \frac{1-B^{-\xi\newalpha}}{\newalpha ~ln B} &\text{if $\xi \in
          \left[0,1\right]$},
        \\
        \frac{1-B^{-\newalpha}}{\newalpha ~ln B} B^{-(\xi-1)\newalpha} & \text{if
          $\xi>1$}.
      \end{cases}
      &%
      \rput(0,.5){\psset{xunit=1cm,yunit=5cm}
\begin{pspicture}(-0.5,-0.1)(1.6,.28)
  {\footnotesize
    \psset{linewidth=.4pt,tickstyle=bottom}
    \psaxes{->}(0,0)(-0.5,0)(1.6,.28)
    \psline(0,\lnormoneatone)(-3pt,\lnormoneatone)
    \uput[180](-3pt,0.2054729137){%
      $\scriptstyle\frac{1-B^{-\gamma}}{\gamma ~ln B}$}
  }%
  \psset{linecolor=cyan}
  \psplot{-0.4}{0}{0}
  \psplot{0}{1}{%
    1 2.718281828 x -\glnB mul exp sub \glnB div}
  \psplot{1}{1.4}{1 2.718281828 -\glnB exp sub \glnB
    div 2.718281828 x 1 sub -\glnB mul exp mul}
\end{pspicture}
}
    \end{align*}
    shows that for $k=1$ all claims hold with $\xi^{k}_{\frac12} = 1$.  So
    in the remainder of this proof we assume $k>1$.

    Next, compute $\Kapprox[~norm]{k}<k>$ inductively:
    \begin{align*}
      \Kapprox[~norm]{k}<k>
      &= B^{-k\newalpha} \frint{ B^{\rho \newalpha} \frint{
          \widetilde{m}^{k-1}(\rho-\rho_{k}) }{\rho_{k}}{0}{1}
      }{\rho}{0}{k}
      \\
      &= \underbrace{ B^{-\newalpha} \frint{ B^{\rho_{k} \newalpha}
        }{\rho_{k}}{0}{1} }_{=\frac{1-B^{-\newalpha}}{\newalpha ~ln B}} \cdot
      \underbrace{ B^{-(k-1)\newalpha} \frint{ B^{\tau \newalpha}
          \widetilde{m}^{k-1}(\tau) }{\tau}{0}{k-1}
      }_{=\Kapprox[~norm]{k-1}<k-1>}
      = \left( \frac{1-B^{-\newalpha}}{\newalpha ~ln B} \right)^{k} .
    \end{align*}
    Here we have substituted $\tau=\rho-\rho_{k}$ and collapsed the new
    integration interval $[-\rho_{k},k-\rho_{k}]$ for $\tau$ to
    $[0,k-1]$ since $\widetilde{m}^{k-1}$ vanishes on the difference
    Moreover, based on $\Kapprox[~norm]{k}<\xi> = B^{-\xi\newalpha} \frint{
      B^{\rho\newalpha} \widetilde{m}^{k}(\rho) }{\rho}{0}{\xi}$ we obtain
    \begin{align}
      \label{res:D-Kapproxnorm}
      \mathcal{D}_{\xi} \Kapprox[~norm]{k}<\xi> = -\newalpha ~ln B \cdot
      \Kapprox[~norm]{k}<\xi> + \widetilde{m}^{k}(\xi),
    \end{align}
    and infer that $\mathcal{D}_{\xi} \Kapprox[~norm]{k}<k> = -\newalpha
    ~ln B \left( \frac{1-B^{-\newalpha}}{\newalpha ~ln B} \right)^{k}$ is
    negative.
    
    Finally, compute the derivate of $\Kapprox[~norm]{k}$ differently
    \begin{align*}
      \mathcal{D}_{\xi}\Kapprox[~norm]{k}<\xi> &= \mathcal{D}_{\xi}
      \left( \frint{ B^{-\rho\newalpha}
          \widetilde{m}^{k}(\xi-\rho)}{\rho}{0}{\xi} \right)
      = B^{-\xi\newalpha} \frint{ B^{\rho\newalpha}
        \mathcal{D}_{\xi}\widetilde{m}^{k}(\rho) }{\rho}{0}{\xi}.
    \end{align*}
    The integral kernel $B^{\rho\newalpha}
    \mathcal{D}_{\xi}\widetilde{m}^{k}(\rho)$ is positive on
    $]0,\frac{k}{2}[$ and negative on $]\frac{k}{2},k[$.  Thus $\frint{
      B^{\rho\newalpha} \mathcal{D}_{\xi}\widetilde{m}^{k}(\rho)
    }{\rho}{0}{\xi}$ increases on $]0,\frac{k}{2}[$ and decreases on
    $]\frac{k}{2},k[$.  Since this term starts at zero, it is positive
    for some time, begins to decrease at $\frac{k}{2}$, traverses zero
    at some point $\xi^{k}_{\frac12}$ recalling that at $k$ the value is
    negative, and stays negative til $\xi = k$ since it continues to
    decrease.  Thus the sign of $\mathcal{D}_{\xi}
    \Kapprox[~norm]{k}<\xi>$ is positive on $]0,\xi^{k}_{\frac12}[$ and
    negative on $]\xi^{k}_{\frac12},k[$, and so $\Kapprox[~norm]{k}<\xi>$ is
    \weaklyisotone{} till $\xi^{k}_{\frac12}$ and \weaklyantitone{}
    afterwards.
  \end{proof}
  \begin{figure}[hbt]
    \centering
    \psset{xunit=1cm,yunit=10cm}
\begin{pspicture}(-3em,-3ex)(7.6,.23)
  {\footnotesize
    \psset{linewidth=.4pt,tickstyle=bottom}
    \psaxes{->}(0,0)(-0.6,0)(7.6,.23)
    \psline(0,\lnormoneatone)(-3pt,\lnormoneatone)
    \uput[180](-3pt,0.2054729137){%
      $\frac{1-B^{-\gamma}}{\gamma ~ln B}$}
    \uput[90](5.6,3pt){$\xi$}
  }%
  \psset{linecolor=blue}
  \def\cutexpdef{
    /cutexp {
      dup 1 lt {
        pop 0 exch
      } {
        0 exch 
        1 exch 
        -1 add 1 exch 1 exch 
        {
          4 1 roll 
          dup 
          3 1 roll add 
          4 1 roll 
          exch dup 
          5 1 roll 
          mul exch div 
        } for 
        add exch 
      } ifelse 
      2.718281828 exch exp 
      exch sub 
    } def }%
  \def\lambdanormdef{%
    /lambdanorm {
      0 1 4 2 roll 
      exch dup floor 
      0 exch 1 exch 
      { 
        4 copy 
        3 2 roll 
        2 copy 
        4 1 roll exch 
        sub exch 1 add div 
        4 -1 roll mul 
        8 1 roll 
        7 2 roll 
        sub \glnB mul neg exch cutexp 
        mul add 
        4 1 roll exch 
      } for 
      pop exch pop 
      \glnB neg exch exp div 
    } def }%
  \psplot[linecolor=cyan]{-0.5}{1.5}{\cutexpdef \lambdanormdef x 1 lambdanorm}
  \psplot[linecolor=green]{-0.5}{2.5}{\cutexpdef \lambdanormdef x 2 lambdanorm}
  \psplot[linecolor=brown]{-0.5}{3.5}{\cutexpdef \lambdanormdef x 3 lambdanorm}
  \psplot[linecolor=orange]{-0.5}{4.5}{\cutexpdef \lambdanormdef x 4 lambdanorm}
  \psplot[linecolor=red]{-0.5}{5.5}{\cutexpdef \lambdanormdef x 5 lambdanorm}
  \psplot[linecolor=magenta]{-0.5}{6.5}{\cutexpdef \lambdanormdef x 6 lambdanorm}
  \psplot[linecolor=blue]{-0.5}{7.5}{\cutexpdef \lambdanormdef x 7 lambdanorm}
\end{pspicture}
    \caption{$\Kapprox[~norm]{k}<\xi>$ for $\xi \in
      \left[-\frac12,k+\frac12\right]$ and $k=1,2,3,4,5,6,7$}
    \label{fig:lambdanorm}
  \end{figure}

  \pagebreak[3]
  \global\let\oldbeta\tau
  \begin{lemma}
    \label{res:Kapprox-leftlowerbound}
    Assume $k\geq 2$, $\newalpha ~ln B \geq \frac{~ln 16}{k}$,
    and $\xi \in \left[0,1\right]$.  Then we have
    \begin{align*}
      \Kapprox[~norm]{k}<\xi> & \geq \frac{ ~exp\treatparameter( -\frac{
          ~ln^{2} 4 }{ \newalpha ~ln B } ) }{\newalpha ~ln B} \cdot
      \frac{\xi^{k}}{k!}  .
    \end{align*}
  \end{lemma}
  The assumptions are already true for $C=2B$ when $k\geq 4$.  Note that
  this is rather sharp as for $\xi\in\left[0,1\right]$ we have
  $\Kapprox[~norm]{k}<\xi> \leq \frac{\xi^{k}}{k!}$.\endnote{%
    \begin{lemma*}
      Fix any $\theta_{1}>0$, $\theta_{2}>0$, let $\oldbeta_{1} =
      \frac{1-~exp(-\theta_{1})}{\theta_{1}}$, $\oldbeta_{2} =
      \frac{1-~exp(-\theta_{2})}{\theta_{2}}$, and
      assume that
      \begin{enumerate}
      \item\label{res:Kapprox-leftlowerbound-i} $0\leq \xi \leq
        \frac{k}{2}$, $\xi \leq 1$, and $\xi \leq \oldbeta_{1} \theta_{1} k
        = (1-~exp(-\theta_{1})) k$,
      \item\label{res:Kapprox-leftlowerbound-ii} $\oldbeta_{1} k \newalpha
        ~ln B \geq \frac{1}{1-~exp(-\theta_{2})} \geq 1$.
      \end{enumerate}
      Then we have
      $$
      \Kapprox[~norm]{k}<\xi> \geq \frac
      {~exp\treatparameter(-\frac{1}{\oldbeta_{2} \oldbeta_{1} \newalpha ~ln B})}
      {\newalpha ~ln B} \cdot \frac{\xi^{k}}{k!}  .
      $$

      Let $\theta_{1} = - ~ln \treatparameter(1-\frac{1}{k})$ and
      $\theta_{2} = - ~ln\treatparameter({ 1-\frac{\theta_{1}}{\newalpha
          ~ln B} })$.
    \end{lemma*}
  }

  \begin{proof}
    \endnote{Compare
      \url{Brough+Csmooth-lookforlowerbound-various.mn}.}
    We use the integral representation from \ref{res:order-approx} and
    estimate the polynomial hill part $\widetilde{m}^{k}(\xi-\rho)$ of
    its kernel by a simple piecewise constant function as indicated in
    the picture.  We obtain
    \begin{align*}
      &
      \begin{aligned}
        \Kapprox[~norm]{k}<\xi> &= \frint{ B^{-\rho\newalpha}
          \widetilde{m}^{k}(\xi-\rho)}{\rho}{0}{\xi}
        \\
        &\geq \frint{ B^{-\rho\newalpha}}{\rho}{0}{\epsilon} \cdot
        \widetilde{m}^{k}(\xi-\epsilon)
        \\
        &= \frac{1}{\newalpha ~ln B} \left(\big.1-
          ~exp\treatparameter({-\epsilon\newalpha ~ln B})\right)
        \widetilde{m}^{k}(\xi-\epsilon) .
      \end{aligned}
      &
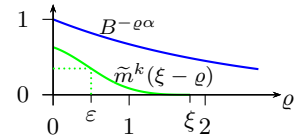
      \begin{pspicture}(-.5,-.5)(2.2,1.2)%
        \footnotesize%
        \psplot[linecolor=blue]{0}{2.7}{%
          2.71828 x -0.4 mul exp }
        \rput[lb](0.6,0.8){\scriptsize$B^{-\rho\newalpha}$}
        \psplot[linecolor=green]{0.8}{1.8}{%
          1.8 x sub 3 exp 6 div } \psplot[linecolor=green]{0}{0.8}{%
          1.8 x sub dup dup 1 -2 div mul 2 add mul -2 add mul 2 3 div
          add }
        \rput[lb](0.8,3pt){\scriptsize$\widetilde{m}^{k}(\xi-\rho)$}%
        {%
          \psset{linewidth=0.4pt,tickstyle=bottom}%
          \psaxes{->}(3,1.2)%
          \rput[lt](3,0){$\rho$}%
          \psline(1.8,0)(1.8,-3pt)%
          \rput[t](1.8,-6pt){$\xi$}%
          \psline(0.5,0)(0.5,-3pt)%
          \rput[t](0.5,-6pt){$\epsilon$}%
        }%
        \SpecialCoor
        \psline[linestyle=dotted,dotsep=0.04,linecolor=green](0.5,0)
        (!0.5 1.8 0.5 sub dup dup 1 -2 div mul 2 add mul -2 add mul 2 3
        div add) (!0.0 1.8 0.5 sub dup dup 1 -2 div mul 2 add mul -2 add
        mul 2 3 div add)
      \end{pspicture}
    \end{align*}
    This holds for any $\epsilon \in \left[0,\xi\right]$ since $0\leq
    \xi\leq \frac{k}{2}$ ensures that $\widetilde{m}^{k}$ is
    \weaklyisotone{}.  So we can optimize $\epsilon$ depending on $\xi$.
    We obtain a suitable value when setting $\epsilon$ by
    \begin{align}
      \label{def:leftlower-epsilon}
      1 - ~exp\treatparameter({-\epsilon\newalpha ~ln B}) &= \frac{\xi}{k}.
    \end{align}
    To make sure that now $\epsilon\leq\xi$ we use the following simple
    fact.
    \begin{fact}
      \label{res:map-1-exp-}
      For any $\theta>0$ and $\oldbeta = \frac{1-~exp(-\theta)}{\theta}$
      the map
      \begin{align*}
        &\map{[0,\theta]}{[0,\oldbeta\theta]}{z}{1-~exp(-z),} &&
        \rput(1.2,0.4){%
          \psset{unit=7.5mm}
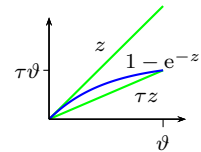
          \begin{pspicture}(-.7,-.5)(2.4,1.8) \scriptsize%
            \psplot[linecolor=green]{0}{2.0}{x}%
            \psplot[linecolor=green]{0}{2.0}{x 2 2.718281828 exch neg
              exp neg 1 add 2 div mul}%
            \psplot[linecolor=blue]{0}{2.0}{x 2.718281828 exch neg exp
              neg 1 add}%
            \psset{linewidth=0.4pt,labels=none,tickstyle=bottom,ticks=none}%
            \psaxes{->}(2.4,1.8)%
            \psline(2,0)(2,-3pt) \rput[t](2,-6pt){$\theta$}%
            \SpecialCoor%
            \psline(! 0 2 2.718281828 exch neg exp neg 1 add)%
            (! -0.1 2 2.718281828 exch neg exp neg 1 add)%
            \rput[r](! -0.15 2 2.718281828 exch neg exp neg 1
            add){$\oldbeta \theta$}%
            \rput[b](! 2 dup 2.718281828 exch neg exp neg 1
            add){$1-\eulernumber^{-z}$}%
            \rput[lt](! 1.5 dup 1 2.718281828 exch neg exp neg 1 add 2
            div mul){$\oldbeta z$}%
            \rput[rb](1.0,1.2){$z$}%
          \end{pspicture}
        }
      \end{align*}
      is bijective and \weaklyisotone{} and for $z \in
      \left[0,\theta\right]$ we have $\oldbeta z \leq 1-~exp(-z) \leq z$.
      \pushright{$\triangle$}
    \end{fact}
    Let $\theta_{1} > 0$ be such that $1-~exp(-\theta_{1}) =
    \frac{1}{k}$.  Namely, $\theta_{1} = -
    ~ln\treatparameter({1-\frac{1}{k}})$.  With $\oldbeta_{1} =
    \frac{1-~exp(-\theta_{1})}{\theta_{1}}$ then $1 = \oldbeta_{1}
    \theta_{1} k$ and $\frac{\xi}{k} \leq \oldbeta_{1} \theta_{1}$.
    Thus we have $\epsilon \newalpha ~ln B \in \left[0,\theta_{1}\right]$ so that
    \begin{gather}
      \label{eq:bounds-for-epsilon}
      \oldbeta_{1} \epsilon \newalpha ~ln B \leq 1 - ~exp(-\epsilon \newalpha ~ln
      B) = \frac{\xi}{k} \leq \epsilon \newalpha ~ln B.
    \end{gather}
    In particular, $\epsilon \leq \xi$ follows from $k \oldbeta_{1} \newalpha
    ~ln B = \frac{1}{\theta_{1}} \newalpha ~ln B \geq 1$.
    By \ref{res:map-1-exp-} with $\theta = ~ln 2$ we obtain $\oldbeta =
    \frac{1}{2 ~ln 2}$ and $ \left( 1-\frac{1}{k} \right)^{k} \geq
    ~exp(-\frac{1}{\oldbeta}) = ~exp(- ~ln 4)$ for all $k \geq 2$.  Thus $k
    \theta_{1} = -k ~ln \treatparameter( 1-\frac{1}{k} ) \leq ~ln 4$.
    Further, $\newalpha ~ln B \geq \frac{~ln 16}{k} > \frac{~ln 4}{k} \geq
    \theta_{1}$.
    This now implies $\epsilon \leq \xi$.
    
    Since $\xi \in \left[0,1\right]$ we have the explicit expression
    $\widetilde{m}^{k}(\xi-\epsilon) =
    \frac{(\xi-\epsilon)^{k-1}}{(k-1)!}$ and so
    \begin{align*}
      \Kapprox[~norm]{k}<\xi> &\geq \frac{\xi}{k \newalpha ~ln B}
      \widetilde{m}^{k}(\xi-\epsilon) =
      \frac{(1-\epsilon/\xi)^{k-1}}{\newalpha ~ln B} \cdot
      \frac{\xi^{k}}{k!}  .
    \end{align*}
    Since $\theta_{1} < \newalpha ~ln B$ we can define $\theta_{2}$ by
    $1-~exp(-\theta_{2}) = \frac{\theta_{1}}{\newalpha ~ln B}$, and
    according to \ref{res:map-1-exp-} let $\oldbeta_{2} =
    \frac{1-~exp(-\theta_{2})}{\theta_{2}}$.  Combining with
    \ref{eq:bounds-for-epsilon} gives us $ 
    \left( 1-\frac{\epsilon}{\xi} \right)^{k} 
    \geq \left( 1-\frac{1}{k \oldbeta_{1}\newalpha ~ln B} \right)^{k}
    = ~exp\treatparameter(-\frac{1}{\oldbeta_{2}\oldbeta_{1}\newalpha ~ln B})
    $ 
    and thus simplifies our above inequality to
    \begin{align}
      \label{res:Kapprox-leftlowerbound-stronger}
      \Kapprox[~norm]{k}<\xi> &\geq
      \frac{~exp\treatparameter(-\frac{1}{\oldbeta_{2} \oldbeta_{1} \newalpha ~ln
          B})}{\newalpha ~ln B} \cdot \frac{\xi^{k}}{k!}  .
    \end{align}
    By our choices
    \begin{align*}
      \frac{1}{\oldbeta_{1} \oldbeta_{2} \newalpha ~ln B} &=
      \let\gtrapprox\geq%
      k \theta_{2} \gtrapprox \frac{k\theta_{1}}{\newalpha ~ln B}
      \gtrapprox \frac{1}{\newalpha ~ln B}.
    \end{align*}
    Though these inequalities get equalities with $k\to\infty$, we need
    a precise estimate.
    Substituting $k$ in $-k ~ln\treatparameter(1-\frac{1}{k}) \leq ~ln
    4$ with $\frac{k \newalpha ~ln B}{~ln 4}$ yields
    \begin{align*}
      \frac{1}{\oldbeta_{2} \oldbeta_{1}} = k \newalpha ~ln B \cdot \theta_{2}
      &=
      - k \newalpha ~ln B ~ln\treatparameter({ 1 -
        \frac{k\theta_{1}}{k\newalpha ~ln B} })
      \\
      &\leq - \frac{k \newalpha ~ln B}{~ln 4} ~ln\treatparameter({ 1 -
        \frac{~ln 4}{k \newalpha ~ln B} }) ~ln 4
      \leq ~ln^{2} 4
    \end{align*}
    provided $\newalpha ~ln B \geq \frac{~ln 16}{k}$.
  \end{proof}


  Though we know the value of $\Kapprox[~norm]{k}<k>$, it is orders
  smaller than the above left lower bound.\endnote{%
    Together with the monotonicity behaviour of $\Kapprox[~norm]{k}$ we
    thus have that for $\xi \in \left[1,k\right]$
    $$
    \Kapprox[~norm]{k}<\xi> \geq ~min\treatparameter(
    \frac{~exp\treatparameter(-\frac{ 4 ~ln^{2} 2 } { \newalpha ~ln B})}{k!
      \newalpha ~ln B} , {\left( \frac{1-B^{-\newalpha}}{\newalpha ~ln B}
      \right)^{k}} ).
    $$
    Clearly, for large $k$ the first term is smaller.  With the
    parameter from our application the crossover point is between $k=8$
    and $k=9$.  Unfortunately, we are particularly interested in the
    behaviour when $B^{\newalpha}$ gets large, and a lower bound of order
    $\Omega\treatparameter({\frac{1}{(\newalpha ~ln B)^{k}}})$ is much
    weaker than we hoped for.} %
  Thus let us consider $\Kapprox[~norm]{k}$ on $[k-1,k]$.
  \begin{lemma}
    \label{res:Kapprox-rightlowerbound}
    If $k\geq 3$ then for $\xi \in \left[k-1,k\right]$ we have
    $$
    \Kapprox[~norm]{k}<\xi> \geq \frac{\left( 1-B^{-\newalpha}
      \right)^{k}}{\newalpha ~ln B} \cdot \frac{(k-\xi)^{k-1}}{(k-1)!}.
    $$
  \end{lemma}

  \begin{proof}
    By definition we have $\Kapprox[~norm]{k}<\xi> =
    \frac{\Kapprox{k}<\xi>}{\newalpha^{k} B^{k+\xi\newalpha}}$.
    \ref{res:order-approx}'s third description expresses
    $\Kapprox[~norm]{k}$ on $[k-1,k]$ as a sum of $k$ terms.  Adding the
    missing term $i=k$ we obtain $\frac{\Kapprox{k}<k>}{\newalpha^{k}
      B^{k+\xi\newalpha}}$, noting that $\Kapprox{k}$ is constant for
    $\xi>k$:
    \begin{align*}
      \Kapprox[~norm]{k}<\xi> &= \left( \frac{1-B^{-\newalpha}}{\newalpha ~ln
          B} \right)^{k} B^{(k-\xi)\newalpha} - \frac{ ~cutexp_{k} \left(
          (k-\xi)\newalpha ~ln B \right) }{\left(\newalpha ~ln B\right)^{k}}
    \end{align*}
    To check the claimed inequality we substitute $\tau = (k-\xi) \newalpha
    ~ln B \in \left[ 0, \newalpha ~ln B \right]$ (eliminating $\xi$):
    \begin{align*}
      \Kapprox[~norm]{k}<\xi>
      &= \frac{ \presum{l}{0}{k-1} \frac{ \tau^{l} }{l!}  - \left( 1 -
          \left( 1-B^{-\newalpha} \right)^{k} \right) \eulernumber^{\tau} }{ \left(
          \newalpha ~ln B \right)^{k} }.
    \end{align*}
    We have to show that this is at least
    $\frac{(1-B^{-\newalpha})^{k}}{\left(\newalpha ~ln B\right)^{k}}
    \frac{\tau^{k-1}}{(k-1)!}$
    which we rewrite to
    \begin{align*}
      \presum{l}{0}{k-2} \frac{ \tau^{l} }{l!}  \geq \left( 1-
        (1-B^{-\newalpha})^{k} \right) \left( \eulernumber^{\tau} -
        \frac{\tau^{k-1}}{(k-1)!}  \right) .
    \end{align*}
    Obviously $(1-B^{-\newalpha})^{k} \geq (1-\eulernumber^{-\tau})^{k}$ in our
    situation, with equality for $\tau = \newalpha ~ln B$ or $\xi=k-1$.
    Thus it suffices to show for any $\tau > 0$
    \begin{align}
      \label{eq:right-lower-extreme}
      \presum{l}{0}{k-2} \frac{ \tau^{l} }{l!}  \geq
      (1-(1-\eulernumber^{-\tau})^k) \left( \eulernumber^{\tau} - \frac{\tau^{k-1}}{(k-1)!}
      \right) .
    \end{align}
    The remaining proof proceeds in four steps:
    \begin{itemize}\parskip=0pt
    \item High case: $\presum{l}{0}{k-2}\frac{\tau^{l}}{l!} \geq k$.
    \item Low case: $\frac{1-\eulernumber^{-\tau}}{\tau} \geq
      \sqrt[k]{\frac{1+\sigma}{k!}}$, where $\frac{1}{\sigma}+1 \leq
      \frac{\eulernumber^{k-1} (k-1)!}{(k-1)^{k-1}}$.
    \item Covering: Fixing $\sigma := \frac{1}{\sqrt{2\pi(k-1)}-1}$
      these cases cover all $\tau > 0$ if $k\geq 4$.
    \item Brute-force: Prove \ref{eq:right-lower-extreme} for
      $k=3$. (Actually, for $k \in \Set{3,4,5,6,7,8,9}$.)
    \end{itemize}
    
  \item[High case:] Since $\eulernumber^{\tau} \geq
    \presum{l}{0}{k-2}\frac{\tau^{l}}{l!} \geq k$ in this case we have
    $\tau \geq ~ln k$.  Employing Bernoulli's inequality $(1-\Tau)^{k}
    \geq 1-k \Tau$ for $\Tau = \eulernumber^{-\tau}\leq 1$ we obtain
    \begin{align*}
      (1-(1-\eulernumber^{-\tau})^k) \left( \eulernumber^{\tau} - \frac{\tau^{k-1}}{(k-1)!}
      \right) &\leq k \eulernumber^{-\tau} \left( \eulernumber^{\tau} -
        \frac{\tau^{k-1}}{(k-1)!}  \right)
      \\
      &= k \left( 1 - \frac{\tau^{k-1}}{(k-1)!} \eulernumber^{-\tau} \right)
      \leq k \leq \presum{l}{0}{k-2}\frac{\tau^{l}}{l!}.
    \end{align*}
    As Bernoulli's inequality is good for $\Tau$ close to $0$ only, it
    is not surprising that this only gives a sufficient result for
    $\tau$ large enough.

  \item[Low case:] Assume $\frac{1-\eulernumber^{-\tau}}{\tau} \geq
    \sqrt[k]{\frac{1+\sigma}{k!}}$ where $\sigma$ is chosen such that
    $\frac{1}{\sigma}+1 < \frac{\eulernumber^{k-1}
      (k-1)!}{(k-1)^{k-1}}$.  (Since
    $\frac{1-\eulernumber^{-\tau}}{\tau}$ is
    \weaklyantitone{}\endnote{%
      From $\eulernumber^{-\tau} \leq 1 - \tau \eulernumber^{-\tau}$ we
      obtain $\mathcal{D}_{\tau} \frac{1-\eulernumber^{-\tau}}{\tau} =
      \frac{\eulernumber^{-\tau}+\tau \eulernumber^{-\tau} -
        1}{\tau^{2}} \leq 0$.}, this is always true on some interval
    $[0,\tau_{0}]$.)

    The condition on $\sigma$ implies that $\eulernumber^{\tau} -
    \left(\frac{1}{\sigma}+1\right) \frac{\tau^{k-1}}{(k-1)!}$ is
    non-negative for all $\tau > 0$: Consider $f_{0}(\tau) =
    \frac{\eulernumber^{\tau} (k-1)!}{\tau^{k-1}} -
    \left(\frac{1}{\sigma}+1\right)$.  Then $f_{0}'$ vanishes at $\tau =
    k-1$ only and thus $f_{0}$ is minimal there.  The assumption on
    $\sigma$ is precisely $f_{0}(k-1) \geq 0$.

    Further, note that $\frac{\tau^{k-2}}{(k-2)!} +
    \frac{\tau^{k-1}}{(k-1)!} \leq \eulernumber^{\tau} \leq \presum{l}{0}{k-1}
    \frac{\tau^{l}}{l!} + \eulernumber^{\tau} \frac{\tau^{k}}{k!}$.  We use this
    to obtain:
    \begin{align*}
      & \Bigl( 1-\underbrace{(1-\eulernumber^{-\tau})^k}_{\geq
        \frac{(1+\sigma)\tau^{k}}{k!}}  \Bigr) \underbrace{\left(
          \eulernumber^{\tau} - \frac{\tau^{k-1}}{(k-1)!}  \right)}_{\geq 0}
      \\
      &\leq \eulernumber^{\tau} - \frac{\tau^{k-1}}{(k-1)!}  - (1+\sigma)
      \frac{\tau^k}{k!}  \eulernumber^{\tau} + (1+\sigma) \frac{\tau^k}{k!}
      \frac{\tau^{k-1}}{(k-1)!}
      \\
      &\leq \presum{l}{0}{k-2} \frac{ \tau^{l} }{l!}  - \frac{\sigma
        \tau^k}{k!}  \eulernumber^{\tau} + (1+\sigma) \frac{\tau^k}{k!}
      \frac{\tau^{k-1}}{(k-1)!}
      \\
      &= \presum{l}{0}{k-2} \frac{ \tau^{l} }{l!}  - \sigma
      \frac{\tau^{k}}{k!}  \underbrace{ \left( \eulernumber^{\tau} - \left(
            \frac{1}{\sigma}+1 \right) \frac{\tau^{k-1}}{(k-1)!}
        \right) }_{\geq 0}
      \leq \presum{l}{0}{k-2} \frac{ \tau^{l} }{l!}  .
    \end{align*}
    Notice that the value of $\sigma$ only influences the set of values
    of $\tau$ that fall in this case.

  \item[Covering:]

    We choose $\sigma := \frac{1}{\sqrt{2\pi(k-1)}-1}$.  Recall
    Stirling's formula: For any $n>0$ there is a $\theta \in
    \left]0,1\right[$ such that $n! = \left( \frac{n}{\eulernumber{}}
    \right)^{n} \sqrt{2 \pi n} \eulernumber^{\frac{1}{12n+\theta}}$, see
    \cite{rob55}.  We thus have $\frac{1}{\sigma}+1 = \sqrt{2\pi(k-1)} <
    \frac{\eulernumber^{k-1} (k-1)!}{(k-1)^{k-1}}$ as required.  Further, we
    define the value $\tausplit$ where we split between the low and the
    high case:
    \begin{figure}[bt!]
      \centering
      \psset{xunit=2mm,yunit=7mm}
\begin{pspicture*}(-1.5em,-4ex)(21.5,5)
  \psset{linewidth=.4pt,tickstyle=bottom}%
  \footnotesize\psaxes[Dx=5]{->}(0,0)(0,0)(21.5,5)%
  \uput[-45](0,5){$\tau$}%
  \rput(15,4){$\tausplit$}
  \rput[b]{59}(6.1,2.75){low case boundary}
  \rput[tr]{15}(21.5,2.9){high case boundary}
  \uput[135](20,0){$k$}%
  \psset{linewidth=.8pt}%
  \savedata{\highcaseboundaryhangover}[{{2.1, 1011.977682}, {2.2,
      32.64142716}, {2.3, 10.22717027}, {2.4, 5.656905887}, {2.5,
      3.952529446}, {2.6, 3.116719627}, {2.7, 2.639389706}, {2.8,
      2.339761818}, {2.9, 2.139567895}, {3, 2.0}} ]
  \savedata{\highcaseboundary}[{{3, 2.0}, {4, 1.645751311}, {5,
      1.70870069}, {6, 1.831236896}, {7, 1.96123174}, {8,
      2.08510539}, {9, 2.199199703}, {10, 2.303233201}, {11,
      2.39809558}, {12, 2.484965107}, {13, 2.564965516}, {14,
      2.639061573}, {15, 2.708051263}, {16, 2.772588976}, {17,
      2.833213402}, {18, 2.890371771}, {19, 2.944438982}, {20,
      2.995732274}, {21, 3.044522438}, {22, 3.091042453}, {23,
      3.135494216}, {24, 3.17805383}, {25, 3.218875825}, {26,
      3.258096538}, {27, 3.295836866}, {28, 3.33220451}, {29,
      3.36729583}, {30, 3.401197382}}]%
  \savedata{\lowcaseboundaryhangover}[{{1.855829291, 0}, {2,
      0.1869959802}, {3, 1.06766762}}]
  \savedata{\lowcaseboundary}[{{2, 0.1869959802},
    {3, 1.06766762}, {4, 1.690982344}, {5, 2.221589919}, {6,
      2.703242291}, {7, 3.154711682}, {8, 3.586029552}, {9,
      4.003251427}, {10, 4.410314764}, {11, 4.809905402}, {12,
      5.203914923}, {13, 5.593704953}, {14, 5.980270704}, {15,
      6.364347575}, {16, 6.746483446}, {17, 7.127089173}, {18,
      7.506474627}, {19, 7.884874805}, {20, 8.262468958}, {21,
      8.639394685}, {22, 9.015758384}, {23, 9.391643001}, {24,
      9.767113812}, {25, 10.14222272}, {26, 10.51701145}, {27,
      10.89151397}, {28, 11.26575826}, {29, 11.63976768}, {30,
      12.01356194}}]%
  \savedata{\tausplitofkleft}[{{4.0, 1.645751311}, {5.0, 2.0}, {6.0, 2.0},
    {7.0, 2.0}, {8.0, 2.158883083}}]
  \savedata{\tausplitofk}[{8.0, 2.158883083}, {9.0, 2.394449155}, {10.0,
    2.605170186}, {11.0, 2.795790546}, {12.0, 2.9698133}, {13.0,
    3.129898715}, {14.0, 3.278114659}, {15.0, 3.416100402}, {16.0,
    3.545177444}, {17.0, 3.666426688}, {18.0, 3.780743516}, {19.0,
    3.888877958}, {20.0, 3.991464547}, {21.0, 4.089044875}, {22.0,
    4.182084907}, {23.0, 4.270988432}, {24.0, 4.356107661}, {25.0,
    4.43775165}, {26.0, 4.516193076}, {27.0, 4.591673732}, {28.0,
    4.66440902}, {29.0, 4.73459166}, {30.0, 4.802394763}}]%
\psset{plotstyle=curve,showpoints=true,dotsize=2pt}
\dataplot[linecolor=red,showpoints=false]{\highcaseboundaryhangover}
\dataplot[linecolor=red]{\highcaseboundary}
\dataplot[linecolor=orange,showpoints=false]{\lowcaseboundaryhangover}
\dataplot[linecolor=orange]{\lowcaseboundary}
\dataplot[dotstyle=o,linecolor=blue,plotstyle=line]{\tausplitofkleft}
\dataplot[dotstyle=o,linecolor=blue]{\tausplitofk}
\end{pspicture*}
      \caption{Case coverage and $\tausplit$}
      \label{fig:casecoverage}
    \end{figure}
    \begin{align*}
      \tausplit &:=
      \begin{cases}
        2 ~ln\frac{k}{\eulernumber{}} & \text{if $k\geq 8$},\\
        2 & \text{if $5\leq k\leq 7$},\\
        \sqrt{7}-1 & \text{if $k=4$}.
      \end{cases}
    \end{align*}
    We start with the treatment of the cases $k\geq 8$.
    
    We claim that for $\tau \leq \tausplit$ we are in the low case, ie.\
    \begin{align}
      \label{eq:low-case-cover}
      \frac{1-\eulernumber^{-\tau}}{\tau} \geq \sqrt[k]{\frac{1+\sigma}{k!}} =:
      \theta.
    \end{align}
    Consider $f_{2}(\tau) = 1-\eulernumber^{-\tau} - \theta \tau$.  It obviously
    vanishes at $\tau = 0$. The derivative of $f_{2}$ shows that there
    is exactly one maximum at $\tau = - ~ln(\theta)$, which is roughly
    $~ln k$.  To prove \ref{eq:low-case-cover} for $\tau \in \left] 0,
      \tausplit\right]$ it is thus sufficient to prove that $f_{2}$ is at
    least $0$ at the right boundary.  First, note that by Stirling's
    formula we have $k!  \geq \left( \frac{k}{\eulernumber{}}
    \right)^{k} \sqrt{2\pi}$ and so we can estimate $\theta$ by
    $\frac{\eulernumber{}}{k}$:
    \begin{align*}
      \theta &= \sqrt[k]{\frac{1+\sigma}{k!}}
      \leq \frac{\eulernumber{}}{k}
      \sqrt[k]{\smashedunderbrace{\frac{2}{\sqrt{2\pi}}}_{\leq 1}} \leq
      \frac{\eulernumber{}}{k}.
    \end{align*}
    Well, now we have
    \begin{align*}
      k \cdot f_{2}\treatparameter({\tausplit})%
      &= k-\frac{\eulernumber^{2}}{k} - k \theta \tausplit
      \\
      &\geq k-\frac{\eulernumber^{2}}{k} - 2 \eulernumber{} ~ln
      \frac{k}{\eulernumber{}} =: f_{3}(k).
    \end{align*}
    Checking that $f_{3}(\eulernumber{}) = 0$ and the derivative
    $f_{3}'(k) = (1-\frac{\eulernumber{}}{k})^{2}$ is positive for $k >
    \eulernumber{}$ shows that $f_{3}(k) > 0$ for all $k\geq 3$.  Thus
    for $\tau \leq \tausplit$ we are in the low case with the above
    choice of $\sigma$.

    It remains to check that for $\tau \geq \tausplit$ we are in the high
    case.  We use the Lagrange remainder estimate of the power series of
    the exponential function and again Stirling's formula to obtain
    \begin{align*}
      \presum{l}{0}{k-2} \frac{\tau^{l}}{l!}  &\geq \eulernumber^{\tau} \left( 1 -
        \frac{\tau^{k-1}}{(k-1)!}  \right) \geq \eulernumber^{\tau} \left( 1 -
        \left( \frac{\eulernumber{} \tau}{k-1} \right)^{k-1} \right)
      =:f_{5}(k,\tau).
    \end{align*}
    As the left hand side is \weaklyisotone{} in $\tau>0$ we consider the
    smallest $\tau$ in question: let $f_{6}(k) := f_{5}( k, \tausplit ) /
    k$.  Therein,
    \begin{align*}
      \left( \frac{\eulernumber{} \tausplit}{k-1} \right)^{k-1} &=
      ~exp\treatparameter({ -\smashedunderbrace{(k-1)}_{\text{(I)}}
        \smashedunderbrace{\left( -~ln 2 - 1 -
            ~ln~ln\frac{k}{\eulernumber{}}+~ln(k-1)
          \right)}_{\text{(II)}} }).  \smashedheights
    \end{align*}
    Observe that the term (I) is obviously positive and \weaklyisotone{} for
    $k\geq 8$.  The same is true for the term (II): its derivative
    $\frac{1}{k-1} - \frac{1}{k \left(~ln{k} - 1\right)}$ is
    positive for $k \geq \eulernumber{}^{2} \approx
    7.\rounded{38}{9056}$ and the value of term (II) at
    $\eulernumber{}^{2}$ is positive.  Using this we infer that $f_{6}$
    is \weaklyisotone{} and positive.  Checking $f_{6}(10) > 1$ ($f_{6}(10) \in
    \left]1.19,1.20\right[$) now proves $\presum{l}{0}{k-2} \frac{
      \left( 2 ~ln \frac{k}{\eulernumber{}} \right)^{l}}{l!} \geq k$ for
    $k\geq 10$.  For $k=8$ and $k=9$ we just verify this inequality
    directly.


    It remains to consider $4\leq k \leq 7$.  Here we use individual
    seperation positions as defined above:
    \begin{align*}
      \begin{array}[t]{r|rrrrrrl}
        k & 4 & 5 & 6 & 7 & (8) & (9)
        \\[1mm]\hline
        \tausplit & \sqrt{7}-1 & 2 & 2 & 2 & 2
        ~ln\frac{8}{\eulernumber{}} & 2 ~ln\frac{9}{\eulernumber{}}
        \\[0.5mm]\hline\\[-4mm]
        \frac{1-\eulernumber^{-\tausplit}}{\tausplit} &%
        0.\roundeddown{49}{0433} &%
        0.\roundeddown{43}{2332} &%
        0.\roundeddown{43}{2332} &%
        0.\roundeddown{43}{2332} &%
        0.\roundeddown{40}{9723} &%
        0.\roundeddown{37}{9534} &%
        \vbox to 0pt{\vskip1mm\hbox{$\vee$ for low case}\vss}
        \\[2mm]
        \sqrt[k]{\frac{1+\sigma}{k!}} &%
        0.\roundedup{47}{7642} &%
        0.\roundedup{39}{9239} &%
        0.\roundedup{34}{4069} &%
        0.\roundedup{30}{2845} &%
        0.\roundedup{27}{0747} &%
        0.\roundedup{24}{4982}
        \\[2mm]\hline\\[-4mm]
        \sum\limits_{0\leq l \leq k-2} \frac{\tausplit^{l}}{l!} &%
        4 &%
        6.\rounded{33}{3333} &%
        7.\rounded{00}{0000} &%
        7.\rounded{26}{6666} &%
        8.\rounded{60}{2828} &%
        10.\rounded{92}{6059} &%
        \vbox{\hbox{$\geq k$ for high case}}
      \end{array}
    \end{align*}
    Just check that for this $\tausplit$ the low and the high case
    conditions are both fulfilled.  Summing up: the claim is proved for
    $k\geq 4$.

  \item[$k<10$:] For $k=3$ we \emph{need} an explicit check as the
    estimates done in the low and the high case are too sloppy.  As the
    following actually is a general computational way to verify the
    inequality \ref{eq:right-lower-extreme} we do describe it in
    general, show computational results for $3\leq k\leq 10$ and make
    the critical case $k=3$ hand-checkable at the end.  For the
    verification we use a small trick and brute force: First, we
    substitute occurrences of $\eulernumber^{-\tau}$ with a new variable
    $\Tau$.  The task turns into showing that the bivariate polynomial
    \begin{align*}
      F_{k}(\tau,\Tau) := \presum{l}{0}{k-2} \frac{ \tau^{l} }{l!}  -
      \frac{(1-(1-\Tau)^k)}{\Tau} \left( 1 - \frac{\tau^{k-1}}{(k-1)!}
        \Tau \right)
    \end{align*}
    is non-negative at $\tau = -~ln \Tau$ for $\Tau \in
    \left]0,1\right]$.  Now, observe that $F_{k}$ is \weaklyisotone{} in
    $\tau>0$ for fixed $\Tau \in \left[0,1\right]$.  If we thus replace
    $- ~ln \Tau$ with a lower bound and we can show that the resulting
    term is still non-negative then we are done.  For $\Tau \in
    \left]0,1\right]$ we have $- ~ln \Tau \geq \presum{l}{1}{s}
    \frac{(1-\Tau)^{l}}{l}$.  This lower bound even converges to $-~ln
    \Tau$, which actually ensures that we can always find some $s$ that
    allows the following reduction.  We consider the univariate
    polynomial\endnote{Check whether $g_{k,s}(1-S)$ can be nicely
      expressed and analyzed\dots}
    \begin{align*}
      g_{k,s}(\Tau) &:= F_{k}\treatparameter({ \presum{l}{1}{s}
        \frac{(1-\Tau)^{l}}{l}, \Tau }).
    \end{align*}
    By our reasoning, the claim follows if $\forall
    \Tau\in\left]0,1\right]\colon g_{k,s}(\Tau) \geq 0$ for some $s$.
    This in turn is implied by
    \begin{align}
      \label{eq:brute-force}
      g_{k,s}(0) > 0 \quad\land\quad \text{$g_{k,s}(\Tau)$ has no zero
        for $0 < \Tau < 1$}.
    \end{align}
    The second statement can be checked using Sturm's theorem
    \citep{stu35} by only evaluating certain rational polynomials at
    $\Tau=0$ and $\Tau=1$.%
    \endnote{Let $r_{0} = g_{k,s}/T^{m}/(1-T)^{n}$ for $m$, $n$ as large
      as possible so that $r_{0}$ remains polynomial, $r_{1} = r_{0}'$,
      and recursively $r_{i+1} = -r_{i-1} \bmod r_{i}$ while
      $r_{i}\ne0$.  Assume that $r_{l}\ne0$, $r_{l+1} = 0$.  Then let
      $p_{i} = r_{i}/r_{l}$ and evaluate $p_{i}$ at $\Tau=0$ and
      $\Tau=1$.  Then the difference of the number of sign changes in
      $\left(p_{i}(0)\right)_{i}$ and $\left(p_{i}(1)\right)_{i}$ is the
      number of zeroes of $r_{0}$ in $\left]0,1\right]$.}  However, this
    only works if the chosen $s$ is large enough.  We have determined
    the smallest~$s$ that make \ref{eq:brute-force} true:
    $$
    \begin{array}{r|cccccccc}
      k & 3 & 4 & 5 & 6 & 7 & (8) & (9) \\\hline
      s & 4 & 3 & 4 & 5 & 6 & 7 & 8\\
      ~deg g_{k,s} & 11 & 13 & 21 & 31 & 43 & 57 & 73\\
      \text{time (sec)} & 0.19 & 0.29 & 0.60 & 2.8 & 24 & 235 & 1704
    \end{array}
    $$
    Though we can always divide out $(1-\Tau)^{k}$ from $g_{k,s}$ the
    degrees are in all cases quite high and the computations better done
    by a computer.  The timings refer to our own (non-optimized)
    MuPad-program\endnote{%
      \url{Brough+Csmooth-guessedlowerbound-lambdanorm-3-case_k_le_9.mn}}
    used to assert \ref{eq:brute-force}.  As $k=3$ is the only case that
    we do not cover otherwise we give $g_{3,4}$ here:
    \begin{align*}
      g_{3,4}(\Tau)/(1-\Tau)^{3} &=
      \begin{aligned}[t]
        & \frac{1}{12} (1-\Tau) + \frac{5}{6} \Tau + \Tau (1-\Tau)
        \biggl(
        \\
        &\mskip20mu \frac{481}{96} (1-\Tau)^{2} + \frac{35}{24} \Tau^{2}
        + \Tau (1-\Tau) \biggl(
        \\
        &\mskip40mu \frac{245}{96} (1-\Tau) + \frac{103}{24} \Tau + \Tau
        (1-\Tau) \biggl(
        \\
        &\mskip60mu
        \frac{119}{144} (1-\Tau)^{2} + \frac{89}{144} \Tau^{2} + \Tau
        (1-\Tau) \cdot \frac{407}{288}
        \; \biggr) \biggr) \biggr)
      \end{aligned}
    \end{align*}
    With this description you can easily see that it is positive on $
    \left[0,1\right[$.
  \end{proof}

  \endnote {\footnotesize%
    Note: $\frint{\widetilde{m}^{k}(\tau)}{\tau}{0}{\xi}$ approximation
    is wrong!\dots

    For a similar lower bound notice that for small positive $\xi$ we
    have\FOOTNOTE{Warning: In this context the formula
      obtained by partial integration is misleading!}
    \begin{align*}
      \Kapprox[~norm]{k}<\xi> &= \sum_{i\geq 0} \frac{(-\newalpha
        ~ln(B))^{i}}{(k+i)!} \xi^{k+i} \quad\in\quad \frac{\xi^{k}}{k!}
      + \mathcal{O}(\xi^{k+1}).
    \end{align*}

    Any lower bound on $\Kapprox[~norm]{k}$ has to cope with the
    $\xi^{k}$ behaviour near $0$.  Since $\Kapprox[~norm]{k}$ is \weaklyisotone{}
    we can thus say that for any positive $y < \Kapprox[~norm]{k}<1>
    \approx \frac{1}{k!}$ there is a bound $\xi_{0}>0$ such that for
    $\xi_{0}<\xi<k$ we have
    $$
    \Kapprox[~norm]{k}<\xi> \geq y.
    $$
    \todo{Find a quantitative description of the trade-off between $y$
      and $\xi_{0}$.}
    More explicitly, for $0<\xi_{0}<\xi\leq1$ we obtain
    \begin{align*}
      \Kapprox[~norm]{k}<\xi> &= \frint{B^{-\rho\newalpha}
        \widetilde{m}^{k}(\xi-\rho)}{\rho}{0}{\infty}
      \\
      &\geq \frint{B^{-\rho\newalpha}
        \frac{1}{(k-1)!}(\xi-\rho)^{k-1}}{\rho}{0}{\xi-\xi_{0}}
      \\
      &\geq \frac{\xi_{0}^{k-1} (1-B^{-\newalpha(\xi-\xi_{0})})}{(k-1)!}
      \frac{1}{\newalpha ~ln B}.
    \end{align*}
    Note that $\Kapprox[~norm]{k}$ is \weaklyisotone{} so that restricting to
    $\xi\leq 1$ does not cause any damage.  For $k>2$ setting $\xi_{0} =
    \left( 1-\frac{1}{k-1+\xi \newalpha ~ln(B)} \right) \xi$ we obtain a
    value close to the optimum:
    $$
    \Kapprox{k}<\xi> \geq x \left( \xi^{k} \frac{\newalpha^{k-1}}{~ln(B)}
      \cdot \frac{\newalpha ~ln B}{(k-1)! (k-1)} \left( 1- \frac{1}{k-1}
      \right)^{k-1} +\mathcal{O}(\xi^{k+1}) \right).
    $$
    If we thus require $\xi\geq 1-c/k$ for some fixed constant $c$ we
    can bound $\Kapprox{k}<\xi> / x$ away from $0$.  \todo{Check this!
      See \texttt{Brough+Csmooth-find\_hatlambda\_lowerbound2.mn}}

    \begin{enumerate}
    \item We can switch between the representations
      $$
      \Kapprox[~norm]{k}<\xi> = \frint{ B^{-\rho\newalpha}
        \widetilde{m}^{k}(\xi-\rho)}{\rho}{0}{\xi} = B^{-\xi\newalpha}
      \frint{ B^{\rho\newalpha} \widetilde{m}^{k}(\rho)}{\rho}{0}{\xi}
      $$
    \item To obtain a different description which is valid for all $\xi$
      we integrate by parts again:
      \begin{align*}
        \Kapprox[~norm]{k}<\xi> &= B^{-\xi\newalpha} \frint{ B^{\rho\newalpha}
          \widetilde{m}^{k}(\rho)}{\rho}{0}{\xi}
        \\
        &= \frint{\widetilde{m}^{k}(\tau)}{\tau}{0}{\xi} - \newalpha ~ln B
        \frint{ B^{-(\xi-\rho)\newalpha}
          \frint{\widetilde{m}^{k}(\tau)}{\tau}{0}{\rho} }{\rho}{0}{\xi}
        \\
        &= \left( 1 - \newalpha ~ln B \frint{ \frac{ B^{\rho\newalpha}
              \frint{\widetilde{m}^{k}(\tau)}{\tau}{0}{\rho} }{
              B^{\xi\newalpha}
              \frint{\widetilde{m}^{k}(\tau)}{\tau}{0}{\xi} }
          }{\rho}{0}{\xi} \right)
        \frint{\widetilde{m}^{k}(\tau)}{\tau}{0}{\xi}
      \end{align*}
      Estimating the antiderivative
      $\frint{\widetilde{m}^{k}(\tau)}{\tau}{0}{\rho}$ by its final
      value leads to the estimate:
      \begin{align*}
        \Kapprox[~norm]{k}<\xi> &= \left( 1 - \underbrace{ \frac{
              \frint{ B^{\rho\newalpha}
                \frint{\widetilde{m}^{k}(\tau)}{\tau}{0}{\rho}
              }{\rho}{0}{\xi} }{ \frint{ B^{\rho\newalpha} }{\rho}{0}{\xi}
              \frint{\widetilde{m}^{k}(\tau)}{\tau}{0}{\xi} }}_{=:c}
          \underbrace{ \newalpha ~ln B \frint{ B^{-(\xi-\rho)\newalpha}
            }{\rho}{0}{\xi} }_{=1-B^{-\xi\newalpha}} \right)
        \frint{\widetilde{m}^{k}(\tau)}{\tau}{0}{\xi}
        \\
        &= \left( 1 - c (1-B^{-\xi\newalpha}) \right)
        \frint{\widetilde{m}^{k}(\tau)}{\tau}{0}{\xi}
      \end{align*}
      with $0\leq c<1$. Unfortunately, $c$ still depends on $\xi$.
    \item A further option is to use part of the power series of
      $B^{-\rho\newalpha}$ or $B^{\rho\newalpha}$, respectively:
      \begin{align*}
        \Kapprox[~norm]{k}<\xi> &=
        \underbrace{\vphantom{\displaystyle\left(\presum{i}{0}{s}
              \frac{(\rho \newalpha ~ln B)^{i}}{i!}\right)} B^{-\xi\newalpha}
        }_{\text{develop!}}  \presum{i}{0}{s} \frint{ \frac{(\rho \newalpha
            ~ln B)^{i}}{i!}  \widetilde{m}^{k}(\rho) }{\rho}{0}{\xi} +
        \frint{ B^{-\xi\newalpha} \underbrace{ \left( B^{\rho\newalpha} -
              \presum{i}{0}{s} \frac{(\rho \newalpha ~ln B)^{i}}{i!}
            \right) }_{\text{small positive}} \widetilde{m}^{k}(\rho)
        }{\rho}{0}{\xi}
        \\
        &= \presum{i}{0}{s} \frint{ \frac{(-\rho \newalpha ~ln B)^{i}}{i!}
          \underbrace{\vphantom{\displaystyle\left(\presum{i}{0}{s}
                \frac{(\rho \newalpha ~ln B)^{i}}{i!}\right)}
            \widetilde{m}^{k}(\xi-\rho) }_{\text{develop?}}
        }{\rho}{0}{\xi} + \frint{ \underbrace{ \left( B^{-\rho\newalpha} -
              \presum{i}{0}{s} \frac{(-\rho \newalpha ~ln B)^{i}}{i!}
            \right) }_{\text{small}} \widetilde{m}^{k}(\xi-\rho)
        }{\rho}{0}{\xi}
      \end{align*}
    \item Another option is to shrink the interval and then use a good
      constant bound for some part there.
      \begin{center}
        \begin{pspicture}(-.5,-.5)(3.2,1.2)
          \psplot[linecolor=blue]{0}{2.7}{%
            2.71828 x -0.4 mul exp }
          \rput[lb](0.6,0.8){\scriptsize$B^{-\rho\newalpha}$}
          \psplot[linecolor=green]{0.8}{1.8}{%
            1.8 x sub 3 exp 6 div } \psplot[linecolor=green]{0}{0.8}{%
            1.8 x sub dup dup 1 -2 div mul 2 add mul -2 add mul 2 3 div
            add }
          \rput[lb](0.8,3pt){\scriptsize$\widetilde{m}^{k}(\xi-\rho)$}
          {%
            \psset{linewidth=0.4pt,tickstyle=bottom} \psaxes{->}(3,1.2)
            \rput[lt](3,0){$\rho$} \psline(1.8,0)(1.8,-3pt)
            \rput[t](1.8,-6pt){$\xi$} \psline(0.5,0)(0.5,-3pt)
            \rput[t](0.5,-6pt){$\epsilon$} }%
          \SpecialCoor
          \psline[linestyle=dotted,dotsep=0.04,linecolor=green](0.5,0)
          (!0.5 1.8 0.5 sub dup dup 1 -2 div mul 2 add mul -2 add mul 2
          3 div add) (!0.0 1.8 0.5 sub dup dup 1 -2 div mul 2 add mul -2
          add mul 2 3 div add)
        \end{pspicture}
      \end{center}
      We estimate $\widetilde{m}^{k}(\xi-\rho)$ as indicated in the
      picture:
      \begin{align*}
        \Kapprox[~norm]{k}<\xi> &= \frint{ B^{-\rho\newalpha}
          \widetilde{m}^{k}(\xi-\rho)}{\rho}{0}{\xi}
        \\
        &\geq \frint{ B^{-\rho\newalpha}}{\rho}{0}{\epsilon} \cdot
        \widetilde{m}^{k}(\xi-\epsilon)
        \\
        &= \frac{1}{\newalpha ~ln B} (1- B^{-\epsilon\newalpha})
        \widetilde{m}^{k}(\xi-\epsilon)
      \end{align*}
      This holds for any $\epsilon \in \left[0,\xi\right]$ as long as
      $0<\xi\leq \frac{k}{2}$.  So we can optimize $\epsilon$ depending
      on $\xi$.  We obtain a good value when setting $\epsilon$ by
      $$
      \frac{1}{\newalpha ~ln B} (1- B^{-\epsilon\newalpha}) = c \xi
      $$
      with some chosen $0<c<\frac{2}{k \newalpha ~ln B} <
      \frac{1}{\xi\newalpha ~ln B}$.  That leads to $\epsilon =
      -\frac{~ln(1-c\xi\newalpha~ln B)}{\newalpha ~ln B}$, which is roughly $c
      \xi$ for small $\xi$.  We obtain
      \begin{align*}
        \Kapprox[~norm]{k}<\xi> &\geq c \xi
        \widetilde{m}^{k}(\xi-\epsilon).
      \end{align*}
      Since we already know that the leading term in the power series is
      $\frac{\xi^{k}}{k!}$ we try $c = \frac{c_{1}}{k}$\dots
    \item How to formulate the trade-off?  Close to $\xi=0$, for $\delta
      \approx \epsilon^{k}$ we expect $\Kapprox[~norm]{k}<\xi> \geq
      \delta/k!$ for $\epsilon\leq \xi \leq 1$.  Since
      $\Kapprox[~norm]{k}$ is monotone, this extends to the entire range
      $\epsilon \leq \xi \leq k$.  We would however want an
      approximation that describes the shape (a smoothed Heaviside)
      better; this description would be expected to be good for
      $\epsilon \leq \xi \leq k-\epsilon$.
    \item For $0<x<k+1$ we have $~cutexp_{k}(x) \leq
      \frac{|x|^{k}}{k!}$.  [Just note that the absolute values of the
      summands decrease monotonically for $x/(k+1)<1$.]  Thus we have
      \begin{align}
        \Kapprox[~norm]{k}<\xi> &= \frac{1}{(-\newalpha ~ln B)^{k}}
        ~cutexp_{k}\treatparameter(-\xi \newalpha ~ln B)
        \nonumber\\
        \label{res:Kapprox-lowerbound6}
        &\geq \frac{\xi^{k}}{k!}  - \frac{1}{(\newalpha ~ln B)^{k}} \cdot
        \frac{( \xi \newalpha ~ln B )^{k+1}}{(k+1)!}
        \\
        &= \frac{\xi^{k}}{k!}  \left( 1 - \frac{\xi \newalpha ~ln B}{k+1}
        \right) \nonumber
      \end{align}
      for $0 < \xi < ~min\treatparameter(1,\frac{k+1}{\newalpha ~ln B})$.
      Plugging in $\xi = \frac{k}{\newalpha ~ln B}$, provided this is at
      most $1$, maximizes the right hand side and we obtain:
      \begin{align*}
        \Kapprox[~norm]{k}<\xi> &\geq \frac{k^{k}}{k! (k+1)}
        \frac{1}{(\newalpha ~ln B)^{k}}
      \end{align*}
      for any $\xi\geq \frac{k}{\newalpha ~ln B}$ (noting that
      $\Kapprox[~norm]{k}$ is \weaklyisotone{}).  A strong form of Stirling's
      formula says that
      $$
      k! = (k/e)^{k} \sqrt{2 \pi k}
      ~exp\treatparameter(\frac{1}{12k+\epsilon})
      $$
      for some $\epsilon \in \left]0,1\right[$, so we can simplify the
      bound:
      \begin{align*}
        \Kapprox[~norm]{k}<\xi> &\geq \frac{6\sqrt{2 k}}{\sqrt{\pi}(k+1)
          (12k-1)} \left(\frac{~e}{\newalpha ~ln B}\right)^{k}
      \end{align*}
      as long as $\frac{k}{\newalpha ~ln B}<1$.  Otherwise we plug in
      $\xi=1$ to \ref{res:Kapprox-lowerbound6}:
      \begin{align*}
        \Kapprox[~norm]{k}<\xi> &\geq \frac{k+1-\newalpha ~ln B}{(k+1)!}
      \end{align*}
      for $\xi\geq 1$.
    \item Let's calculate $\Kapprox[~norm]{k}<\xi>$ for $\xi\geq k$:
      \begin{align*}
        B^{\xi\newalpha} \Kapprox[~norm]{k}<\xi> &= \frint{ B^{\rho \newalpha}
          \widetilde{m}^{k}(\rho) }{\rho}{0}{k}
        \\
        &= \frint{ B^{\rho \newalpha} \frint{
            \widetilde{m}^{k-1}(\rho-\rho_{k}) }{\rho_{k}}{0}{1}
        }{\rho}{0}{k}
        \\
        &= \frint{ B^{\rho_{k} \newalpha} \frint{ B^{\rho \newalpha}
            \widetilde{m}^{k-1}(\rho) }{\rho}{-\rho_{k}}{k-\rho_{k}}
        }{\rho_{k}}{0}{1}
        \\
        &= \underbrace{ \frint{ B^{\rho_{k} \newalpha} }{\rho_{k}}{0}{1}
        }_{=\frac{B^{\newalpha}-1}{\newalpha ~ln B}} \frint{ B^{\rho \newalpha}
          \widetilde{m}^{k-1}(\rho) }{\rho}{0}{k-1}
      \end{align*}
      as long as $k\geq 2$.  For $k=1$ we obtain $ \frint{ B^{\rho
          \newalpha} \widetilde{m}^{1}(\rho) }{\rho}{0}{1} =
      \frac{B^{\newalpha}-1}{\newalpha ~ln B} $.  So by induction
      $$
      \Kapprox[~norm]{k}<\xi> = B^{-(\xi-k)\newalpha} \left(
        \frac{1-B^{-\newalpha}}{\newalpha ~ln B} \right)^{k}.
      $$
      This should be compared to
      \begin{align*}
        \frac{\Qapprox{k}<\xi>}{\newalpha^{k} B^{k+\xi\newalpha}} &\sim \left(
          \frac{C}{~ln C} - \frac{B}{~ln B} \right)^{k} \newalpha^{-k}
        B^{-k-\xi\newalpha}
        \\
        &= \left( \frac{B^{1+\newalpha}}{(1+\newalpha)~ln B} - \frac{B}{~ln B}
        \right)^{k} \newalpha^{-k} B^{-k-\xi\newalpha}
        \\
        &= \frac{B^{-(\xi-k)\newalpha}}{(1+\newalpha)^{k}} \left(
          \frac{1-(1+\newalpha)B^{-\newalpha}}{\newalpha~ln B} \right)^{k}.
      \end{align*}
      As $\Eapprox[~norm]{k}<\xi>$ is merely $\Kapprox[~norm]{k}<\xi>$
      with $\newalpha$ replaced by $\newalpha-\frac{~ln(1+\newalpha)}{~ln B}$ we
      find
      \begin{align*}
        \frac{\Eapprox{k}<\xi>}{\newalpha^{k} B^{k+\xi\newalpha}} =
        \frac{\Eapprox[~norm]{k}<\xi>}{(1+\newalpha)^{\xi}} &=
        \frac{B^{-(\xi-k)\left(\newalpha-\frac{~ln(1+\newalpha)}{~ln
                B}\right)}} {(1+\newalpha)^{\xi}} \left(
          \frac{1-B^{-\left(\newalpha-\frac{~ln(1+\newalpha)}{~ln B}\right)}}
          {\left(\newalpha-\frac{~ln(1+\newalpha)}{~ln B}\right) ~ln B}
        \right)^{k}
        \\
        &= \frac{B^{-(\xi-k)\newalpha}}{(1+\newalpha)^{k}} \left(
          \frac{1-(1+\newalpha)B^{-\newalpha}} {\newalpha ~ln B} \right)^{k}
        \left(1-\frac{~ln(1+\newalpha)}{\newalpha ~ln B}\right)^{-k}
      \end{align*}
      which looks much better.  (In our test case we have $
      \frac{~ln(1+\newalpha)}{\newalpha ~ln B} \approx 0.04319903823$, and
      $\frac{\Eapprox{k}<k>}{\Qapprox{k}<k>} \sim
      (1-\frac{~ln(1+\newalpha)}{\newalpha ~ln B})^{4} \approx 1.193200941$,
      whereas $\frac{\Kapprox{k}<k>}{\Qapprox{k}<k>} \approx
      2.320373988$.)
    \item Seeking for the behavior for $\Kapprox[~norm]{k}<\xi>$:
      \begin{enumerate}
      \item $\Kapprox[~norm]{k}<0> = 0$,
      \item $\Kapprox[~norm]{k}<\xi>$ gets positive for small $\xi>0$,
      \item $\Kapprox[~norm]{k}<k> = \left( \frac{1-B^{-\newalpha}}{\newalpha
            ~ln B} \right)^{k} > 0$,
      \item To learn more we differentiate:
        \begin{align*}
          \mathcal{D}_{\xi}\Kapprox[~norm]{k}<\xi> &= \mathcal{D}_{\xi}
          \left( B^{-\xi\newalpha} \frint{ B^{\rho\newalpha}
              \widetilde{m}^{k}(\rho)}{\rho}{0}{\xi} \right)
          \\
          &= -\newalpha ~ln B \cdot \Kapprox[~norm]{k}<\xi> +
          \widetilde{m}^{k}(\xi)
          \\
          &= \mathcal{D}_{\xi} \left( \frint{ B^{-\rho\newalpha}
              \widetilde{m}^{k}(\xi-\rho)}{\rho}{0}{\xi} \right)
          \\
          &= \frint{ B^{-\rho\newalpha}
            \mathcal{D}_{\xi}\widetilde{m}^{k}(\xi-\rho) }{\rho}{0}{\xi}
          +
          \begin{cases}
            B^{-\xi\newalpha} & \text{if $k=1$,}
            \\
            0 & \text{otherwiese.}
          \end{cases}
        \end{align*}
        Iterating these we obtain for $l<k$
        \begin{align*}
          \mathcal{D}_{\xi}^{l}\Kapprox[~norm]{k}<\xi> &= \left( -\newalpha
            ~ln B \right)^{l}\cdot \Kapprox[~norm]{k}<\xi> +
          \presum{i}{0}{l-1} (-\newalpha ~ln B)^{l-1-i}
          \mathcal{D}_{\xi}^{i} \widetilde{m}^{k}(\xi)
          \\
          &= \frint{ B^{-\rho\newalpha}
            \mathcal{D}_{\xi}^{l}\widetilde{m}^{k}(\xi-\rho)
          }{\rho}{0}{\xi} +
          \begin{cases}
            B^{-\xi\newalpha} & \text{if $k=l$,}
            \\
            0 & \text{otherwiese.}
          \end{cases}
        \end{align*}

      \item $\mathcal{D}_{\xi}\Kapprox[~norm]{k}<0> =
        \begin{cases}
          0 & \text{if $k>1$,}
          \\
          1 & \text{if $k=1$.}
        \end{cases}
        $

      \item
        \begin{align*}
          \mathcal{D}_{\xi}\Kapprox[~norm]{k}<1> &=
          \frac{~cutexp_{k}(-\newalpha~ln B)}{(-\newalpha ~ln B)^{k}}
          \\
          &= \presum{l}{k}{\infty} \frac{1}{l!} (-\newalpha~ln B)^{l-k}
          \\
          &= \frac{1}{k!}  + \frac{1}{(k+1)!} (-\newalpha~ln B) +
          \frac{1}{(k+2)!} (-\newalpha~ln B)^{2} + \dots
          \\
          &\approx \frac{1}{k!} \left(1-\frac{\newalpha ~ln B}{k+1}\right).
        \end{align*}

      \item $\mathcal{D}_{\xi}\Kapprox[~norm]{k}<k> = -\newalpha ~ln B
        \cdot \Kapprox[~norm]{k}<k> = -\newalpha ~ln B \left(
          \frac{1-B^{-\newalpha}}{\newalpha ~ln B} \right)^{k} < 0$.
        
      \item $\Kapprox[~norm]{k}$ starts \strictlyisotone{}:
        $\mathcal{D}_{\xi}^{l}\Kapprox[~norm]{k}<0> =
        \begin{cases}
          0 & \text{if $l<k$,}
          \\
          1 & \text{if $l=k$.}
        \end{cases}
        $

      \end{enumerate}
    \end{enumerate}
    \todo{Finish this\dots} %
  }%
  \ifdraft\textsuperscript{,}\fi%
  \endnote {\footnotesize \todo{Do the right hand side $[k-1,k]$!}
    \begin{enumerate}
    \item From \ref{res:order-approx} we obtain the following precise
      formula:
      \begin{align*}
        \Kapprox[~norm]{k}<\xi> &= \left( \frac{1-B^{-\newalpha}}{\newalpha
            ~ln B} \right)^{k} B^{(k-\xi)\newalpha} - \frac{ ~cutexp_{k}
          \left( (k-\xi)\newalpha ~ln B \right) }{\left(\newalpha ~ln
            B\right)^{k}}
        \\
        &= \frac{ \left( \left( 1-B^{-\newalpha} \right)^{k} -1 \right)
          B^{(k-\xi)\newalpha} + \presum{l}{0}{k-1} \frac{ ((k-\xi)\newalpha
            ~ln B)^{l} }{l!}  }{\left(\newalpha ~ln B\right)^{k}}
        \stackrel?  \geq \frac{
          \frac{ ((k-\xi)\newalpha ~ln B)^{k-1} }{(k-1)!}  }{\left(\newalpha
            ~ln B\right)^{k}}
        \\
        &= \sum_{l\geq 0} \left( \left( 1-B^{-\newalpha} \right)^{k}
          -\chi_{l\geq k} \right) \frac{ (\newalpha~ln B)^{l-k} }{l!}
        (k-\xi)^{l}
        \\
        &= \frac{ B^{(k-\xi)\newalpha} B^{-1} \presum{l}{1}{k} \binom{k}{l}
          (-1)^{i} B^{-(l-1) \newalpha} + \presum{l}{0}{k-1} \frac{
            ((k-\xi)\newalpha ~ln B)^{l} }{l!}  }{\left(\newalpha ~ln
            B\right)^{k}}
        \\
        &= \frac{ \left( \left( 1-B^{-\newalpha} \right)^{k} -1 \right)
          \left( \presum{l}{0}{k-1} \frac{ ((k-\xi)\newalpha ~ln B)^{l}
            }{l!}  + B^{(k-\xi_{~mid})\newalpha} \frac{ ((k-\xi)\newalpha ~ln
              B)^{k} }{k!}  \right) + \presum{l}{0}{k-1} \frac{
            ((k-\xi)\newalpha ~ln B)^{l} }{l!}  }{\left(\newalpha ~ln
            B\right)^{k}}
        \\
        &\stackrel?= \frac{ \left( \left( 1-B^{-\newalpha} \right)^{k} -1
          \right) \left( \presum{l}{0}{k-1} \frac{ ((k-\xi)\newalpha ~ln
              B)^{l} }{l!}  + B^{(k-\xi_{~mid})\newalpha} \frac{
              ((k-\xi_{~mid})\newalpha ~ln B)^{k-1} }{(k-1)!}  (k-\xi)
            \newalpha ~ln B \right) + \presum{l}{0}{k-1} \frac{
            ((k-\xi)\newalpha ~ln B)^{l} }{l!}  }{\left(\newalpha ~ln
            B\right)^{k}}
      \end{align*}
      for some $\xi_{~mid} \in \left]\xi,k\right[$.  Since the $1$ in $\left(
        1-B^{-\newalpha} \right)^{k} -1$ cancels the overwhelming term for
      large $B^{\newalpha}$ stems from the sum's summand for $l = k-1$.  So
      we can be sure that $\Kapprox[~norm]{k}<k-1>$ is in
      $\Omega\treatparameter(\frac{1}{\newalpha ~ln B})$.

    \item We adapt the proof of $\Kapprox[~norm]{k}<k> = \left(
        \frac{1-B^{-\newalpha}}{\newalpha ~ln B} \right)^{k}$ to produce a
      lower bound for values in $[k-1,k]$:

    \item Fix $\epsilon \in \left[k-1,k\right]$ compute
      $\Kapprox[~norm]{k}<k-1+\epsilon>$ inductively:
      \begin{align*}
        \Kapprox[~norm]{k}<k-1+\epsilon> &= \frint{
          B^{-(k-1+\epsilon-\rho) \newalpha} \frint{
            \widetilde{m}^{k-1}(\underbrace{\rho-\rho_{k}}_{=:\tau})
          }{\rho_{k}}{0}{1} }{\rho}{0}{k-1+\epsilon}
        \\
        &= \frint{ B^{-(1-\rho_{k})\newalpha} \frint{
            B^{-(k-2+\epsilon-\tau) \newalpha} \widetilde{m}^{k-1}(\tau)
          }{\tau}{-\rho_{k}}{k-2+\epsilon+(1-\rho_{k})}
        }{\rho_{k}}{0}{1}
        \\
        &\geq \underbrace{ \frint{ B^{-(1-\rho_{k})\newalpha}
          }{\rho_{k}}{0}{1} }_{=\frac{1-B^{-\newalpha}}{\newalpha ~ln B}}
        \cdot \underbrace{ \frint{ B^{-(k-2+\epsilon-\tau) \newalpha}
            \widetilde{m}^{k-1}(\tau) }{\tau}{0}{k-2+\epsilon}
        }_{=\Kapprox[~norm]{k-1}<k-2+\epsilon>}
        \\
        &\geq \left( \frac{1-B^{-\newalpha}}{\newalpha ~ln B} \right)^{k-1}
        \Kapprox[~norm]{1}<\epsilon>
        = \left( \frac{1-B^{-\newalpha}}{\newalpha ~ln B} \right)^{k-1}
        \frac{1-B^{-\epsilon\newalpha}}{\newalpha ~ln B} .
      \end{align*}
      Well that's nothing. It increases with $\epsilon$ whereas
      $\Kapprox[~norm]{k}<k-1+\epsilon>$ decreases and so the best bound
      obtained is plugging in $\epsilon=1$ which just reproduces the
      value $\Kapprox[~norm]{k}<k>$, which we already know.  But basing
      the recursion at $\Kapprox[~norm]{2}$ instead, we obtain a good
      bound:
      \begin{align*}
        \Kapprox[~norm]{k}<k-1+\epsilon> &\geq \left(
          \frac{1-B^{-\newalpha}}{\newalpha ~ln B} \right)^{k-2}
        \Kapprox[~norm]{2}<1+\epsilon>
        \\
        &= \left( \frac{1-B^{-\newalpha}}{\newalpha ~ln B} \right)^{k-2}
        \left( \frac{1-\epsilon}{\newalpha ~ln B} - \frac{1-2
            B^{-\epsilon\newalpha} + B^{-(1+\epsilon)\newalpha}}{\newalpha^{2}
            ~ln^{2} B} \right) .
      \end{align*}
      This provides a value much better than the previous though more
      complex.  We could even raise the basis to get even better bounds.
      Unfortunately, all this is still quite far from optimal: The
      difference to $\Kapprox[~norm]{k}$ is only $\mathcal{O}(k-\xi)$,
      rather than $\mathcal{O}(k-\xi)^{k}$ as the left lower bound on
      $\left[0,1\right]$.

    \item So we try to be more flexible\dots\ Instead of cutting the
      integral over $[-\rho_{k},k-2+\epsilon+(1-\rho_{k})]$ down to
      $[0,k-2+\epsilon]$, cut down to $[0,k-2+\epsilon+\delta]$ only and
      shorten $\left[0,1\right]$ to $\left[0,1-\delta\right]$ instead.
      That idea didn't work: both factors decrease with increasing
      $\delta$, so best $\delta=0$.

    \item In the first inequality we loose:
      \begin{align*}
        & \frint{ B^{-(1-\rho_{k})\newalpha} \frint{
            B^{-(k-2+\epsilon-\tau) \newalpha} \widetilde{m}^{k-1}(\tau)
          }{\tau}{k-2+\epsilon}{k-2+\epsilon+(1-\rho_{k})}
        }{\rho_{k}}{0}{1}
        \\
        &= \frint{ B^{-(1-\rho_{k})\newalpha} \frint{ B^{(1-\epsilon-\rho)
              \newalpha} \widetilde{m}^{k-1}(\rho)
          }{\rho}{\rho_{k}-\epsilon}{1-\epsilon} }{\rho_{k}}{0}{1}
      \end{align*}

    \item Switch between representations:
      \begin{align*}
        \Kapprox[~norm]{k}<\xi> &= \frint{ B^{-\rho\newalpha}
          \widetilde{m}^{k}(\xi-\rho)}{\rho}{0}{\xi}
        = \frint{ B^{-\rho\newalpha}
          \widetilde{m}^{k}(\xi-\rho)}{\rho}{0}{k} - \frint{
          B^{-\rho\newalpha} \widetilde{m}^{k}(\xi-\rho)}{\rho}{\xi}{k}
        \\
        &= B^{-\xi\newalpha} \frint{ B^{\rho\newalpha}
          \widetilde{m}^{k}(\rho)}{\rho}{0}{\xi}
        = B^{(k-\xi)\newalpha} \Kapprox[~norm]{k}<k> - B^{-\xi\newalpha}
        \frint{ B^{\rho\newalpha} \widetilde{m}^{k}(\rho)}{\rho}{\xi}{k}
        \\
        &= \left( \frac{1-B^{-\newalpha}}{\newalpha ~ln B} \right)^{k}
        B^{(k-\xi)\newalpha} - \frac{ ~cutexp_{k} \left( (k-\xi)\newalpha ~ln
            B \right) }{\left(\newalpha ~ln B\right)^{k}}
        \\
        &= \frac{ \presum{l}{0}{k-1} \frac{ \left( (k-\xi)\newalpha ~ln B
            \right)^{l}}{l!}  }{\left(\newalpha ~ln B\right)^{k}} - \frac{
          1 - \left( 1-B^{-\newalpha} \right)^{k} }{ \left( \newalpha ~ln B
          \right)^{k} } B^{(k-\xi)\newalpha}
        \\
        &= \presum{l}{0}{k-1} \frac{(1-B^{-\newalpha})^{k}}{(\newalpha ~ln
          B)^{k-l}} \frac{(k-\xi)^{l}}{l!}  - \presum{l}{k}{\infty}
        (1-(1-B^{-\newalpha})^{k})(\newalpha ~ln B)^{l-k}
        \frac{(k-\xi)^{l}}{l!}
        \stackrel?\geq \frac{(1-B^{-\newalpha})^{k}}{\newalpha ~ln B}
        \frac{(k-\xi)^{k-1}}{(k-1)!}
      \end{align*}
    \item Rewrite in terms of $\tau = (k-\xi) \newalpha ~ln B$, $A =
      B^{\newalpha} = C/B$:
      \begin{align*}
        \Kapprox[~norm]{k}<\xi> &= \frint{ A^{-\rho}
          \widetilde{m}^{k}(\tau/~ln A+\rho)}{\rho}{0}{(k-\tau/~ln A)}
        = \frint{ A^{-\rho} \widetilde{m}^{k}(\tau/~ln
          A+\rho)}{\rho}{0}{k} - \frint{ A^{-\rho}
          \widetilde{m}^{k}(\tau/~ln A+\rho)}{\rho}{k-\tau/~ln A}{k}
        \\
        &= A^{\tau/~ln A-k} \frint{ A^{\rho}
          \widetilde{m}^{k}(\rho)}{\rho}{0}{k-\tau/~ln A}
        = \eulernumber^{\tau} \Kapprox[~norm]{k}<k> - A^{-(k-\tau/~ln A)} \frint{
          A^{\rho} \widetilde{m}^{k}(\rho)}{\rho}{k-\tau/~ln A}{k}
        \\
        &= \left( \frac{1-A^{-1}}{~ln A} \right)^{k} \eulernumber^{\tau} - \frac{
          ~cutexp_{k} \left( \tau \right) }{~ln^{k} A}
        \\
        &= \frac{ \presum{l}{0}{k-1} \frac{ \tau^{l} }{l!}  }{\left(~ln
            A\right)^{k}} - \frac{ 1 - \left( 1-A^{-1} \right)^{k} }{
          \left( ~ln A \right)^{k} } \eulernumber^{\tau}
        \\
        &= \presum{l}{0}{k-1} \frac{(1-A^{-1})^{k}}{~ln^{k} A}
        \frac{\tau^{l}}{l!}  - \presum{l}{k}{\infty}
        \frac{(1-(1-A^{-1})^{k})}{~ln^{k} A} \frac{\tau^{l}}{l!}
        \stackrel?\geq \frac{(1-A^{-1})^{k}}{~ln^{k} A}
        \frac{\tau^{k-1}}{(k-1)!}
      \end{align*}
    \item Possible round-up for small $k$:
      \begin{itemize}
      \item for $k=2$ that does \emph{not} always hold:
        \begin{align*}
          1 &\geq (1-(1-\eulernumber^{-\tau})^2) \left(
            \eulernumber^{\tau} - \tau \right)
          \\
          &= (2 - \eulernumber^{-\tau}) \left( 1 - \tau
            \eulernumber^{-\tau} \right) .
        \end{align*}
      \item for $k=3$ we get
        \begin{align*}
          1 + \tau \geq (3-3\eulernumber^{-\tau}+\eulernumber^{-2\tau})
          \left( 1 - \frac{\tau^{2}}{2}\eulernumber^{-\tau} \right) .
        \end{align*}
        This is obviously true for $\tau\geq 3$
      \end{itemize}
      \todo{CONTINUE HERE.}
    \end{enumerate}
  } 
  Finally, we put together the upper bound on $\Kapprox[~norm]{k}$,
  \ref{res:Kapprox-qualitativebehavior},
  \ref{res:Kapprox-leftlowerbound}, and
  \ref{res:Kapprox-rightlowerbound} in the following theorem.
  \begin{theorem}
    \label{res:Kapprox-bounds}
    For any $k \geq 1$ we have
    $$
    \Kapprox[~norm]{k}<\xi> \leq \capprox{k} :=
    \begin{cases}
      \frac{1}{\newalpha ~ln B} & \text{if $k>0$},
      \\
      1 & \text{if $k=0$}.
    \end{cases}
    $$
    Assume $\newalpha ~ln B \geq ~max\treatparameter( ~ln 2, \frac{~ln
      16}{k} )$.  Then for any $\epsilon \in \left] 0,1 \right]$ with
    $\epsilon \leq k- \xi^{k}_{\frac12}$ or $k<3$ there is a
    $\clower{k}>0$ such that for $\xi \in \left[ \epsilon,
      k-\epsilon \right]$
    $$
    \Kapprox[~norm]{k}<\xi> \geq \frac{\clower{k}}{\newalpha ~ln B}.
    $$
    Here we can choose 
    \begin{align*}
      \clower{k} = ~min\left(
        \quad
        2^{-4}
        \cdot
        \frac{\epsilon^{k}}{k!},
        \quad
        2^{-k}
        \cdot
        \frac{\epsilon^{k-1}}{(k-1)!}
        \quad
      \right).
    \end{align*}
  \end{theorem}
  Note that $\xi^{k}_{\frac12}$ is the maximum of $\Kapprox[~norm]{k}$
  which in our experiments is always less then $k-1$ for $k\geq 3$.
  However, proving that would result in a stronger version of
  \ref{res:Kapprox-rightlowerbound}, which even in the given form
  required quite some effort.  However, we just want that to work for
  some small $\epsilon$ and that is always granted.
  \begin{proof}
    \let\oldalpha\theta
    The upper bound being proven we consider the lower bound.  First,
    keeping in mind that $\newalpha ~ln B \geq ~ln 2$ and $\newalpha ~ln B
    \geq \frac{~ln 16}{k}$, we consider the cases $k\geq 3$.  We know
    by \ref{res:Kapprox-qualitativebehavior} that $\Kapprox[~norm]{k}$
    on $[\epsilon,k-\epsilon]$ attains its minimal value at one of the
    boundaries since $\xi^{k}_{\frac12} \in [\epsilon,k-\epsilon]$ (by
    assumption).
    We thus only need to consider its values at $\xi = \epsilon$ and at
    $\xi = k-\epsilon$.

    On the left hand side \ref{res:Kapprox-leftlowerbound} gives
    $$
    \newalpha ~ln B \cdot \Kapprox[~norm]{k}<\epsilon> \geq
    ~exp\left(-\frac{~ln^2 4}{\newalpha ~ln B}\right) \frac{\epsilon^k}{k!}
    \geq ~exp\left(-\frac{~ln^2 4}{~ln 2}\right) \frac{\epsilon^k}{k!}
    = 2^{-4} \cdot \frac{\epsilon^k}{k!}
    $$
    using the conditions on $\newalpha ~ln B$.

    On the right hand side by \ref{res:Kapprox-rightlowerbound} we find
    $$
    \newalpha ~ln B \cdot \Kapprox[~norm]{k}<k-\epsilon> \geq
    (1-~exp\treatparameter(-\newalpha ~ln B))^k
    \frac{\epsilon^{k-1}}{(k-1)!}  \geq 2^{-k}
    \frac{\epsilon^{k-1}}{(k-1)!}
    $$
    using again $\newalpha ~ln B \geq ~ln 2$.  This completes the proof for
    the cases $k\geq 3$.

    We will now show corresponding lower bounds for the cases
    $k=1,2$. For $k=1$ we have for $\epsilon \in \left]0,1\right]$ using
    $\newalpha ~ln B \geq ~ln 2$:
    $$
    \newalpha ~ln B \cdot \Kapprox[~norm]{1}<\epsilon> =
    1-~exp\treatparameter(- \epsilon \newalpha ~ln B) \geq
    1-~exp\treatparameter(-\epsilon ~ln 2).
    $$
    Using \ref{res:map-1-exp-} and $\epsilon \leq 1$, we have
    $$
    1-~exp(-\epsilon ~ln 2) \geq \frac{1-~exp(- ~ln 2)}{~ln 2} \epsilon
    = \frac{1}{2 ~ln 2} \epsilon \geq \frac{1}{2} \epsilon.
    $$

    For $k=2$ we again apply \ref{res:Kapprox-leftlowerbound} for the
    left hand side as in the general case.
    \def\myc{2^{-2}}\def\oneminusmyc{\frac{3}{4}}%
    \def\myc{\left(1-\frac{1}{2 ~ln
          2}\right)}\def\oneminusmyc{\frac{1}{2 ~ln 2}}%
    For the right hand side we show that
    $$
    \newalpha^{2} ~ln^{2} B \cdot \Kapprox[~norm]{2}<2-\epsilon> \geq \myc
    \epsilon \newalpha ~ln B
    $$
    for $0 \leq \epsilon \newalpha ~ln B \leq \newalpha ~ln B$, $~ln 2 \leq
    \newalpha ~ln B$.  Since $\myc \geq 2^{-2}$ this proves the claim.
    Now, using \ref{res:order-approx} we write the left hand side minus
    the right hand side as
    $f_{5}( \epsilon \newalpha ~ln B, \newalpha ~ln B)$ with $f_{5}( \tau,
    \oldalpha ) = ~exp(\tau-2\oldalpha) - 2 ~exp(\tau-\oldalpha) + \oneminusmyc
    \tau + 1$.  We have to show that $f_{5}$ is non-negative if $0 \leq
    \tau \leq \oldalpha$ and $\oldalpha \geq ~ln 2$.  The $\oldalpha$-derivative
    of $f_{5}$,
    $$
    \frac{\partial f_{5}}{\partial \oldalpha} ( \tau, \oldalpha ) = 2
    ~exp(-\oldalpha) \left( 1 - ~exp(-\oldalpha) \right) ~exp(\tau),
    $$
    is positive for $\oldalpha > 0$.  Thus it suffices to show that $f_{5}(
    \tau, \oldalpha ) \geq 0$ for the smallest allowed $\oldalpha$, which is
    the larger of $\tau$ and $~ln 2$.  If $\tau \geq ~ln 2$ then we
    consider $f_{5}(\tau,\tau) = ~exp(-\tau) + \oneminusmyc \tau - 1$.
    This expression is increasing in this case (even for $\tau \geq ~ln
    2 + ~ln ~ln 2$) and so it is greater than or equal to $f_{5}( ~ln 2,
    ~ln 2 ) = 0$.  If otherwise $0 \leq \tau \leq ~ln 2$ then we
    consider $f_{5}(\tau,~ln 2) = -\frac{3}{4} ~exp(\tau) + \oneminusmyc
    \tau + 1$ which is \weaklyantitone{} even for $\tau\geq 0$ and so it
    is greater than or equal to $f_{5}( ~ln 2, ~ln 2 ) = 0$.%
    \endnote{%
      For $k=2$ the task is slightly more
      work. Consider now $\epsilon \in \left]0,1\right]$.  And let $f_1(z)
      := ~exp(z) - z - 1$. Then we have
      $$
      \newalpha^2 ~ln^2 B \cdot \Kapprox[~norm]{k}<\epsilon> = f_1(- \epsilon
      \newalpha ~ln B)
      $$
      We are looking for a constant $c > 0$ such that
      $$
      f_1(- \epsilon \newalpha ~ln B) \geq c \cdot \epsilon^2 \newalpha^2 ~ln^2
      B = c \cdot (- \epsilon \newalpha ~ln B)^2.
      $$
      Differently formulated, we need to find a lower bound for the
      function $f_2(z) := \frac{f_1(z)}{z^2}$. This function is isotone
      for all $z \in \R_{<0}$: Consider its derivative
      $$
      f_2'(z) = \frac{z-2}{z^3}\left( ~exp(z)-\frac{2+z}{2-z}\right).
      $$
      We have to show that for all $z\in \R_{<0}$ we have
      $$f_{21}(z) := f_2'(z) \cdot \frac{z^3}{z-2} = ~exp(z) - \frac{2+z}{2-z} \geq 0.$$
      Note that $f_{21}(0) = 0$. Thus it suffices to show that
      $$f_{21}'(z) = ~exp(z) - \frac{1}{2-z} - \frac{z+2}{(2-z)^2} < 0$$
      or equivalently that
      $$f_{22}(z) := f_{21}'(z) (2-z)^2 = ~exp(z) (2-z)^2 - 4 < 0$$
      for all $z < 0$. This, however is clear, since the derivative
      $f_{22}'(z) = z ~exp(z) (z - 2) > 0$ for all $z < 0$ and $f_{22}(0)
      = 0$. This shows that $f_2(z)$ is isotone.

      Thus for a lower bound for $f_2(z)$, we can use the minimal value
      for $z$, i.e. the value $z = -\newalpha ~ln B$.  We now find
      $$
      \Kapprox[~norm]{k}<\epsilon> \geq \frac{~exp(-\newalpha ~ln B) - 1 +
        \newalpha ~ln B}{\newalpha^2 ~ln^2 B} \cdot \epsilon^2.
      $$
      Since this expression is still dependent on the value of $\newalpha ~ln
      B$, we are going to estimate it as well.  Thus let $f_3(z) :=
      \frac{~exp(z) - 1 - z}{-z}$. We have to analyze the function
      $$f_3(-\newalpha ~ln B) \cdot \frac{\epsilon^2}{\newalpha ~ln B}.$$
      The function $f_3(z)$ is antitone for all $z \in \R$ since
      $$
      f_4(z) := -z^2 f_3'(z) = (~exp(z) (z-1) + 1) \geq 0.
      $$
      The reason is that the only zero of $f_4'(z) = z ~exp(z)$ is at $z =
      0$ and from the signchange of $f_4'(z)$ we see that we have a global
      minimum. This implies $f_4'(z) > 0$ for all $z\neq 0$.

      Thus we will obtain a lower bound for $f_3(z)$ if we use the bound
      $\newalpha ~ln B \geq ~ln 2$:
      $$
      f_3(-\newalpha ~ln B) \cdot \frac{\epsilon^2}{\newalpha ~ln B} \geq
      \frac{~exp(-\newalpha ~ln B) - 1 + \newalpha ~ln B}{\newalpha ~ln B} \cdot
      \frac{\epsilon^2}{\newalpha ~ln B} \geq \left(1 - \frac{1}{~ln
          4}\right) \cdot \frac{\epsilon^2}{\newalpha ~ln B}.
      $$
      This implies the slightly weaker lower bound given in the claim.

      We will now show the right lower bound for the case $k=2$. We again
      consider $\epsilon \in \left]0,1\right]$.  Let $f_5(z) := ~exp(z-2
      \newalpha ~ln B) - 2 ~exp(z - \newalpha ~ln B) + z + 1$.  By
      \ref{res:order-approx} we have
      $$
      \newalpha^2 ~ln^2 B \cdot \Kapprox[~norm]{k}<2-\epsilon> = f_5(\epsilon
      \newalpha ~ln B).
      $$
      Our goal is to provide a constant $c > 0$ such that
      $$
      f_5(\epsilon \newalpha ~ln B) \geq c \cdot \epsilon \newalpha ~ln B
      $$
      or equivalently that the function $f_6(z) := \frac{f_5(z)}{z}$ can
      be bounded from below by a constant. We now prove that $f_6(z)$ is
      antitone for all $z > 0$, since its derivative is
      $$f_6'(z) = \frac{~exp(z-\newalpha ~ln B) (z-1)(~exp(-\newalpha ~ln B) - 2) - 1}{z^2}.$$
      Equivalently we can show that for the function $f_7(z) := z^2
      f_6'(z)$ we have that $f_7(z) < 0$ for all $z > 0$.  Consider its
      derivative
      $$f_7'(z) = (\underbrace{~exp(-2 \newalpha ~ln B) - 2 ~exp(-\newalpha ~ln B)}_{< 0}) \underbrace{z ~exp(z)}_{> 0} < 0.$$
      Thus the function $f_7(z)$ is antitone and we have $f_7(0) = -
      (~exp(-\newalpha ~ln B) - 1)^2 < 0$.  This proves that $f_6(z)$ is
      antitone for all $z > 0$.
      We thus obtain a lower bound for $f_6(z)$ if we plug the largest
      possible value for $z$, i.e. $z = \newalpha ~ln B$.  We obtain:
      $$
      \Kapprox[~norm]{k}<2-\epsilon> \geq \frac{~exp(-\newalpha ~ln B) +
        \newalpha ~ln B - 1}{\newalpha ~ln B} \cdot \frac{\epsilon}{\newalpha ~ln
        B} = f_3(-\newalpha ~ln B) \frac{\epsilon}{\newalpha ~ln B}
      $$
      and thus
      $$
      \Kapprox[~norm]{k}<2-\epsilon> \geq \left(1-\frac{1}{~ln 4}\right)
      \frac{\epsilon}{\newalpha ~ln B} > 2^{-2} \frac{\epsilon}{\newalpha ~ln
        B}.
      $$
    } %
  \end{proof}
  Summing up we obtain:
  \begin{corollary}
    For any $\epsilon \in \left]0,1\right]$ and any $k\geq 1$ we have
    uniformly for $\xi \in \left[\epsilon,k-\epsilon\right]$
    $$
    \Kapprox[~norm]{k}<\xi> \in
    \Theta\treatparameter(\frac{1}{\newalpha ~ln B}).
    \qed
    $$
  \end{corollary}

  \section{Estimating the estimate $\Kerror{k}$}
  \label{sec:estimate-Kerror}

  The recurrence \ref{def:Lerror} for $\Kerror{k}$ is more complex than
  the one for $\Kapprox{k}$, so instead of solving it we estimate
  it.  
  We consider also here the normed version $\Kerror[~norm]{k}<\xi> :=
  \frac{\Kerror{k}<\xi>}{\newalpha^{k}B^{k+\xi\newalpha}}$.  To
  better 
  understand how the error behaves we compute it for~$k=1$:
  \begin{align*}
    \Kerror[~norm]{1}<\xi> &=
    \begin{cases}
      0
      &
      \text{if $\xi\in \left]- \infty, 0\right[$},
      \\
      \frac{1+\xi \newalpha}{8 \pi \newalpha}
      \frac{~ln B}{B^{ \frac{1}{2} +  \frac{\xi \newalpha}{2}}}
      +
      \frac{1}{8 \pi \newalpha}
      \frac{~ln B}{B^{ \frac{1}{2} + \xi \newalpha}}
      &
      \text{if $\xi\in \left[0, 1\right[$},
      \\
      \frac{1+\newalpha}{8 \pi \newalpha}
      \frac{~ln B}{B^{\frac{1}{2}-\frac{\newalpha}{2} + \xi \newalpha}}
      +
      \frac{1}{8 \pi \newalpha}
      \frac{~ln B}{B^{\frac{1}{2} + \xi \newalpha}}
      &
      \text{if $\xi\in \left[1, \infty\right[$}.
      \\
    \end{cases}
  \end{align*}
  From this we estimate $\Kerror[~norm]{1}$ directly:
  \begin{align*}
    \Kerror[~norm]{1}<\xi> &\leq
    \begin{cases}
      0 & \text{if $\xi\in \left]- \infty, 0\right[$},
      \\
      \frac{(2+\newalpha) ~ln B}{8\pi \newalpha} B^{-\frac{1+\xi\newalpha}{2}} &
      \text{if $\xi\in \left[0, 1\right[$},
      \\
      \frac{(1+\newalpha) ~ln B}{8\pi \newalpha} 
      {B^{-\frac{1}{2}+\frac{\newalpha}{2}-\xi\newalpha}} +
      \frac{~ln B}{8\pi\newalpha} \cdot {B^{-\frac{1}{2}-\xi\newalpha}} &
      \text{if $\xi\in \left[1, \infty\right[$}.
      \\
    \end{cases}
  \end{align*}
  We have also looked at precise expressions for larger $k$, yet they
  are huge and do not give rise to better bounds.%
  \endnote{%
    \label{en:lambdahat}
    For $k=2$ we obtain the following tremendous expression
    \begin{align*}
      \Kerror{2}<\xi>/x &=
      \begin{cases}
        0
        &
        \text{if $\xi\in \left]- \infty, 0\right[$},
        \\
        \begin{array}{@{}r@{}}
          \left(\frac{3}{4 \pi}\right) \frac{B^{- \frac{1}{2}}}{1}
          \\
          +
          \frac{1}{\pi} \frac{B^{- \frac{1}{2}}}{~ln B}
          \\
          +
          \left( - \frac{1}{2 \pi} - \frac{\xi \newalpha}{2 \pi}\right) \frac{B^{ - \frac{\xi \newalpha}{2} - \frac{1}{2}}}{1}
          \\
          +
          \left(-\frac{1}{\pi}\right) \frac{B^{ - \frac{\xi \newalpha}{2} - \frac{1}{2}}}{~ln B}
          \\
          +
          \left(\frac{{\newalpha}^{3} {\xi}^{3}}{768 {\pi}^{2}} + \frac{{\newalpha}^{2} {\xi}^{2}}{128 {\pi}^{2}} + \frac{\xi \newalpha}{128 {\pi}^{2}}\right) \frac{B^{ - \frac{\xi \newalpha}{2} - 1}}{~ln^{-3} B}
          \\
          +
          \left(\frac{{\newalpha}^{2} {\xi}^{2}}{128 {\pi}^{2}} + \frac{\xi \newalpha}{16 {\pi}^{2}} + \frac{3}{64 {\pi}^{2}}\right) \frac{B^{ - \frac{\xi \newalpha}{2} - 1}}{~ln^{-2} B}
          \\
          +
          \left(-\frac{1}{4 \pi}\right) \frac{B^{ - \xi \newalpha - \frac{1}{2}}}{1}
          \\
          +
          \left(\frac{1}{64 {\pi}^{2}}\right) \frac{B^{ - \xi \newalpha - 1}}{~ln^{-2} B}
        \end{array}
        &
        \text{if $\xi\in \left[0, 1\right[$},
        \\
        \begin{array}{@{}r@{}}
          \left( - \frac{1}{4 \pi} - \frac{\newalpha}{4 \pi}\right) \frac{B^{ - \frac{\newalpha}{2} - \frac{1}{2}}}{1}
          \\
          +
          \left(-\frac{1}{\pi}\right) \frac{B^{ - \frac{\newalpha}{2} - \frac{1}{2}}}{~ln B}
          \\
          +
          \left(\frac{3}{4 \pi} - \frac{\newalpha}{2 \pi} - \frac{\xi}{4 \pi} + \frac{3 \xi \newalpha}{4 \pi}\right) \frac{B^{\frac{\newalpha}{2} - \frac{\xi \newalpha}{2} - \frac{1}{2}}}{1}
          \\
          +
          \frac{1}{\pi} \frac{B^{\frac{\newalpha}{2} - \frac{\xi \newalpha}{2} - \frac{1}{2}}}{~ln B}
          \\
          +
          \left(\frac{\xi}{4 \pi} - \frac{1}{4 \pi} - \frac{\xi \newalpha}{4 \pi}\right) \frac{B^{ - \frac{\newalpha}{2} - \frac{\xi \newalpha}{2} - \frac{1}{2}}}{1}
          \\
          +
          \left(\frac{1}{4 \pi}\right) \frac{B^{\newalpha - \xi \newalpha - \frac{1}{2}}}{1}
          \\
          +
          \left(\frac{{\newalpha}^{3} \xi}{128 {\pi}^{2}} - \frac{{\newalpha}^{3} {\xi}^{3}}{768 {\pi}^{2}} - \frac{{\newalpha}^{3}}{192 {\pi}^{2}} - \frac{{\newalpha}^{2} {\xi}^{2}}{128 {\pi}^{2}} + \frac{{\newalpha}^{2} \xi}{64 {\pi}^{2}} - \frac{\xi \newalpha}{128 {\pi}^{2}} + \frac{\newalpha}{64 {\pi}^{2}}\right) \frac{B^{ - \frac{\xi \newalpha}{2} - 1}}{~ln^{-3} B}
          \\
          +
          \left(\frac{{\newalpha}^{2} \xi}{16 {\pi}^{2}} - \frac{{\newalpha}^{2} {\xi}^{2}}{128 {\pi}^{2}} - \frac{{\newalpha}^{2}}{64 {\pi}^{2}} + \frac{\newalpha}{16 {\pi}^{2}} - \frac{\xi}{32 {\pi}^{2}} + \frac{3}{64 {\pi}^{2}}\right) \frac{B^{ - \frac{\xi \newalpha}{2} - 1}}{~ln^{-2} B}
          \\
          +
          \left( - \frac{1}{4 \pi} - \frac{\newalpha}{4 \pi}\right) \frac{B^{\frac{\newalpha}{2} - \xi \newalpha - \frac{1}{2}}}{1}
          \\
          +
          \left(\frac{1}{32 {\pi}^{2}} - \frac{\xi}{32 {\pi}^{2}} + \frac{\xi \newalpha}{32 {\pi}^{2}}\right) \frac{B^{ - \frac{\newalpha}{2} - \frac{\xi \newalpha}{2} - 1}}{~ln^{-2} B}
          \\
          +
          \left(-\frac{1}{4 \pi}\right) \frac{B^{ - \xi \newalpha - \frac{1}{2}}}{1}
          \\
          +
          \left(\frac{1}{32 {\pi}^{2}} + \frac{\newalpha}{32 {\pi}^{2}}\right) \frac{B^{\frac{\newalpha}{2} - \xi \newalpha - 1}}{~ln^{-2} B}
          \\
          +
          \left(\frac{1}{64 {\pi}^{2}}\right) \frac{B^{ - \xi \newalpha - 1}}{~ln^{-2} B}
        \end{array}
        &
        \text{if $\xi\in \left[1, 2\right[$},
        \\
        \begin{array}{@{}r@{}}
          \left(\frac{1}{4 \pi} + \frac{\newalpha}{4 \pi}\right) \frac{B^{\frac{3 \newalpha}{2} - \xi \newalpha - \frac{1}{2}}}{1}
          \\
          +
          \left(\frac{1}{4 \pi}\right) \frac{B^{\newalpha - \xi \newalpha - \frac{1}{2}}}{1}
          \\
          +
          \left( - \frac{1}{4 \pi} - \frac{\newalpha}{4 \pi}\right) \frac{B^{\frac{\newalpha}{2} - \xi \newalpha - \frac{1}{2}}}{1}
          \\
          +
          \left(\frac{{\newalpha}^{2}}{64 {\pi}^{2}} + \frac{\newalpha}{32 {\pi}^{2}} + \frac{1}{64 {\pi}^{2}}\right) \frac{B^{\newalpha - \xi \newalpha - 1}}{~ln^{-2} B}
          \\
          +
          \left(-\frac{1}{4 \pi}\right) \frac{B^{ - \xi \newalpha - \frac{1}{2}}}{1}
          \\
          +
          \left(\frac{1}{32 {\pi}^{2}} + \frac{\newalpha}{32 {\pi}^{2}}\right) \frac{B^{\frac{\newalpha}{2} - \xi \newalpha - 1}}{~ln^{-2} B}
          \\
          +
          \left(\frac{1}{64 {\pi}^{2}}\right) \frac{B^{ - \xi \newalpha - 1}}{~ln^{-2} B}
        \end{array}
        &
        \text{if $\xi\in \left[2, \infty\right[$},
        \\
      \end{cases}
    \end{align*}
  }
  \begin{theorem}\label{res:Kerror-upperbound}
    Define values $\cerror{k}$ recursively by
    $$
    \cerror{k}
    := 
    \cerror{k-1}
    +
    \frac{4+3~ln B}{8\pi\newalpha\sqrt{B}}
    \cboth{k-1}
    $$
    for $k\geq 3$ based on $\cerror{0} := 0$, $\cerror{1} :=
    \frac{(2+\newalpha) ~ln B}{8\pi\newalpha\sqrt{B}}$, and
    \begin{align*}
      \cerror{2} &:= 
      \frac{6+3\newalpha}{8\pi\newalpha^{2}}
      \frac{1}{\sqrt{B}}
      +
      \frac{4+3~ln B}{8\pi\newalpha\sqrt{B}}
      \cboth{1}
      \\
      &=
      \frac{9+3\newalpha}{8\pi\newalpha^{2}}
      \frac{1}{\sqrt{B}}
      + \frac{1}{2\pi\newalpha^{2}\sqrt{B} ~ln B}
      +
      \frac{(2+\newalpha) (4+3~ln B) ~ln B}{64\pi^{2}\newalpha^{2}B}
      .
    \end{align*}
    Then for any $k$ and $\xi \in \R$ we have
    $$
    \Kerror[~norm]{k}<\xi> \leq \cerror{k}
    $$
    If $\alpha \geq \frac{~ln B}{\sqrt{B}}$ we have for $k\geq 2$ and
    large $B$ the inequality
    $$
    \cerror{k} \leq \frac{(2^k-1) (1+\newalpha)}{\newalpha^{2} \sqrt{B}}.
    $$
    For $k=1$ the order of $\cerror{1}$ is necessarily slightly
    larger. More precisely, we have for large $B$ that
    $$\cerror{1} \leq \frac{(1+\alpha) ~ln B}{\newalpha \sqrt{B}}.$$
  \end{theorem}
  Instead of defining $\cerror{2}$ and $\cerror{1}$ we could have left
  that to the recursion.  But the given values are smaller than the ones
  derived from the recursion based on $\cerror{0}$ only.  For $k=1$ the
  recursion would give $\frac{4+3~ln B}{8\pi\newalpha\sqrt{B}}$.  For
  $k\geq 2$ however the improvement due to these explicit settings is a
  factor of order $~ln B$.

  \begin{proof}
    We first show that the value $\cerror{k}$ is bounded as claimed.
    For $k = 1$ the claim follows directly from the definition, since
    we have for $B > ~exp(1)$ that $(2+\newalpha) ~ln B \leq 2
    (1+\newalpha) ~ln B$ as $\newalpha$ is positive. For $k=2$ we have
    for $~ln^2 B \leq\sqrt{B}$ the inequality
    \begin{align*}
      \cerror{2} &= \frac{6+3\newalpha}{8\pi\newalpha^{2}}
      \frac{1}{\sqrt{B}} + \frac{4+3~ln B}{8\pi\newalpha\sqrt{B}}
      \cboth{1}\\
      &\leq \frac{1+\newalpha}{\newalpha^{2} \sqrt{B}} +
      \frac{1}{\newalpha^2 \sqrt{B}} + \frac{(1+\newalpha) ~ln^2 B}{\newalpha^2 B}\\
      &\leq \frac{3 (1+\newalpha)}{\newalpha^2 \sqrt{B}}.
    \end{align*}
    For $k \geq 2$ we proceed inductively. We have
    \begin{align*}
      \cerror{k} &= \cerror{k-1} + \frac{4+3~ln B}{8\pi\newalpha\sqrt{B}} \cboth{k-1}\\
      &\leq \frac{(2^{k-1}-1) (1+\newalpha)}{\newalpha^2 \sqrt{B}} + \frac{~ln B}{\newalpha \sqrt{B}} \left( \frac{1}{\newalpha ~ln B} + \frac{(2^{k-1}-1) (1+\newalpha)}{\newalpha ~ln B} \right)\\
      &= \frac{(2^k-1) (1+\newalpha)}{\newalpha^2 \sqrt{B}}.
    \end{align*}

    Now we prove the remaining estimate by induction on $k$.  The case
    $k=0$ is true by definition of $\Kerror{0}$ (with equality).  For
    $k=1$ the inspection above proves the claim.  The explicit
    calculation of $\Kerror{1}$ also shows that a bound of order
    $\mathcal{O} \left( \Frac{x}{\sqrt{B}} \right)$ is impossible.  We
    defer the case $k=2$ to the end of the proof as most of it will be
    as in the general case.  So assume $k\geq 3$.  Using the
    definition of $\Kerror{k}$ from \ref{def:Lerror} we split
    $\Kerror{k}$ into three summands:
    \begin{align*}
      \Kerror[~norm]{k}<\xi>&=
      \begin{aligned}[t]
        &\frint{ \Kerror[~norm]{k-1}<\xi-\rho> }{\rho}{0}{1}
        \\
        &+ \frac{2 \widehat{E}(B)}{\newalpha B} \cdot
        \Kboth[~norm]{k-1}<\xi>
        \\
        &+ \frint{ \Kboth[~norm]{k-1}<\xi-\rho>
          \widehat{E}'(B^{1+\rho\newalpha}) ~ln B}{\rho}{0}{1}.
      \end{aligned}
    \end{align*}
    For $\widehat{E}(x) = \frac{1}{8\pi} \sqrt{x} ~ln x$ we calculate
    as a preparative
    \begin{align}
      \label{res:error-estimate-B-schoenfeld}
      \frac{2\widehat{E}(B)}{\newalpha B} &= \frac{~ln B}{4\pi
        \newalpha \sqrt{B}},
      \\
      \label{res:error-estimate-C-schoenfeld}
      \frint{ \widehat{E}'(B^{1+\rho\newalpha}) ~ln B}{\rho}{0}{1} &=
      \frac{4+~ln B}{8\pi\newalpha\sqrt{B}} - \frac{4+~ln
        C}{8\pi\newalpha\sqrt{C}}.
    \end{align}
    The first summand of $\Kerror[~norm]{k}$ is at most $\cerror{k-1}$
    by induction hypothesis.
    The second summand we estimate using
    \ref{res:error-estimate-B-schoenfeld} by
    $$
    \frac{~ln B}{4 \pi \newalpha \sqrt{B}} \cboth{k-1}.
    $$
    The third summand is bounded by
    $$
    \cboth{k-1} \frint{ \widehat{E}'(B^{1+\rho\newalpha}) ~ln B
    }{\rho}{0}{1}.
    $$
    By \ref{res:error-estimate-C-schoenfeld} the third summand is at
    most
    $$
    \frac{4+~ln B}{8\pi\newalpha\sqrt{B}} \cboth{k-1}.
    $$
    This completes the proof of the case $k\geq 3$:
    \begin{align*}
      \Kerror[~norm]{k}<\xi> &\leq
      \cerror{k-1}
      + \frac{~ln B}{4 \pi \newalpha \sqrt{B}} \cboth{k-1}
      + \frac{4+~ln B}{8\pi\newalpha\sqrt{B}} \cboth{k-1}
      \\
      &= \cerror{k}.
    \end{align*}

    For $k=2$ we have to improve the estimate of the first summand
    only, since the other terms anyways are of order at most
    $\mathcal{O}\left( \Frac{1}{\sqrt{B}} \right)$.  For $\xi<0$ there
    is nothing to prove.  For $\xi \in \left[0,1\right[$ we find for
    this first summand
    \begin{align*}
      \frint{\Kerror[~norm]{1}<\xi-\rho> }{\rho}{0}{\xi} &\leq
      \frint{ \frac{(2+\newalpha)~ln B}{8\pi \newalpha }
        B^{-\frac{1}{2}-\frac{(\xi-\rho)\newalpha}{2}} }{\rho}{0}{\xi}
      \\
      &= \frac{(2+\newalpha)~ln B}{8\pi\newalpha} B^{-\frac{1}{2}}
      \underbrace{ B^{-\frac{\xi\newalpha}{2}} \frint{
          B^{\frac{\rho\newalpha}{2}} }{\rho}{0}{\xi} }_{=
        \frac{2}{\newalpha~ln B} \left( 1 -
          B^{-\frac{\xi\newalpha}{2}} \right)}
      \\
      &\leq \frac{2+\newalpha}{4\pi\newalpha^{2}\sqrt{B}}.
    \end{align*}
    For $\xi\in \left[1,2\right]$ we find
    \begin{align*}
      \frint{\Kerror[~norm]{1}<\xi-\rho>}{\rho}{0}{1} &\leq
      \begin{aligned}[t]
        &\frint{ \frac{(2+\newalpha)~ln B}{8\pi\newalpha}
          B^{-\frac{1}{2}-\frac{(\xi-\rho)\newalpha}{2}}}{\rho}{\xi-1}{1}
        \\
        &+ \frint{ 
          { \Kerror[~norm]{1}<1> B^{1+\newalpha}
            B^{-(1+(\xi-\rho)\newalpha)}
          }
        }{\rho}{0}{\xi-1}
      \end{aligned}
      \\\allowpagebreak &=
      \begin{aligned}[t]
        & \frac{(2+\newalpha)~ln B}{8\pi\newalpha}
        B^{-\frac{1}{2}-\frac{(\xi-1)\newalpha}{2}} \underbrace{
          B^{-\frac{\newalpha}{2}} \frint{ B^{\frac{\rho\newalpha}{2}}
          }{\rho}{\xi-1}{1} }_{= \frac{2}{\newalpha~ln B} \left( 1 -
            B^{-\frac{(2-\xi)\newalpha}{2}} \right)}
        \\
        &+ \Kerror[~norm]{1}<1> \underbrace{ B^{-(\xi-1)\newalpha}
          \frint{B^{\rho\newalpha}}{\rho}{0}{\xi-1}
        }_{=\frac{1}{\newalpha~ln B} \left( 1 - B^{-(\xi-1)\newalpha}
          \right)}
      \end{aligned}
      \\\allowpagebreak &\leq \frac{2+\newalpha}{4\pi\newalpha^{2}}
      B^{-\frac{1}{2}-\frac{(\xi-1)\newalpha}{2}}
      +
      \frac{1+\newalpha}{8 \pi \newalpha^{2}}
      B^{-\frac{1}{2}-\frac{\newalpha}{2}} + \frac{1}{8 \pi
        \newalpha^{2}} B^{-\frac{1}{2} - \newalpha}
      \\\allowpagebreak &\leq
      \frac{6+3\newalpha}{8\pi\newalpha^{2}\sqrt{B}}
    \end{align*}
    As $\Kerror[~norm]{1}$ decreases for $\xi\geq 1$ this bound also
    holds for $\xi\geq 2$.  Putting everything together the above
    defined value $\cerror{2}$ bounds $\Kerror[~norm]{2}<\xi>$ as
    claimed.
  \end{proof}

  It is tempting to guess that we can save more $~ln B$ factors for
  larger $k$.  However, inspecting $\Kerror[~norm]{2}$ shows that, say,
  $\Kerror[~norm]{3}<\frac12> \in
  \Omega\treatparameter(\frac{1}{\sqrt{B}})$.
  (For $\xi\in \left[ 0,1 \right[$ we find $\Kerror[~norm]{2}<\xi> =
  \frac{3}{4\pi\newalpha \sqrt{B}} +
  \mathcal{O}\treatparameter(\frac{1}{\sqrt{B} ~ln B})$.)

  \section{Reestimating $\Kerror{k}$ without Riemann}

  If you do not want to assume the Riemann hypothesis then only weaker
  bounds $\widehat{E}(x)$ on $|\pi(x)-~Li(x)|$ can be used.  In
  \cite{for02,for02a} we found the following explicit bounds, the first
  one he attributes to a paper by Y. Cheng which we could not find.
  \begin{fact}
    \label{res:pnt-ford}
    \begin{itemize}
    \item For $x>10$ we have
      $$
      \left| \pi(x) - ~Li(x) \right| \leq 11.88
      \, x (~ln x)^{\frac{3}{5}} 
      ~exp \left( -\frac{1}{57}
        (~ln x)^{\frac{3}{5}}
        (~ln ~ln x)^{-\frac{1}{5}}
      \right).
      $$
    \item There is a constant $C$ and a frontier $x_{0}$ such that for
      $x>x_{0}$ we have
      $$
      \left| \pi(x) - ~Li(x) \right| \leq
      C \, x 
      ~exp \left( -0.2098
        (~ln x)^{\frac{3}{5}}
        (~ln ~ln x)^{-\frac{1}{5}}
      \right).
      $$
    \end{itemize}
  \end{fact}
  Admittedly, these bounds only start to be meaningful at large values
  of $x$ (eg.\ the first statement around
  $10^{\spacednumber{159299}}$).  All those bounds are of the form:
  For all $x>x_{0}$
  $$
  \left| \pi(x) - ~Li(x) \right| \leq 
  \underbrace{
    C
    \, x \, (~ln x)^{c_{0}} 
    ~exp \left(
      -A \,
      (~ln x)^{c_{1}}
      (~ln ~ln x)^{-c_{2}}
    \right)
  }_{=: \widehat{E}(x)}
  $$
  holds.  Here, $C>0$, $x_{0}>0$, $c_{0}\in\R$, $c_{1}>0$, $c_{2}>0$ and
  $A>0$ are given parameters (which are not always known).
  Note that we have that
  $$
  \Frac{\widehat{E}(x)}{x}
  $$
  is \weaklyantitone{} for large $x$.  Actually, with the parameter sets
  from \ref{res:pnt-ford} this is already true for $x\geq 5$.  Moreover,
  the quotient of the relative errors at~$x^{1+\newalpha}$ and at~$x$
  $$
  \frac{
    \widehat{E}(x^{1+\newalpha}) \frac{~ln x^{1+\newalpha}}{x^{1+\newalpha}}
  }{
    \widehat{E}(x) \frac{~ln x}{x}
  }
  =
  \frac{
    (1+\newalpha) \widehat{E}(x^{1+\newalpha})
  }{
    x^{\newalpha} \widehat{E}(x)
  }
  $$
  is bounded (or even tends to zero) with $x\to\infty$ for any
  $\newalpha>0$.  This follows from $\Frac{\widehat{E}(x)}{x}$ decreasing
  when $\newalpha$ is constant, but you may also consider values for
  $\newalpha$ that increase when $x$ grows.

  Revisiting the proof of \ref{res:Kerror-upperbound} shows that only
  \ref{res:error-estimate-B-schoenfeld},
  \ref{res:error-estimate-C-schoenfeld}, and the initial values
  $\cerror{1}$ and $\cerror{2}$ depend on the specific bound
  $\widehat{E}$.  We now use the following recursion for the bounds:
  $$
  \cerror{k}
  := 
  \cerror{k-1}
  +
  \underbrace{\left(
      \frac{2 \widehat{E}(B)}{\newalpha B}
      + \frint{ \widehat{E}'(B^{1+\rho\newalpha}) ~ln B }{\rho}{0}{1}
    \right)}_{=:u}
  \cboth{k-1}
  $$
  for $k\geq 1$ based on $\cerror{0} = 0$, and possibly values for
  $\cerror{1}$ and $\cerror{2}$.

  To bound $u$ tightly the trickiest step is bounding the integral.  As
  our interests lie elsewhere we take the easy way out.  We integrate by
  parts and use that $\Frac{\widehat{E}(x)}{x}$ is \weaklyantitone{} for the
  following rough estimate
  \begin{align*}
    \frint{ \widehat{E}'(B^{1+\rho\newalpha}) ~ln B }{\rho}{0}{1}
    &=
    \frac{\widehat{E}(B^{1+\newalpha})}{\newalpha B^{1+\newalpha}}
    -
    \frac{\widehat{E}(B)}{\newalpha B}
    +
    ~ln B
    \underbrace{
      \frint{\frac{\widehat{E}(B^{1+\rho\newalpha})}{B^{1+\rho\newalpha}}}{\rho}{0}{1}
    }_{\leq \frac{\widehat{E}(B)}{B}}
    .
  \end{align*}
  Thus $u$ is bounded by
  \begin{align*}
    u 
    &\leq
    \left(
      1 + \frac{1}{\newalpha~ln B} 
      \left(
        1+ 
        \smashedunderbrace{
          \frac{\widehat{E}(B^{1+\newalpha})}{B^{\newalpha} \widehat{E}(B)}
        }_{\text{bounded}}
      \right)
    \right)
    \frac{\widehat{E}(B) ~ln B}{B}.
    \smashedheights
  \end{align*}
  In the following we neglect the bounded term, as we can compensate its
  effect for example by a small additional factor.  Since $u$ is small
  for large $B$, we expect $\cerror{k}$ to be dominated by $\cerror{1} =
  u$.
  Precisely, for $k>1$ we have $\cerror{k} = (1+u) \cerror{k-1} +
  \frac{u}{\newalpha~ln B}$, thus
  \begin{align*}
    \cerror{k}
    &= (1+u)^{k-1} u
    +
    \frac{(1+u)^{k-1}-1}{\newalpha~ln B}
    \sim
    \left( 1+\frac{k-1}{\newalpha ~ln B} \right)
    u 
    .
  \end{align*}

  \begin{theorem}
    \label{res:order-error-woRiemann}
    Assume that $\widehat{E}(x)$ bounds $ \left| \pi(x) - ~Li(x)
    \right|$ and $\widehat{E}(x)/x$ is \weaklyantitone{} for $x>x_{0}$
    and the relative error decreases fast, ie.
    $$
    \frac{
      \widehat{E}(x^{1+\newalpha}) \frac{~ln x^{1+\newalpha}}{x^{1+\newalpha}}
    }{
      \widehat{E}(x) \frac{~ln x}{x}
    }
    =
    \frac{
      (1+\newalpha) \widehat{E}(x^{1+\newalpha})
    }{
      x^{\newalpha} \widehat{E}(x)
    }
    $$
    is bounded under the chosen behavior of $\newalpha$.  Then for any
    $k\geq 2$ and $B$ large we have
    $$
    \Kerror[~norm]{k}<\xi> \in
    \mathcal{O}\treatparameter({
      \left( 1+\frac{k-1}{\newalpha ~ln B} \right)
      \left( 1+\frac{1}{\newalpha ~ln B} \right)
      \frac{\widehat{E}(B) ~ln B}{B}
    })
    $$
    for $k\geq 2$.
    \qed
  \end{theorem}
  This is close to optimal, we only loose a factor of order $~ln B$ in
  the relative error compared to the used error bound in the prime
  number theorem:
  $$
  \frac{\Kerror[~norm]{k}<\xi>}{\Kapprox[~norm]{k}<\xi>}
  \leq
  \frac{
    \left( 1+\frac{k-1}{\newalpha ~ln B} \right)
    \left( 1+\frac{1}{\newalpha ~ln B} \right)
    \frac{\widehat{E}(B) ~ln B}{B}
  }{
    \frac{\clower{k}}{\newalpha ~ln B}
  }
  \sim
  \frac{\newalpha}{\clower{k}} ~ln B\cdot
  \frac{\widehat{E}(B) ~ln B}{B}.
  $$
  The assumptions on $\widehat{E}$ also hold for explicit error bounds
  with $\widehat{E}(x) \in \mathcal{O}\treatparameter( \frac{x}{~ln^{l}
    x} )$.  Due to the lost $~ln B$ the result is only meaningful if
  $l\geq 3$, so \cite{rossch62} does not suffice.  From \cite{dus98} we
  can use $\widehat{E}(x) = 2.3854 \frac{x}{~ln^{3}x}$ for
  $x>\spacednumber{355991}$, and obtain
  $$
  \Kerror[~norm]{k}<\xi> \in \mathcal{O}\treatparameter(\frac{1}{~ln B}).
  $$

  \section{Improvements}
  \label{sec:improvements}

  We have a look at the quality of \ref{res:order-basis} when applied to
  our inspiring application.  There $B=1100\cdot10^{6}$, $C=2^{37}-1$,
  $\newalpha = ~ln(C)/~ln(B)-1 = 0.\rounded{231}{902}$, and the
  largest $k$ of interest is $k=4$.  For these parameters we find
  $$
  \left[ \frac{1}{(1+\newalpha)^{k}}, 1 \right]
  \subset
  \left[ 0.434, 1 \right].
  $$
  That is a great loss when we try to enclose the function of interest
  in a small interval.  Actually, we can improve the theorem for the
  price of a slightly more complicated recursion.  The present result
  was based on approximating $\mathring{\kappa}(\rho) = \frac{1}{~ln p}$
  for $p = B^{1+\rho\newalpha} \in \left]B,C\right]$ by
  $\mathring{\lambda}(\rho) =\frac{1}{~ln B}$ in the recursion of
  \ref{def:Qapprox}:
  $$
  \Qapprox{k}<\xi> = \newalpha ~ln B \frint{\Qapprox{k-1}<\xi-\rho> \cdot
    \mathring{\kappa}(\rho) B^{1+\rho\newalpha}}{\rho}{B}{C}.
  $$
  We get a better bound by using
  $$
  \mathring{\nu}(\rho) = \frac{1-\rho}{~ln B} + \frac{\rho}{~ln C}
  $$
  with $p=B^{1+\rho\newalpha}$ instead.
  \begin{definition}
    \label{def:Napprox}
    For $x\geq 0$ we let $\Napprox{0} := \Q{0}$, and
    recursively for $k>0$
    \begin{align*}
      \Napprox{k}<\xi> &:= 
      \newalpha ~ln B
      \frint{\Napprox{k-1}<\xi-\rho>
        \left(
          \frac{1-\rho}{~ln B} + \frac{\rho}{~ln C}
        \right)
        B^{1+\rho\newalpha}
      }{\rho}{0}{1}.
    \end{align*} 
  \end{definition}
  \begin{theorem}
    \label{res:Norder-basis}
    Write $C=B^{1+\newalpha}$ and fix $k \in \N_{>0}$.  Then for $x \in
    \R_{>0}$ we have
    \begin{align*}
      \Qapprox{k}(x) &\in 
      \left[\left(\frac{1+\newalpha}{(1+\frac{\newalpha}{2})^{2}}\right)^{k},1\right]
      \Napprox{k}(x)
      .
      \qed
    \end{align*}
  \end{theorem}
  In the light of the inspiring application we now find
  $$
  \left[\left(\frac{1+\newalpha}{(1+\frac{\newalpha}{2})^{2}}\right)^{k},1\right]
  \subset [ 0.957, 1 ] .
  $$
  When we started to think about solving the recursion in
  \ref{def:Napprox} our first trial was to reuse the polynomial hills
  $\widetilde{m}^{k}$.  Yet, that didn't want to fit nicely.  Instead we
  learned from our calculations that an exponential density instead of a
  linear one would be easier to connect to the polynomial hills.  So we
  tried to approximate like this and got
  $$
  \mathring{\kappa}(\rho)
  =
  \frac{1}{~ln p} = \frac{1}{(1-\rho)~ln B + \rho ~ln C} \approx
  \myexp{ -~ln ~ln B - \rho ~ln(1+\newalpha) } =: \mathring{\eta}(\rho).
  $$
  The exponent in $\mathring{\eta}$ is chosen such that for $\rho=0$ and
  $\rho=1$ we have equality.  It turns out that this approximation is
  even better than the one before and at the same time easier to handle.
  Thus we replace the functions $\Napprox{k}$ with another family
  $\Eapprox{k}$:
  \begin{definition}
    \label{def:Eapprox}
    For $\xi\in\R$ we let $\Eapprox{0} := \Q{0}$, and
    recursively for $k>0$
    \begin{align*}
      \Eapprox{k}<\xi> &:= 
      \newalpha ~ln B
      \frint{\Eapprox{k-1}<\xi-\rho>
        \underbrace{\frac{(1+\newalpha)^{-\rho}}{~ln B}}
        _{= \mathring{\eta}(\rho)}
        B^{1+\rho\newalpha}
      }{\rho}{0}{1}.
    \end{align*} 
  \end{definition}
  \begin{theorem}
    \label{res:Eorder-basis}
    Write $C=B^{1+\newalpha}$ and fix $k \in \N_{>0}$.  Then for $x \in
    \R_{>0}$ we have
    \begin{align*}
      \Qapprox{k}(x) &\in 
      \left[
        \left(
          \frac{~ln(1+\newalpha)}{\newalpha}
          (1+\newalpha)^{\frac{1}{~ln(1+\newalpha)} - \frac{1}{\newalpha}}
        \right)^{k},
        1
      \right]
      \Eapprox{k}(x).
    \end{align*}
  \end{theorem}
  \begin{proof}
    To prove this we have to relate the function
    $\map[\mathring{\kappa}]{\left[0,1\right]}{\R_{>
        0}}{\rho}{1/~ln\left( B^{1+\rho\newalpha} \right)}$ occurring in the
    definition of $\Qapprox{k}$ to the function
    $\map[\mathring{\eta}]{\left[0,1\right]}{\R_{>
        0}}{\rho}{\mathring{\eta}(\rho)}$ replacing it in the definition
    of $\Eapprox{k}$.  Routine calculus shows that the function
    $\mathring{\kappa} / \mathring{\eta}$ is at most $1$, namely at
    $\rho=0$ and $\rho=1$, and assumes its minimum value
    $\frac{~ln(1+\newalpha)}{\newalpha}
    \left(1+\newalpha\right)^{\frac{1}{~ln(1+\newalpha)} - \frac{1}{\newalpha}}$
    at $\rho = \frac{\newalpha-~ln(1+\newalpha)}{\newalpha+~ln(1+\newalpha)}$.
  \end{proof}
  Testing this with the parameters $\newalpha = 0.\rounded{231}{902}$
  and $k=4$ from our inspiring application we obtain
  $$
  \left[
    \left(
      \frac{~ln(1+\newalpha)}{\newalpha}
      (1+\newalpha)^{\frac{1}{~ln(1+\newalpha)} - \frac{1}{\newalpha}}
    \right)^{k},
    1
  \right]
  \subset
  [0.978, 1].
  $$
  To get a better impression we have plotted the lower interval boundary
  as a function of $\newalpha$ for all three cases in \ref{fig:better}.
  \begin{figure}[tb!]
    \centering
      \psset{tickstyle=bottom,unit=3cm}
  \def\g{0.2319020119 }%
  \def\g{3 }%
  \def\legendwidth{0.17}%
  \def\ringkappa{ x \g mul 1 add 1 exch div }%
  \def\ringlambda{ 1 }%
  \def\ringnu{ x \g mul 1 \g add div 1 exch sub }%
  \def\ringeta{ x 1 \g add exch neg exp }%
  \begin{pspicture}(-0.1,-1.5\baselineskip)(1.6,1.1)%
    \footnotesize
    \psset{linewidth=.4pt}
    \psaxes[Dy=2]{->}(0,0)(-0.1,0)(1.1,1.1)%
    \psline(0,1)(-3pt,1)
    \uput[180](-3pt,1){$\frac{1}{~ln B}$}
    \uput[120](1.1,0){$\rho$}
    \psset{linewidth=.8pt}
    \rput(0.2,0.2){$\gamma=\g$}
    \psset{linecolor=black}
    \psplot{0}{1}{\ringkappa}
    \rput(0.7,0.51){%
      \psline(0,0)(\legendwidth,0)
      \uput[0](\legendwidth,0){$\mathring\kappa$}}
    \psset{linecolor=blue}
    \psplot{0}{1}{\ringlambda}
    \rput(0.7,0.90){%
      \psline(0,0)(\legendwidth,0)
      \uput[0](\legendwidth,0){$\mathring\lambda$}}
    \psset{linecolor=red}
    \psplot{0}{1}{\ringnu}
    \rput(0.7,0.77){%
      \psline(0,0)(\legendwidth,0)
      \uput[0](\legendwidth,0){$\mathring\nu$}}
    \psset{linecolor=green}
    \psplot{0}{1}{\ringeta}
    \rput(0.7,0.64){%
      \psline(0,0)(\legendwidth,0)
      \uput[0](\legendwidth,0){$\mathring\eta$}}
  \end{pspicture}
  \begin{pspicture}(-0.1,-1.5\baselineskip)(1.1,1.1)%
    \footnotesize
    \psaxes[linewidth=.4pt]{->}(0,0)(-0.1,0)(1.1,1.1)%
    \uput[120](1.1,0){$\rho$}
    \rput(0.8,0.5){$\gamma=\g$}
    \psset{linecolor=blue}
    \psplot{0}{1}{\ringkappa \ringlambda div}
    \rput(0.05,0.07){%
      \psline(0,0)(\legendwidth,0)
      \uput[0](\legendwidth,0){$\mathring\kappa/\mathring\lambda$}}
    \psset{linecolor=red}
    \psplot{0}{1}{\ringkappa \ringnu div}
    \rput(0.05,0.20){%
      \psline(0,0)(\legendwidth,0)
      \uput[0](\legendwidth,0){$\mathring\kappa/\mathring\nu$}}
    \psset{linecolor=green}
    \psplot{0}{1}{\ringkappa \ringeta div}
    \rput(0.05,0.33){%
      \psline(0,0)(\legendwidth,0)
      \uput[0](\legendwidth,0){$\mathring\kappa/\mathring\eta$}}
  \end{pspicture}
  \psset{xunit=0.375\psyunit}
  \begin{pspicture}(-0.4,-1.5\baselineskip)(7.5,1.1)
    \def\drawall#1#2{%
      \psplot[linecolor=blue]{#1}{#2}{%
        x 1 exch 1 add div
      }%
      \psplot[linecolor=red]{#1}{#2}{%
        x 1 add 4 mul x 2 add dup mul div
      }%
      \psplot[linecolor=green]{#1}{#2}{%
        x 1 add ln 
        dup 
        1 exch div 1 x div sub 
        x 1 add exch exp 
        mul 
        x div 
      }%
      {\SpecialCoor
        \rput(!#1 #2 add 2 div /x exch def x
        x 1 exch 1 add div){\rnode{blue}{}}%
        \rput(!#1 #2 add 2 div /x exch def x
        x 1 add 4 mul x 2 add dup mul div){\rnode{red}{}}%
        \rput(!#1 #2 add 2 div /x exch def x
        x 1 add ln 
        dup 
        1 exch div 1 x div sub 
        x 1 add exch exp 
        mul 
        x div 
        ){\rnode{green}{}}%
      }%
    }
    \footnotesize
    \def\xmin{0.001}\def\xmax{3.3}%
    \rput[bl](0,-1.5\baselineskip){%
      \begin{pspicture}(0,-1.5\baselineskip)(\xmax,1)
        \expandafter\expandafter\expandafter\drawall\expandafter\expandafter\expandafter{\expandafter\xmin\expandafter}\expandafter{\expandafter\xmax}%
        \psset{linewidth=.4pt}%
        \psaxes{-}(\xmax,1.1)\psline{->}(0,1.099)(0,1.1)%
        \psline(0.232,0)(0.232,-3pt)
      \end{pspicture}%
    }%
    \let\oxmax\xmax
    \def\xmin{0.001}\def\xmax{35}%
    \rput[bl](\oxmax,-1.5\baselineskip){%
      \psset{Ox=0,Dx=10,xunit=1.5mm}%
      \begin{pspicture*}(\oxmax,-1.5\baselineskip)(\xmax,1)
        \expandafter\expandafter\expandafter\drawall\expandafter\expandafter\expandafter{\expandafter\xmin\expandafter}\expandafter{\expandafter\xmax}%
        \psset{linewidth=.4pt}
        \psaxes{->}(\xmax,1.1)\psline[linestyle=dashed](\oxmax,-6pt)(\oxmax,1)
        \uput[-120](\xmax,0){$\gamma$}
        \psset{linewidth=.8pt}
        \rput(5,0.9){%
          \psset{xunit=\psyunit}
          \psline[linecolor=green](0,0)(\legendwidth,0)
          \uput[0](\legendwidth,0){$
            \frac{~ln(1+\gamma)}{\gamma}
            (1+\gamma)^{\frac{1}{~ln(1+\gamma)}-\frac{1}{\gamma}}
            $}
        }
        \rput(10,0.7){%
          \psset{xunit=\psyunit}
          \psline[linecolor=red](0,0)(\legendwidth,0)
          \uput[0](\legendwidth,0){$
            \frac{1+\gamma}{(1+\frac{\gamma}{2})^{2}}
            = ~min \Frac{\mathring\kappa}{\mathring\nu}
            $}
        }
        \rput(15,0.5){%
          \psset{xunit=\psyunit}
          \psline[linecolor=blue](0,0)(\legendwidth,0)
          \uput[0](\legendwidth,0){$
            \frac{1}{1+\gamma}
            = ~min \Frac{\mathring\kappa}{\mathring\lambda}
            $}
        }
      \end{pspicture*}
    }
  \end{pspicture}

%
    \caption{Comparing the quality of $\Kapprox{}$, $\Napprox{}$,
      $\Eapprox{}$: the top pictures show the integral kernels and their
      ratios for a specific $\newalpha$, the lower pictures show the
      minimum of the ratios as a function of $\newalpha$}
    \label{fig:better}
  \end{figure}
  We see that for small values of $\newalpha$ we obtain good approximations
  of $\Qapprox{}$ and all our attempts give only weak results for large
  $\newalpha$, but the one with $\mathring{\eta}$ is always best.

  If we now rewrite the recursion to one for
  $$
  \Eapprox[~norm]{k}<\xi>
  = \frac{\Eapprox{k}<\xi>}{\newalpha^{k} B^{k+\xi\newalpha} (1+\newalpha)^{-\xi}}
  = \frac{\Eapprox{k}<\xi>}{\newalpha^{k} 
    B^{k+\xi \left(\newalpha  - \frac{~ln(1+\newalpha)}{~ln B}\right)}
  }
  $$
  we find that $\Eapprox[~norm]{k} = \mathcal{M} \Eapprox[~norm]{k-1}$.
  So we'll obtain the solution from the polynomial hills as in
  \ref{res:order-approx} for $\Kapprox{k}$:
  $$
  \Eapprox[~norm]{k}<\xi> = \frint{ \myexp{-\rho (\newalpha ~ln B -
      ~ln(1+\newalpha)) } \widetilde{m}^{k}(\xi-\rho) }{\rho}{0}{\infty}.
  $$
  The only difference is that instead of $\newalpha$ we have
  $\newalpha-\frac{~ln(1+\newalpha)}{~ln B}$.  With this replacement
  \ref{res:order-approx} becomes:
  \begin{theorem}
    \label{res:order-approx-exp}
    For any $\xi \in \R$ we have
    \begin{align*}
      \Eapprox[~norm]{k}<\xi>
      &=
      \frint{B^{-\rho\left(\newalpha-\frac{~ln(1+\newalpha)}{~ln B}\right)} \widetilde{m}^{k}(\xi-\rho)}{\rho}{0}{\xi}
      \\&=
      B^{-\xi\left(\newalpha-\frac{~ln(1+\newalpha)}{~ln B}\right)}
      \frint{B^{\rho\left(\newalpha-\frac{~ln(1+\newalpha)}{~ln B}\right)} \widetilde{m}^{k}(\rho)}{\rho}{0}{\xi},
      \\
      &=
      \frint{ 
        \left(\frac{B^{\newalpha}}{1+\newalpha}\right)^{-\rho}
        \widetilde{m}^{k}(\xi-\rho)}{\rho}{0}{\xi}
      \\&=
      B^{-\xi\left(\newalpha-\frac{~ln(1+\newalpha)}{~ln B}\right)}
      \frint{ \left(\smashedunderbrace{\frac{B^{\newalpha}}{1+\newalpha}}
          _{=\frac{C/~ln C}{B/~ln B}}
        \right)^{\rho} \widetilde{m}^{k}(\rho)}{\rho}{0}{\xi},
      \smashedheights
      \\\allowpagebreak
      \Eapprox[~norm]{k}<\xi>
      &=
      \frac{1}{\left(-\newalpha ~ln B+~ln(1+\newalpha)\right)^{k}}
      \\&\phantom{=\;}
      \begin{aligned}[t]
        &
        \left(
          \presum{i}{0}{\floor{\xi}}
          \binom{k}{i}
          (-1)^{i} B^{-(\xi-i)\left(\newalpha-\frac{~ln(1+\newalpha)}{~ln B}\right)}
        \right.
        \\
        &\hphantom{\left(\vphantom{\presum{i}{0}{\floor{\xi}}}\right.}\left.-
          \presum{l}{0}{k-1}
          \left( -\newalpha ~ln B+~ln(1+\newalpha) \right)^{l}
          \cdot
          \mathcal{D}_{\xi}^{k-l-1}\widetilde{m}^{k}(\xi)
        \right)\mathop{\vphantom{\biggr)}}\limits_{\displaystyle\vphantom{\int}},
      \end{aligned}  
      \\\allowpagebreak
      \Eapprox[~norm]{k}<\xi>
      &=
      \begin{aligned}[t]
        &
        \presum{i}{0}{\floor{\xi}}
        \binom{k}{i}
        (-1)^{i}
        \frac{
          ~cutexp_{k}
          \left(
            -(\xi-i)\left(\newalpha ~ln B-~ln(1+\newalpha)\right)
          \right)
        }{\left(-\newalpha ~ln B+~ln(1+\newalpha)\right)^{k}},
      \end{aligned}
    \end{align*}
    where $~cutexp_{k}(\zeta) = ~exp(\zeta) - \presum{l}{0}{k-1} \frac{
      \zeta^{l} }{l!} = \sum_{l\geq k} \frac{ \zeta^{l} }{l!}$.
  \end{theorem}

  \section{Non-squarefree numbers are negligible}
  \label{sec:squareful-negligible}

  Considering $\Q{k}(x)$ we immediately observe that the ordering of the
  counted prime lists are not important and we can group together many
  such elements.  To get a precise picture, we define the \emph{sorting}
  of a tuple $P = (p_{1}, \ldots, p_{k})$ in the following way:
  \begin{gather*}
    S(P) := \left(
      \vphantom{\Bigm.}
      \Set{j\leq k; ~rank_{j} P = i }
      \vphantom{\Bigm.}
    \right)_{i \leq k},
  \end{gather*}
  where $~rank_{j} P = \#\Set{ p_{l}; l \leq k \wedge p_{l} \leq p_{j}}$.
  Now the count $\Q{k}(x)$ can be partitioned using the sets
  $$
  A_{S}(x) :=     \Set{
    P = ( p_{1}, \dots, p_{k} ) \in \left(\P{B,C}\right)^{k};
    \begin{array}{@{}l@{}}
      n = p_{1} \dots p_{k} \leq x,\\
      S(P) = S
    \end{array}
  },
  $$
  namely $\Q{k}(x) = \sum_{S} \# A_{S}(x)$ where $S$ runs over all
  possible sortings.  The above intuition would imply that many of these
  sets are essentially equal.  We group them by their \emph{type}
  $$
  T(S) := \left( \#S_{i} \right)_{i \leq k}.
  $$
  Given any type $T = (T_{1}, \ldots, T_{r})$, there are exactly
  $\binom{k}{T} = \frac{k!}{T_{1}!  \cdots T_{r}!}$ different
  sortings~$S$ of type~$T$.  This corresponds to possible reorderings
  of a specific vector $P \in A_{S}(x)$ for a sorting~$S$ of type~$T$.
  The type of such a vector $P$ is defined to be $T(S)$. It is clear 
  that the type of~$P$ is invariant under permutations, yet not its 
  sorting.
  \begin{lemma}
    Let $T$ be a type for $k$ elements.
    \begin{enumerate}
    \item There exists a sorting $S(T)$ of type $T$ such that all
      vectors in $A_{S(T)}(x)$ are \weaklyisotone{}.
    \item If $T(S) = T$ then there is a permutation $\sigma$ of~$k$
      elements such that for all $x$ we have $A_{S}(x) =
      A_{S(T)}(x)^{\sigma}$.
    \item More precisely, for any sorting $S$ of~$k$ elements the
      following are equivalent:
      \begin{enumerate}
      \item $T(S) = T$.
      \item $\exists \sigma\colon S = S(T)^{\sigma}$.
      \item $\exists \sigma\colon \forall x\colon A_{S}(x) =
        A_{S(T)}(x)^{\sigma}$.
        \qed
      \end{enumerate}
    \end{enumerate}
  \end{lemma}
  Noting that $\#\Set{S; T(S)=T} = \binom{k}{T}$ we have
  \begin{align*}
    \Q{k}(x)
    &=
    \sum_{T} \sum_{S\colon T(S) = T} \# A_{S}(x) = \sum_{T} \binom{k}{T} \# A_{S(T)}(x)
  \end{align*}
  On the other hand we have
  $
  \Qclean{k}(x)
  =
  \sum_{T}
  1 \cdot
  \# A_{S(T)}(x)
  $.
  In particular, we can deduce
  $$
  \Qclean{k}(x) < \Q{k}(x) \leq k! \cdot \Qclean{k}(x).
  $$
  Actually, for large $B$ (and $C$ and $x$) we have $\Qclean{k}(x) \sim
  \frac{1}{k!}\Q{k}(x)$.  This stems from the following fact that
  $\#A_{S}(x)$ is asymptotically much smaller than
  $\#A_{S(1,\dots,1)}(x)$ for any sorting $S$ of $k$ elements of type
  different from $(1,\dots,1)$.
  \begin{lemma}
    \label{res:squarefulterms}
    For any sorting $S$ of $k$ elements of type different from
    $(1,\dots,1)$ there is a sorting $S'$ of $k-1$ elements such that
    we have $\#A_{S}(x) \leq \#A_{S'}(x/B) \leq \frac{x}{B}$.
  \end{lemma}
  \begin{proof}
    Take $S$ as specified.  Let $t$ be a position which does not occur
    as a singleton in $S$.  Further, say $t \in S_{\tau}$, and let $r
    = \#S_{\tau} \geq 2$.  Let $S^{-}$ be the sorting with $S_{\tau}$
    removed, and $S'$ the sorting with $t$ removed.  (Retaining the
    old indexing is easier, yet then indices run over $\Set{1,\dots,k}
    \setminus S_{\tau}$, or $\Set{1,\dots,k} \setminus \Set{t}$,
    respectively.)  Then
    $$
    \#A_{S}(x) = \primesum{ \# A'_{S^{-}}(x/p_{t}^{r}) }{p_{t}}{B}{C}
    \leq \primesum{ \# A'_{S^{-}}(x/p_{t}B) }{p_{t}}{B}{C} =
    \#A_{S'}(x/B).
    $$
    Here, $A'$ denotes a variant of $A$ consisting of prime vectors
    where at the positions $j$ which should have smaller primes than
    $p_{t}$ originally still satisfy that, and same for positions with
    primes larger than $p_{t}$.
    For the inequality note that $S^{-} = (S')^{-}$.  Since obviously
    $\#A_{S}(z) \leq z$ we are done.
  \end{proof}
  Combining this with $\sum_{S} \#A_{S}(x) = \Q{k}(x) \sim
  \Qapprox{k}(x) \in \Theta \left( \frac{x}{~ln{B}} \right)$, shows
  that there must be a large summand, which can be only
  $\#A_{S(1,\dots,1)}(x)$.


  The number $s(k) = \sum_{T} \binom{k}{T}$ of sortings of $k$
  elements is called ordered Bell number.  We can also recursively
  define them: $s(0) = 1$, $s(k) = \frsum{ \binom{k}{r} s(r)
  }{r}{0}{k-1}$.  According to \cite{wil94}, page 175f, we have $s(k)
  = \frac{k!}{2 ~ln^{k+1}2} + \mathcal{O}\left((0.16)^{k}k!\right)$.
  In particular, $s(k)$ is small in comparison to $2^{k-1} k!$.  Using
  \ref{res:squarefulterms} for a comparison yields the ---~now
  immediate~--- following
  \begin{lemma}
    \label{res:pi-kappa-comparison}
    We have 
    $$
    \left| \Qclean{k}(x) - \frac{1}{k!}\Q{k}(x) \right| 
    \leq
    \left(
      2^{k-1} - \frac{s(k)}{k!}
    \right)
    \frac{x}{B}
    < 2^{k-1} \frac{x}{B}.
    $$
  \end{lemma}
  \begin{proof}
    $\left| k! \cdot \Qclean{k}(x) - \Q{k}(x) \right| \leq \sum_{T}
    \left( k! - \binom{k}{T} \right) \frac{x}{B} = \left( k! 2^{k-1} - s(k)
    \right) \frac{x}{B}$.
  \end{proof}
  Compared to the error bound in $\left| \Q{k}(x)-\Qapprox{k}(x)
  \right| \leq \Qerror{k}(x) \in \mathcal{O}(\frac{x}{\sqrt{B}})$ this
  is negligible when $B$ is large.  Here we assume $k\geq 2$ since the
  present observations are irrelevant for $k=1$, namely $\Qclean{1} =
  \Q{1}$.
\end{fullversion}

\section{Results on coarse-grained integers}
We are going to combine the results of the last section with
\ref{res:Kapprox-bounds} and \ref{res:Kerror-upperbound} to finally
arrive at our main theorem.

Analogously to \ref{def:Qapprox}, we define an approximation function
for the function $\Qclean{k}(x)$.
\begin{definition}
  \label{def:Qcleanapprox}%
  For $x\geq 0$ and $k \geq 0$ we define
  $$
  \Qcleanapprox{k}(x) := \frac{1}{k!} \Qapprox{k}(x)
  \text{ and } 
  \Qcleanerror{k}(x) := \frac{1}{k!} \Qerror{k}(x) + 2^{k-1}
  \frac{x}{B}.
  $$  
\end{definition}
Similarly to \ref{def:Qapprox} we can also recursively define
$\Qcleanapprox{k}(x)$ by
\begin{align*}
  \Qcleanapprox{0}(x) &=
  \begin{cases}
    0 & \text{if $x<1$},
    \\
    1 & \text{if $1\leq x$},
  \end{cases}
  &%
  \Qcleanapprox{k}(x) &= \frac{1}{k}
  \frint{\frac{\Qcleanapprox{k-1}(x/p_{k})}{~ln p_{k}}}{p_{k}}{B}{C}.
\end{align*}
It is also possible to define $\Qcleanerror{k}(x)$ similarly based on
\ref{def:Qapprox}. We can now describe the behavior of $\Qclean{k}$
nicely and give our main result.
\begin{theorem}
  \label{res:cleanapprox}
  Given $x \in \R_{> 0}$ and $k \in \N$. Then the inequality
  $$
  \left|\Qclean{k}(x) - \Qcleanapprox{k}(x)\right| \leq
  \Qcleanerror{k}(x)
  $$
  holds.\iffullversion\else\qed\fi
\end{theorem}
\begin{fullversion}
\begin{proof}  
  By \ref{res:pi-kappa-comparison} we have
  $$
  \left| \Qclean{k}(x) - \frac{1}{k!}\Q{k}(x) \right| < 2^{k-1}
  \frac{x}{B}.
  $$
  Thus using the triangle inequality and \ref{res:approx}, we obtain
  $$
  \left| \Qclean{k}(x) - \frac{1}{k!}\Qapprox{k}(x) \right| <
  \frac{1}{k!} \Qerror{k}(x) + 2^{k-1} \frac{x}{B},
  $$
  which proves the claim.
\end{proof}
\end{fullversion}
\begin{theorem}\label{res:final-count-lambda}
  Fix $k \geq 2$. Then for any $\epsilon>0$ and $B$ tending to
  infinity, there are for $x \in \left[ B^{k} (1+\epsilon), C^{k}
    (1-\epsilon) \right]$ values $\widetilde{s}, \widehat{s} \in
  \left[\frac{1}{(1+\newalpha)^{k}},1\right]$ such that
  \begin{align*}
    \Qcleanapprox{k}(x) &= \frac{\widetilde{s}}{k!} \Kapprox{k}(x),
    &
    \Qcleanerror{k}(x) &= \frac{\widehat{s}}{k!} \Kerror{k}(x) +
    2^{k-1} \frac{x}{B}.\iffullversion\else\qed\fi
  \end{align*}
\end{theorem}
\begin{fullversion}
\begin{proof}
  By \ref{def:Qcleanapprox} we have
  $$
  \Qcleanapprox{k}(x) = \frac{1}{k!} \Qapprox{k}(x).
  $$
  \ref{res:order-basis} tells us that there is a value
  $\widetilde{s} \in \left[\frac{1}{(1+\newalpha)^{k}},1\right]$ such
  that
  $$
  \Qapprox{k}(x) = \widetilde{s} \Kapprox{k}(x),
  $$
  implying that for the same value $\widetilde{s}$ we have
  $$
  \Qcleanapprox{k}(x) = \frac{\widetilde{s}}{k!} \Kapprox{k}(x).
  $$
  Considering $\Qcleanerror{k}(x)$ we have by \ref{def:Qcleanapprox}
  that
  $$
  \Qcleanerror{k}(x) = \frac{1}{k!} \Qerror{k}(x) + 2^{k-1}
  \frac{x}{B}.
  $$
  Applying \ref{res:order-basis} gives a value $\widehat{s} \in
  \left[\frac{1}{(1+\newalpha)^{k}},1\right]$ such that
  $$
  \Qerror{k}(x) = \widehat{s} \Kerror{k}(x),
  $$
  which directly gives
  $$
  \Qcleanerror{k}(x) = \frac{\widehat{s}}{k!} \Kerror{k}(x) + 2^{k-1}
  \frac{x}{B}.
  $$
  This proves the theorem.
\end{proof}
\end{fullversion}
Unrolling our results on $\Kapprox{k}(x)$ and $\Kerror{k}(x)$, namely
\ref{res:Kapprox-bounds} and \ref{res:Kerror-upperbound}, gives a
slightly weaker result.
\begin{theorem}
  \label{res:Brough+Csmooth}
  Let $B<C=B^{1+\newalpha}$ with $\alpha \geq \frac{~ln B}{\sqrt{B}}$
  and fix $k\geq 2$.  Then for any (small) $\epsilon>0$ and $B$
  tending to infinity we have for $x \in \left[ B^{k} (1+\epsilon),
    C^{k} (1-\epsilon) \right]$ a value $\widetilde{a} \in \left[
    \frac{\newalpha^{k-1} \clower{k}}{k! (1+\newalpha)^k},
    \frac{1}{k!}  \right]$ with $\clower{k} =
  ~min\left(2^{-4} \frac{\epsilon^{k}}{k!}, 2^{-k}
    \frac{\epsilon^{k-1}}{(k-1)!}\right) $ such that
  $$
  \left|\Qclean{k}(x) - \widetilde{a} \frac{x}{~ln B}\right| \leq
  (2^k-1) \newalpha^{k-2} (1+\newalpha) \cdot \frac{x}{\sqrt{B}} +
  2^{k-1} \frac{x}{B}.
  \iffullversion\else\qed\fi
  $$
\end{theorem}
\begin{fullversion}
\begin{proof}
  By \ref{res:final-count-lambda} we have values $\widetilde{s},
  \widehat{s} \in \left[\frac{1}{(1+\newalpha)^{k}},1\right]$ such
  that
  \begin{align*}
    \Qcleanapprox{k}(x) &= \frac{\widetilde{s}}{k!} \Kapprox{k}(x),
    &
    \Qcleanerror{k}(x) &= \frac{\widehat{s}}{k!} \Kerror{k}(x) +
    2^{k-1} \frac{x}{B}.
  \end{align*}
  By \ref{res:Kapprox-bounds} we have for $\epsilon$ small enough that
  $$
  \Kapprox{k}(x) \in \left[ \clower{k}, 1 \right]
  \frac{\newalpha^{k-1} x}{~ln B}
  $$
  and by \ref{res:Kerror-upperbound} that
  $$
  \Kerror{k}(x) \leq (2^{k} -1) \newalpha^{k-2} (1+\newalpha) \cdot
  \frac{x}{\sqrt{B}}.
  $$
  This gives the claim.
\end{proof}
\end{fullversion}

\endnote{%
  \paragraph{MN20110923: A test case} %
  Let's consider the quantity $\tau_{2}(N,\gamma)$ from Conjecture 7
  in \cite{mau92b} as a test case.  We only allow ourselves to replace
  sometimes replace $\leq$ with $<$ to better fit into our notation,
  which does not change the asymptotics.  In our language,
  $\tau_{2}(N,\gamma) = \Qclean[N^{\gamma},\infty]{2}(N) =
  \Qclean[N^{\gamma},N^{1-\gamma}]{2}(N)$, in particular $k=2$,
  $B=N^{\gamma}$, $C=N^{1-\gamma}$, $x=N$,
  $\alpha=\frac{1}{\gamma}-2$, $\xi=1$.  Directly we see that this
  term vanishes if $\gamma\leq\frac{1}{2}$, 
  so we assume $\gamma < \frac{1}{2}$.  (In Maurer's version
  $\tau_{2}(N,\frac12) = 1$ occasionally, namely if $N$ is a prime
  square.)  We find
  \begin{align*}
    \Qclean[N^{\gamma},N^{1-\gamma}]{2}(N) 
    &= \frac12 \left(
      \Q[N^{\gamma},N^{1-\gamma}]{2}(N)
      -\Q[N^{\gamma},N^{1-\gamma}]{1}(N^{\frac12}) \right),
  \end{align*}
  and
  $\Q[N^{\gamma},N^{1-\gamma}]{2}(N) \in \left[
    \frac1{(1+\alpha)^{2}}, 1
  \right] \Kapprox{2}<1> +
  \left[
    -1,1
  \right]\Kerror{2}<1>$.
  Further, we already computed (see endnote \ref{en:lambdatilde})
  \begin{align*}
      \Kapprox{2}<1>
      &= x \left(
        \alpha \frac{1}{~ln B} - \frac{1}{~ln^{2}B} +
        \frac{B^{-\alpha}}{~ln^{2}B}
      \right)
      \\
      &=
      \left(
        \frac{1}{\gamma^{2}} - \frac{2}{\gamma}
      \right)
      \cdot
      \frac{N}{~ln N}
      - \frac{1}{\gamma^{2}} \cdot \frac{N}{~ln^{2}N} +
      \frac{1}{\gamma^{2}} \cdot \frac{N^{2\gamma}}{~ln^{2} N}.
  \end{align*}
  Further, $\Kerror{2}<1>$ is bounded by $\sqrt{B} =
  N^{\frac{\gamma}{2}}$ for $B$ large enough (see endnote
  \ref{en:lambdahat}).
  Also, $\Q[N^{\gamma},N^{1-\gamma}]{1}(N^{\frac12})$ is bound by
  $N^{\frac12}$.  So these errors are negligible and
  \begin{align*}
    \Qclean[N^{\gamma},N^{1-\gamma}]{2}(N) 
    &\in
    \left[
      \frac1{2(1+\alpha)^{2}}, \frac{1}{2}
    \right] \Kapprox{2}<1> +
    \left[
      -1,1
    \right] N^{\frac{1}{2}}
    \\
    &=
    \begin{aligned}
      &\left[
        \frac{\gamma^{2}}{2(\gamma-1)^{2}}, \frac{1}{2}
      \right] 
      \left(
        \left(
          \frac{1}{\gamma^{2}} - \frac{2}{\gamma}
        \right)
        \cdot
        \frac{N}{~ln N}
        - \frac{1}{\gamma^{2}} \cdot \frac{N}{~ln^{2}N} +
        \frac{1}{\gamma^{2}} \cdot \frac{N^{2\gamma}}{~ln^{2} N}
      \right)+
      \bigO{1} N^{\frac{1}{2}}.
    \end{aligned}
  \end{align*}
  The major term thus is $c \frac{N}{~ln N}$ with $c \in
  \left[
    \frac{1-2\gamma}{2(1-\gamma)^{2}},
    \frac{1-2\gamma}{2\gamma^{2}}
  \right]$.
  For $\gamma=\frac{1-\epsilon}{2}$ the latter intervall is $
  \left[ \frac{2\epsilon}{(1+\epsilon)^{2}},
    \frac{2\epsilon}{(1-\epsilon)^{2}} \right]$.
  If, for example, $\gamma = \frac{1}{2} - \frac{~ln 2}{~ln N}$ then
  this interval is around $\frac{~ln 16}{~ln N}$ and thus the
  asymptotics change here!
  Maurer's conjectured value would be $c = ~ln(1-\gamma)-~ln(\gamma) =
  ~ln\left(\frac{1}{\gamma}-1\right)$ which is in the interval 
  for $\gamma<\frac{1}{2}$.

  \paragraph{Can a better analysis give a more precise answer?}
  
  Well, we can get hold of $\Qclean{2}<1>$ as follows: We know that
  $\Qclean{2}<1>$ is close to $2 \Qapprox{2}<1>$.  And we would like
  to see which asymptotic coefficient has to be put with its leading
  term $\frac{N}{~ln N}$. Consider the double integral $\Qapprox{2}<1>
  = \frint{ \frint{ \frac{1}{~ln p_{1} ~ln p_{2}} }{p_{1}}{B}{\frac{B
        C}{p_{2}}} }{p_{2}}{B}{C} $ with $B=N^{\gamma}$,
  $C=N^{1-\gamma}$.  Then take the derivative of $T := \frac{~ln N}{N}
  \Qapprox{2}<1>$ with respect to $\gamma$.  This quantity has only
  single integrals!  They can all be expressed by logarithmic
  integrals $~Li$.  Next, as approximations replace logarithmic
  integrals with their major term, ie.\ $~Li(t) \approx \frac{t}{~ln
    t}$ for large $t$, and compute the limit for $N \longto \infty$.
  Compare to the derivative of $~ln(\frac{1}{\gamma}-1)$ with respect
  to $\gamma$.  These are equal.
  
  Though this doesn't prove anything, it is an indicator that the
  term given by \cite{mau92b} is indeed the right one.

  \todo{Can the shown path be turned into a proof?}  %
}%
\ifdraft\textsuperscript{,}\fi%
\endnote{
\todo{Repeat the previous with $\Eapprox{k}$ in place of
  $\Kapprox{k}$.} Note that \ref{res:Kapprox-bounds} with $\alpha$
replace by $\alpha-\frac{~ln(1+\alpha)}{~ln B}$ is the appropriate
result for $\Eapprox{k}$.  \todo{So} we have $$ \Eapprox{k}(x) \in
[\clower{k},1] \left(\alpha-\frac{~ln(1+\alpha)}{~ln
    B}\right)^{k-1}\frac{x}{~ln B}.
$$
Hm\dots\ \ref{res:Kapprox-bounds} with the replacement gives
$$
\Eapprox[~norm]{k}(x) \in [\clower{k},1]
\frac{1}{\left(\alpha-\frac{~ln(1+\alpha)}{~ln B}\right) ~ln B}
$$
and $\Eapprox{k}(x) = \alpha^{k} (1+\alpha)^{-\xi} x \cdot 
\Eapprox[~norm]{k}(x)$ with $x=B^{k+\xi\alpha}$ (ie.\ $\xi = \frac{~ln
  x}{\alpha ~ln B}-\frac{k}{\alpha}$ and
$(1+\alpha)^{-\xi} = 
(1+\alpha)^{-\frac{~ln
  x}{\alpha ~ln B}+\frac{k}{\alpha}}
= 
~exp\treatparameter({\frac{k ~ln(1+\alpha)}{\alpha}})
x^{-\frac{~ln (1+\alpha)}{\alpha ~ln B}}$)
$$
\Eapprox{k}(x) \in [\clower{k},1] \frac{\alpha^{k-1}
  ~exp\treatparameter({\frac{k ~ln(1+\alpha)}{\alpha}})
  x^{1-\frac{~ln (1+\alpha)}{\alpha ~ln B}}}{\left(1-\frac{~ln(1+\alpha)}{\alpha~ln B}\right)
  ~ln B}.
$$
\todo{Is this correct?  How can this be?  The exponent of $x$ is
  smaller than $1$ and still the remaining stuff is independent of
  $x$. How far can $x$ and $x^{1-\frac{~ln (1+\alpha)}{\alpha ~ln B}}$
  be apart?  Actually, $x$ is bounded by $\alpha$ somehow, namely $x
  \in [B^{k},B^{k(1+\alpha)}]$.  Is that enough to explain this
  discrepancy? Is the reasoning for \ref{res:order-approx-exp}
  correct?  Is it correct for $k=0$?}
\todo{DL2012/02/02: I think that is is \emph{not} ok to simply take the bounds for $\Kapprox{k}(x)$ and substitute $\newalpha$. Example: In the normalization of this function the value of $\alpha$ is once substituted by $\newalpha-\frac{~ln(1+\newalpha)}{~ln B}$, but another $\newalpha$ is \emph{not}. Furthermore, the proof of the bounds heavily depends on specific properties of $\Kapprox{k}(x)$.}
}

\section{Numeric evaluation}
\label{sec:last-words}

To discuss the quality of our results we consider again the example
parameters $B=1100 \cdot 10^{6}$, $C=2^{37}-1$, $\newalpha =
~ln(C)/~ln(B)-1 = 0.\rounded{231}{902}$, $k = 4$ from our inspiring
application.  For $k=2$, $k=3$ we can give similar pictures; of
course, the errors are even smaller in these cases.  At present we do
not have efficient algorithms for computing $\Q{k}$ itself.  However,
based on our estimates we can compute values for encapsulating
intervals in three variants, listed in increasing quality:
\begin{itemize}
\item $\lambda$~estimate (\ref{res:order-basis}): 
$$
\left[ \frac{1}{(1+\newalpha)^{k}}
  \Kapprox{k}(x)-\cerror{k} \newalpha^{k} x, \quad
  \Kapprox{k}(x)+\cerror{k} \newalpha^{k} x
\right].
$$ 
\item
  $\eta$~estimate (\ref{res:Eorder-basis}): 
  $$
  \left[ \left( \frac{~ln(1+\newalpha)}{\newalpha}
      (1+\newalpha)^{\frac{1}{~ln(1+\newalpha)} - \frac{1}{\newalpha}}
    \right)^{k} \Eapprox{k}(x)-\cerror{k} \newalpha^{k} x, \quad
    \Eapprox{k}(x)+\cerror{k} \newalpha^{k} x
  \right].
  $$ 
\item $\Q[]{}$~estimate (\ref{res:approx}):
  $$
  \left[ \Qapprox{k}(x)-\Qerror{k}(x), \quad
    \Qapprox{k}(x)+\Qerror{k}(x) \right].
  $$ 
\end{itemize}
The $\Q[]{}$~estimate was easiest to obtain and is of course the most
accurate one, however, it is difficult to evaluate.  The
$\lambda$~estimate was easy to obtain and compute.  But it is of
course the least accurate of the three.  The $\eta$~estimate was
slightly more difficult to find, is as easy to evaluate as the prior,
and it is much more accurate.  As usual we write $x = B^{k+\xi
  \newalpha}$ and use $\xi$ as a running parameter.

\newdimen\myy

\begin{figure}[t!]
  \centering%
  \myy=20mm%
  \def\pshlabel#1{\small$#1$}%
  \def\psvlabel#1{\small$#1\cdot{10^{40}}$}%
  \psset{xunit=29mm,yunit=13\myy}%
  \def\myk{4}\def\ymax{0.23}%
\setbox99\hbox{\small\psvlabel{0.3}\hskip10pt}%
\begin{pspicture*}(-\wd99,-1.2\baselineskip)(\myk.2,\ymax)%
  \readdata{\lbdatalambda}{grainyfinalpicsdat-lambdaerrlbk\myk.dat}
  \readdata{\ubdatalambda}{grainyfinalpicsdat-lambdaerrubk\myk.dat}
  \readdata{\lbdataeta}{grainyfinalpicsdat-etaerrlbk\myk.dat}
  \readdata{\ubdataeta}{grainyfinalpicsdat-etaerrubk\myk.dat}
  \readdata{\lbdatakappa}{grainyfinalpicsdat-kappaerrlbk\myk.dat}
  \readdata{\ubdatakappa}{grainyfinalpicsdat-kappaerrubk\myk.dat}
  \psclip{
    \pscustom[linestyle=none]{%
      \listplot[plotstyle=curve,showpoints=false, dotstyle=triangle,linecolor=green]{\lbdatalambda}
      \lineto(\myk,\ymax)
      \lineto(0,\ymax)
    }
    \pscustom[linestyle=none]{%
      \listplot[plotstyle=curve,showpoints=false, dotstyle=triangle,linecolor=red]{\ubdatalambda}
      \lineto(\myk,0)
      \lineto(0,0)
    }
  }
  \psframe*[linecolor=lightgray](0,0)(\myk,\ymax)
  \endpsclip
  \psclip{
    \pscustom[linestyle=none]{%
      \listplot[plotstyle=curve,showpoints=false, dotstyle=triangle,linecolor=yellow]{\lbdataeta}
      \lineto(\myk,\ymax)
      \lineto(0,\ymax)
    }
    \pscustom[linestyle=none]{%
      \listplot[plotstyle=curve,showpoints=false, dotstyle=triangle,linecolor=cyan]{\ubdataeta}
      \lineto(\myk,0)
      \lineto(0,0)
    }
  }
  \psframe*[linecolor=gray](0,0)(\myk,\ymax)
  \endpsclip
  \psclip{
    \pscustom[linestyle=none]{%
      \listplot[plotstyle=curve,showpoints=false, dotstyle=triangle,linecolor=yellow]{\lbdatakappa}
      \lineto(\myk,\ymax)
      \lineto(0,\ymax)
    }
    \pscustom[linestyle=none]{%
      \listplot[plotstyle=curve,showpoints=false, dotstyle=triangle,linecolor=cyan]{\ubdatakappa}
      \lineto(\myk,0)
      \lineto(0,0)
    }
  }
  \psframe*[linecolor=black](0,0)(\myk,\ymax)
  \endpsclip
  \psaxes[Dx=0.5, Oy=0, Dy=0.1, linewidth=0.4pt]{->}(0,0)(0,0)(\myk.2,\ymax)
  \uput[-120]{0}(\myk.2,0){$\xi$}%
\end{pspicture*}%
 
  \caption{Absolute behavior of the estimates for $\Q{k}<\xi>$ .  The
    light gray area shows the $\lambda$~estimate, the dark gray area
    the $\eta$~estimate, and the black area (well, yes) shows the
    $\Q[]{}$~estimate.  The parameters are $B = 1100\cdot10^{6}$, $C =
    2^{37}-1$, $\newalpha = ~ln(C)/~ln(B)-1 =
    0.\protect\rounded{231}{902}$, $k=4$.}
  \label{fig:abs-behavior}
\end{figure}
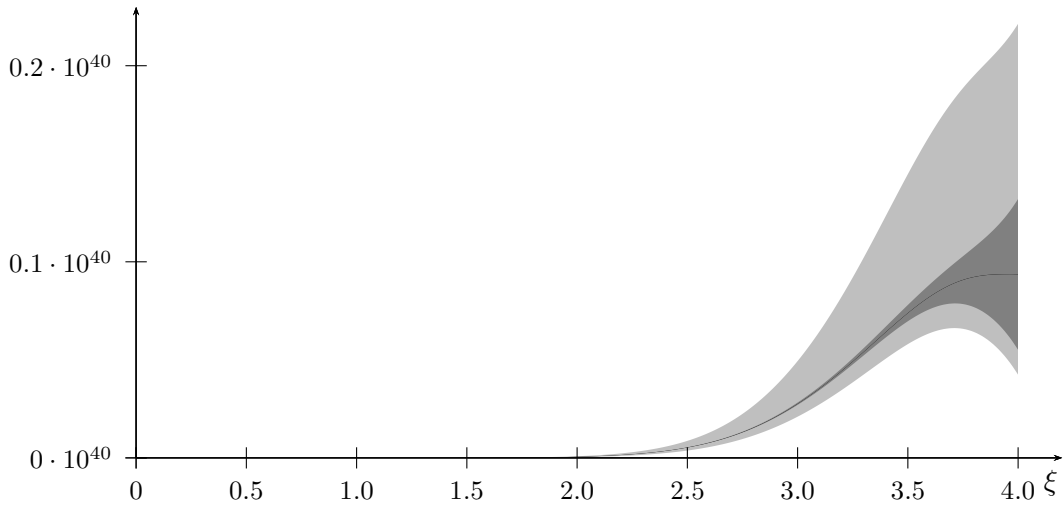

\ref{fig:abs-behavior} shows the absolute behavior of all estimates.
We observe that the absolute errors at the right margin are huge.
This is expected as also the error estimates in the prime number
theorem only bound the relative error.  However, the picture
completely conceals information about the middle and the left part of
the interval $[0,k]$.

\begin{figure}[t!]
  \centering%
  \myy7mm%
  \def\pshlabel#1{\small$#1$}%
  \def\psvlabel#1{\small$#1\cdot{10^{-3}}$}%
  \psset{xunit=29mm,yunit=16\myy}%
  \def\myk{4}\def\ymax{0.4}%
\setbox99\hbox{\small\psvlabel{0.3}\hskip10pt}%
\begin{pspicture*}(-\wd99,-1.2\baselineskip)(\myk.2,\ymax)%
  \readdata{\lbdatalambda}{grainyfinalpicsdat-lambdanormerrlbk\myk.dat}
  \readdata{\ubdatalambda}{grainyfinalpicsdat-lambdanormerrubk\myk.dat}
  \readdata{\lbdataeta}{grainyfinalpicsdat-etanormerrlbk\myk.dat}
  \readdata{\ubdataeta}{grainyfinalpicsdat-etanormerrubk\myk.dat}
  \readdata{\lbdatakappa}{grainyfinalpicsdat-kappanormerrlbk\myk.dat}
  \readdata{\ubdatakappa}{grainyfinalpicsdat-kappanormerrubk\myk.dat}
  \def\drawallthree(#1)(#2){%
  \psclip{
    \pscustom[linestyle=none]{%
      \listplot[plotstyle=curve,showpoints=false, dotstyle=triangle,linecolor=green]{\lbdatalambda}
      \lineto(\myk,\ymax)
      \lineto(0,\ymax)
    }
    \pscustom[linestyle=none]{%
      \listplot[plotstyle=curve,showpoints=false, dotstyle=triangle,linecolor=red]{\ubdatalambda}
      \lineto(\myk,0)
      \lineto(0,0)
    }
  }
  \psframe*[linecolor=lightgray](#1)(#2)
  \endpsclip
  \psclip{
    \pscustom[linestyle=none]{%
      \listplot[plotstyle=curve,showpoints=false, dotstyle=triangle,linecolor=yellow]{\lbdataeta}
      \lineto(\myk,\ymax)
      \lineto(0,\ymax)
    }
    \pscustom[linestyle=none]{%
      \listplot[plotstyle=curve,showpoints=false, dotstyle=triangle,linecolor=cyan]{\ubdataeta}
      \lineto(\myk,0)
      \lineto(0,0)
    }
  }
  \psframe*[linecolor=gray](#1)(#2)
  \endpsclip
  \psclip{
    \pscustom[linestyle=none]{%
      \listplot[plotstyle=curve,showpoints=false, dotstyle=triangle,linecolor=yellow]{\lbdatakappa}
      \lineto(\myk,\ymax)
      \lineto(0,\ymax)
    }
    \pscustom[linestyle=none]{%
      \listplot[plotstyle=curve,showpoints=false, dotstyle=triangle,linecolor=cyan]{\ubdatakappa}
      \lineto(\myk,0)
      \lineto(0,0)
    }
  }
  \psframe*[linecolor=black](#1)(#2)
  \endpsclip
  }
  \drawallthree(0,0)(\myk,\ymax)
  \def\zoom(#1,#2)(#3,#4)%
           (#5,#6)(#7,#8)#9{%
    {\psset{linewidth=.2pt}{%
      \psframe(#5,#6)(#7,#8)
      \psline(#5,#6)(#1,#2)
      \psline(#5,#8)(#1,#4)
      \psline(#7,#6)(#3,#2)
      \psline(#7,#8)(#3,#4)
    }%
    \rput[bl](#1,#2){%
      \psset{xunit=#9\psxunit,yunit=#9\psyunit}
      \begin{pspicture*}(#5,#6)(#7,#8)
        \psframe[fillstyle=solid,fillcolor=white](#5,#6)(#7,#8)
        \psaxes[Dx=0.5, Oy=0, Dy=0.1, linewidth=0.4pt]{->}(0,0)(0,0)(\myk.2,\ymax)
        \drawallthree(#5,#6)(#7,#8)
        \psframe(#5,#6)(#7,#8)
      \end{pspicture*}}%
    }%
  }%
  \zoom(0.1,0.17)(1.1,0.37)(1.975, 0.242)(2.025,0.252){20}
  \zoom(3.1,0.17)(4.1,0.37)(3.951,-0.001)(4.001,0.009){20}
  \psaxes[Dx=0.5, Oy=0, Dy=0.1, linewidth=0.4pt]{->}(0,0)(0,0)(\myk.2,\ymax)
  \uput[-120]{0}(\myk.2,0){$\xi$}%
\end{pspicture*}%
 
  \caption{Behavior of the estimates relative to $x =
    B^{k+\xi\newalpha}$. Colors and parameters are as in
    \ref{fig:abs-behavior}.}
  \label{fig:rel-behavior}
\end{figure}
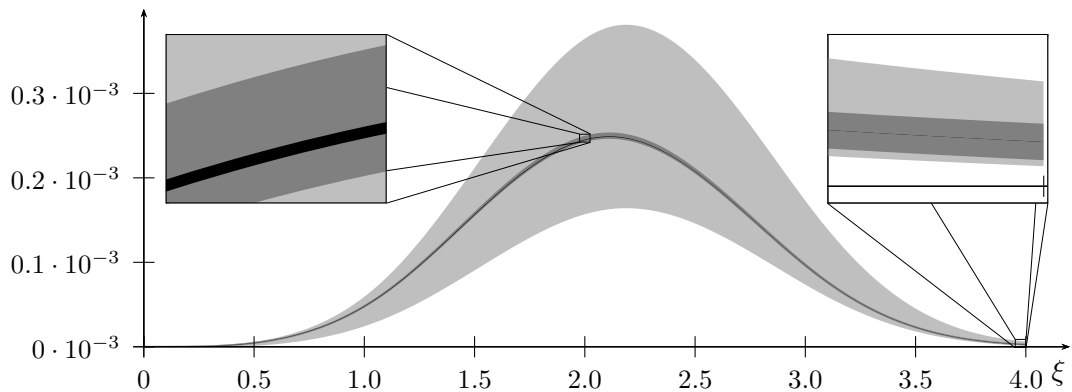

To see more we divide by $x = B^{k+\xi\newalpha}$ and therefore obtain
estimates for the ratio of $\left] B, C \right]$-grained integers $x$
in \ref{fig:rel-behavior}.
This reveals a lot about the quality of our estimates.  The black area
indicates the best that we could hope for, namely the estimate based
merely on \ref{pnt-schoenfeld}.  However, as this is difficult to
evaluate we have to approximate once more.  The $\lambda$~estimate,
shown in light gray, is clearly only of use to get a rough idea.  The
$\eta$~estimate, however, is rather close to the actual behavior and
may well serve as a basis for stochastic fine tuning of algorithms
like the general number field sieve.

\pagebreak[3] Last, \ref{fig:error-terms} illustrates the size of the
various errors terms relative to $\Qapprox[]{k}$: \newpsstyle{basis}{%
  showpoints=false,dotstyle=triangle,plotstyle=curve,%
  linecolor=black}%
\ifbw%
\newpsstyle{errorlambdatilde}{linestyle=solid,linewidth=1.5\pslinewidth}%
\newpsstyle{erroretatilde}{linestyle=solid,linewidth=\pslinewidth}%
\newpsstyle{errorboundlambdahat}{linestyle=dashed,dash=7pt 1pt}%
\newpsstyle{errorlambdahat}{linestyle=dashed,dash=6pt 2pt}%
\newpsstyle{errorkappahat}{linestyle=dashed,dash=4pt 4pt}%
\newpsstyle{errorboundlambdahat-unconditional}{linestyle=dotted,dotsep=1pt}%
\newpsstyle{errorlambdahat-unconditional}{linestyle=dotted,dotsep=1.5pt}%
\newpsstyle{errorkappahat-unconditional}{linestyle=dotted,dotsep=2.5pt}%
\newpsstyle{errornonsquarefree}{linestyle=dashed,dash=2pt
  6pt,linewidth=.75\pslinewidth}%
\else%
\newpsstyle{errorlambdatilde}{linecolor=red}%
\newpsstyle{erroretatilde}{linecolor=green}%
\newpsstyle{errorlambdahat}{linecolor=cyan}%
\newpsstyle{errorboundlambdahat}{linecolor=blue}%
\newpsstyle{errorkappahat}{linecolor=orange}%
\ifdraft
\newpsstyle{dashing}{linestyle=solid,linewidth=.1pt}
\else
\newpsstyle{dashing}{linestyle=dashed,dash=4pt 4pt}
\fi
\newpsstyle{errorlambdahat-unconditional}{style=dashing,linecolor=cyan}%
\newpsstyle{errorboundlambdahat-unconditional}{style=dashing,linecolor=blue}%
\newpsstyle{errorkappahat-unconditional}{style=dashing,linecolor=orange}%
\newpsstyle{errornonsquarefree}{linecolor=gray,linewidth=.4pt}%
\fi%
\begin{figure}[t!]
  \centering \def\pshlabel#1{\small$#1$}%
  \def\psvlabel#1{\small$#1$}%
  \psset{xunit=30mm, yunit= 400mm}
  \def\myk{4}\def\ymax{0.11}%
\setbox99\hbox{\small\psvlabel{0.04}\hskip10pt}%
\begin{pspicture*}(-\wd99,-1.2\baselineskip)(\myk.2,\ymax)%
  \readdata{\errorlambdatilde}{grainyfinalpicsdat-errorlambdatildek\myk.dat}%
  \readdata{\erroretatilde}{grainyfinalpicsdat-erroretatildek\myk.dat}%
  \readdata{\errorlambdahathat}{grainyfinalpicsdat-errorlambdahatk\myk.dat}%
  \readdata{\errorlambdahat}{grainyfinalpicsdat-errorlambdahatbetterk\myk.dat}%
  \readdata{\errorkappahat}{grainyfinalpicsdat-errorkappahatk\myk.dat}%
  \readdata{\errorlambdahathatworh}{grainyfinalpicsdat-errorlambdahatworhk\myk.dat}%
  \readdata{\errorlambdahatworh}{grainyfinalpicsdat-errorlambdahatbetterworhk\myk.dat}%
  \readdata{\errorkappahatworh}{grainyfinalpicsdat-errorkappahatworhk\myk.dat}%
  \readdata{\errornonsquarefree}{grainyfinalpicsdat-errorsquarefulk\myk.dat}%
  \readdata{\errornonsquarefreebetter}{grainyfinalpicsdat-errorsquarefulbetterk\myk.dat}%
  \readdata{\errornonsquarefreebest}{grainyfinalpicsdat-errorsquarefulbestk\myk.dat}%
  \def\drawall{%
    \psclip{%
      \pscustom[linestyle=none]{%
        \psframe*[linecolor=lightgray](0,0)(\myk,\ymax)%
      }%
    }%
    \psset{style=basis}%
    \listplot[style=errorlambdatilde]{\errorlambdatilde}%
    \listplot[style=erroretatilde]{\erroretatilde}%
    \listplot[style=errorboundlambdahat]{\errorlambdahathat}%
    \listplot[style=errorlambdahat]{\errorlambdahat}%
    \listplot[style=errorkappahat]{\errorkappahat}%
    \listplot[style=errorboundlambdahat-unconditional]{\errorlambdahathatworh}%
    \listplot[style=errorlambdahat-unconditional]{\errorlambdahatworh}%
    \listplot[style=errorkappahat-unconditional]{\errorkappahatworh}%
    \listplot[style=errornonsquarefree]{\errornonsquarefree}%
    \listplot[style=errornonsquarefree,linecolor=magenta]{\errornonsquarefreebetter}%
    \endpsclip
  }%
  \psaxes[Dx=0.5, Oy=0, Dy=0.01, linewidth=0.4pt]{->}(0,0)(0,0)(\myk.2,\ymax)
  \uput[-120]{0}(\myk.2,0){$\xi$}%
  \drawall
  \rput[tl](1.5,\ymax){%
    \footnotesize\let\small\tiny\psset{ticksize=-3pt,linewidth=.5\pslinewidth}%
    \def\ymax{2}\psset{xunit=0.25\psxunit,yunit=0.02\psyunit}%
    \begin{pspicture}(0,0)(\myk.2,\ymax)
      \psaxes[Dx=1, Oy=0, Dy=1, linewidth=0.4pt]{->}(0,0)(0,0)(\myk.2,\ymax)
      \drawall
    \end{pspicture}
  }
\end{pspicture*}%
 
  \caption{Errors of the various estimates relative to $\Qapprox{k}$.
    Parameters are as in \ref{fig:abs-behavior}.}
  \label{fig:error-terms}
\end{figure}
\begin{itemize}\itemsep0pt
  \def\legend#1{%
    \psline[style=#1](-0.5em,.3\baselineskip)(-2em,.3\baselineskip)}
\item[]\legend{errorlambdatilde} $\Kapprox{}$~error: $\left( 1-
    \frac{1}{(1+\newalpha)^{k}} \right) \Kapprox{k}(x)$.
\item[]\legend{erroretatilde} $\Eapprox{}$~error: $\left( 1 - \left(
      \frac{~ln(1+\newalpha)}{\newalpha} (1+\newalpha)^{\frac{1}{~ln(1+\newalpha)}
        - \frac{1}{\newalpha}} \right)^{k} \right) \Eapprox{k}(x)$.
\item[]\legend{errorboundlambdahat} $\Kerror{}$~error bound:
  $\cerror{k} \newalpha^{k} x$.
  \hfil\legend{errorboundlambdahat-unconditional} Unconditional
  $\Kerror{}$~error bound: $\cerror{k} \newalpha^{k} x$.
\item[]\legend{errorlambdahat} $\Kerror{}$~error: $\Kerror{k}(x)$.
  \hfil\legend{errorlambdahat-unconditional} Unconditional
  $\Kerror{}$~error: $\Kerror{k}(x)$.
\item[]\legend{errorkappahat} $\Qerror[]{}$~error: $\Qerror{k}(x)$.
  \hfil\legend{errorkappahat-unconditional} Unconditional
  $\Qerror[]{}$~error: $\Qerror{k}(x)$.
\item[]\legend{errornonsquarefree} Non-squarefree error bound: $\left(
    2^{k-1} k!  - s(k) \right) \frac{x}{B}$.
\item[]\legend{errornonsquarefree,linecolor=magenta} Non-squarefree
  error: $\left( 2^{k-1} k!  - s(k) \right) \frac{x}{B}$.
\end{itemize}
For the unconditional errors we use \citeauthor{dus98}'s unconditional
bound on $ \left| \pi(x) - ~Li(x) \right|$ which is given by
$\widehat{E}(x) = 2.3854 \frac{x}{~ln^{3} x}$ for $x \geq
\spacednumber{355991}$.

The figure shows that all the error terms but the $\Kapprox{}$~error
are sufficiently small.  Our best choice is the $\Eapprox{}$~estimate
which is ruled by the $\Eapprox{}$~error and the $\Kerror{}$~error.
Both are fairly less than 3\% of the target value at least in the
middle of the interval.  The estimations are more difficult close to
the boundaries.  It is also positive that, at the parameters of our
interest, most error terms are still comparative in size to the
contributions of the $\Qerror[]{}$~error, which is induced by the
prime number theorem.  Definitely, the $\eta$~estimate, combining the
$\Eapprox{}$~error and the $\Kerror{}$~error, is good enough for
practical purposes, as for example the fine tuning of the general
number field sieve.

\begin{acknowledge}
  This work was funded by \ifanonymous[anonymized].\else the b-it
  foundation and the state of North~Rhine-Westphalia.
  We thank Igor Shparlinski for his remarks.
  \fi
\end{acknowledge}

\bibliography{journals,refs,lncs}

\end{document}